\documentclass[11pt,reqno]{amsart}

\usepackage{amsmath,amssymb,amsthm,mathtools}
\usepackage[margin=2.8cm]{geometry}
\usepackage{enumitem}
\usepackage[hidelinks]{hyperref}
\hypersetup{
  pdftitle={Gopakumar--Vafa Invariant and Macdonald Formula},
  pdfauthor={Lutian Zhao}
}

\newtheorem{theorem}{Theorem}[section]
\newtheorem{proposition}[theorem]{Proposition}
\newtheorem{lemma}[theorem]{Lemma}
\newtheorem{corollary}[theorem]{Corollary}
\newtheorem{conjecture}[theorem]{Conjecture}
\theoremstyle{definition}
\newtheorem{definition}[theorem]{Definition}
\newtheorem{remark}[theorem]{Remark}

\newcommand{\C}{\mathbb C}
\newcommand{\Q}{\mathbb Q}
\newcommand{\Pp}{\mathbb P}
\newcommand{\Ocal}{\mathcal O}
\newcommand{\IC}{\operatorname{IC}}
\newcommand{\Hilb}{\operatorname{Hilb}}
\newcommand{\Supp}{\operatorname{Supp}}
\newcommand{\Fitt}{\operatorname{Fitt}}
\newcommand{\Hom}{\operatorname{Hom}}
\newcommand{\codim}{\operatorname{codim}}
\newcommand{\Tot}{\operatorname{Tot}}
\newcommand{\PT}{\operatorname{PT}}

\newcommand{\red}{\mathrm{red}}
\newcommand{\sm}{\mathrm{sm}}
\newcommand{\nr}{\mathrm{nr}}
\newcommand{\Eff}{\operatorname{Eff}}
\newcommand{\Mac}{\operatorname{Mac}}
\newcommand{\ExpChow}{\operatorname{Exp}_{\mathrm{Chow}}}
\newcommand{\bigboxtimes}{\mathop{\boxtimes}\limits}

\title[Gopakumar--Vafa Invariant and Macdonald Formula]
{Gopakumar--Vafa Invariant and Macdonald Formula}

\author{Lutian Zhao}
\address{Kavli Institute for the Physics and Mathematics of the Universe,
The University of Tokyo, 5-1-5 Kashiwanoha, Kashiwa, Chiba 277-8583,
Japan}
\email{lutian.zhao@ipmu.jp}

\begin{document}

\begin{abstract}
We use a derived constructible Chow exponential and prove that its
coefficients are perverse minimal extensions from reduced cycles.  With
compatible orientations, we formulate the cohomological PT/GV relation by
this exponential.  If the relation holds over reduced cycles, its extension
to the full Chow varieties is equivalent to the absence of nonreduced strict
supports, and it implies the refined and numerical formulas.  For local
$\Pp^2$ and $0\le n\le d+1$, we identify the stable-pair space with the
smooth relative Hilbert scheme and the vanishing-cycle sheaf with its
intersection complex.  We determine the first reducible summand and give
del Pezzo examples.  In degree $(2,4)$, we determine the incidence--ribbon
union, its canonical PT critical germ, and the attachment triangle, and we
formulate the primitive degree-two KKV comparison.
\end{abstract}

\maketitle

\section{Introduction}
Gopakumar and Vafa expressed the all-genus topological-string series of a
Calabi--Yau threefold in terms of integral multiplicities of spinning BPS
states \cite{GopakumarVafaI,GopakumarVafaII}.  Hosono--Saito--Takahashi
proposed extracting the two Lefschetz actions from the intersection
cohomology of moduli spaces of one-dimensional sheaves
\cite{HosonoSaitoTakahashi}.  Pandharipande--Thomas introduced stable pairs
\cite{PT} and proved the strong BPS rationality conjecture for irreducible
curve classes \cite{PTBPS}.  For local $\Pp^2$, Katz--Klemm--Vafa computed
the BPS multiplicities from families of plane curves, Macdonald's formula,
and corrections carried by special curve strata \cite{KKV}.
Choi--Katz--Klemm subsequently compared stable-pair spaces with relative
Hilbert schemes in the low-cokernel range, proposed the two-spin refined
stable-pair/BPS product, and checked it in low degree for local $\Pp^2$ and
local $\Pp^1\times\Pp^1$ \cite{CKK}.

The sheaf-theoretic approach began with the relative-Lefschetz construction
above.  Kiem--Li constructed vanishing-cycle perverse sheaves on fine moduli
spaces of simple sheaves on smooth projective Calabi--Yau threefolds; their
construction uses a square root of the determinant line and, in general, is
made after an \'etale Galois cover \cite{KiemLi}.  The gluing theorem
for oriented $d$-critical loci gives a global perverse sheaf once the required
orientation has been chosen \cite{BBDJS}.

Maulik--Toda
defined Gopakumar--Vafa invariants from such vanishing cycles, conjectured
their agreement with Gromov--Witten and stable-pair theories, and proved the
numerical PT/GV formula for irreducible one-cycles on local surfaces
\cite{MaulikToda}.  Toda extended this construction to BPS sheaves on coarse
moduli spaces of semistable one-dimensional sheaves, proved stability
independence, and established Euler-characteristic independence for
primitive cycles on local surfaces \cite{TodaGVWallCrossing}.  He later
described stable D0--D2--D6 objects on a local surface as a dual obstruction
cone over the surface-pair stack \cite{TodaLocal}.  Kinjo's dimensional
reduction theorem computes the vanishing-cycle theory of the associated
shifted cotangent space from the base \cite{Kinjo}.  Dimensional reduction
does not by itself control strict supports after pushforward to the Chow
variety.

Several recent results clarify neighboring forms of the correspondence.
P\u{a}durariu--Toda proved a categorical DT/PT decomposition for reduced
curve classes on local surfaces \cite{PadurariuToda}.  Pardon proved the
numerical GW/DT correspondence for Calabi--Yau threefolds by reduction to
local curves \cite{Pardon}.  Davison--Koseki determined the degree-two GV
invariants of generic local curves in genus two and proved the corresponding
GV/GW formula \cite{DavisonKoseki}.

After embedding an integral locally planar curve in a sufficiently versal
family, Maulik--Yun and, independently, Migliorini--Shende lifted
Macdonald's formula to a perverse-sheaf identity for relative Hilbert
schemes and compactified Jacobians \cite{MaulikYun,MiglioriniShende}.
Their relative statements use the smoothness of the relevant relative
Hilbert schemes and codimension bounds for the cogenus strata.
Migliorini--Shende--Viviani
extended the support decomposition to reduced reducible curves; the new
terms are formed from connected subcurves and Hilbert schemes of partial
normalizations, and their relative theorem applies to versal or
independently broken H-smooth families \cite{MSV}.

On the sheaf-moduli side, Maulik--Shen proved cohomological
$\chi$-independence, including the perverse and Hodge filtrations, for ample
classes on toric del Pezzo surfaces \cite{MaulikShen}.  Yuan's estimates over
multiple curves extended this result to every del Pezzo surface
\cite{YuanNonreduced}.  Kinjo--Koseki proved $\chi$-independence for
generalized Gopakumar--Vafa invariants of arbitrary local curves and for BPS
cohomology of Higgs bundles \cite{KinjoKoseki}.  Very recently,
Davison--Hennecart--Kinjo--Schiffmann--Vasserot proved Toda's
$\chi$-independence conjecture relative to the Chow variety for
one-dimensional sheaves on quasi-projective symplectic surfaces
\cite{DavisonHennecartKinjoSchiffmannVasserot}.

For compactified Jacobian fibrations of integral locally planar curves,
Maulik--Shen--Yin developed a Fourier-transform theory, proved
multiplicativity of the perverse filtration, and obtained the
$P\supset C$ direction of the $P=C$ conjecture for local $\Pp^2$
\cite{MaulikShenYin}.  Kononov--Pi--Shen formulated $P=C$ and proved it for
degrees $d\le4$ \cite{KononovPiShen}.  Kononov--Lim--Moreira--Pi then
determined the relevant cohomology rings through degree five and verified
the refined GV/PT correspondence for local $\Pp^2$ in degree five
\cite{KononovLimMoreiraPi}.  Pi--Shen--Si--Zhang proved cohomological
stabilization and an asymptotic refined-BPS product formula for local del
Pezzo surfaces \cite{PiShenSiZhang}; Guo--Wu computed the Poincar\'e
polynomials for the plane moduli spaces through degree sixteen
\cite{GuoWu}.

The geometry over nonreduced curves has also developed.  Luan classified
the components and multiplicities of Hilbert schemes of points on
nonreduced plane curves \cite{LuanNonreduced} and later studied punctual
Hilbert schemes of the nonreduced nodal germs $x^uy^v=0$
\cite{LuanNodal}.  Savarese--Viviani described the strata, components, and
tangent spaces of moduli of torsion-free sheaves on ribbons
\cite{SavareseViviani}.  Wu proved a full-support theorem for universal
compactified Jacobians over the moduli of stable curves
\cite{WuUniversalJacobian}.  These results sharpen the sheaf,
compactified-Jacobian, and Hilbert-scheme sides, but they do not identify the
stable-pair vanishing-cycle pushforward over the nonreduced strata of a
complete Chow variety.

The complete family of plane curves is not H-smooth: already the relative
four-point Hilbert scheme of conics is singular, as shown in
Proposition~\ref{prop:conic-correction-node}.  Passing to a finite \`etale
cover labels components and handles their permutations, but does not restore
this smoothness.  We therefore do not apply the global MSV theorem to the
complete plane family.  Instead, we apply its versal calculation at the
generic partition strata and use the node-smoothing transversality of
Lemma~\ref{lem:transverse-nodal-base-change} to pull the resulting summands
to the plane linear system.

There is a second distinction
over the nonreduced locus: the constant sheaf of a singular incidence
variety, its intersection complex, and the vanishing-cycle sheaf on the
stable-pair space need not have the same direct image.  The ordinary Euler
characteristic of the underlying incidence space does not determine its IC
or vanishing-cycle cohomology.

The present paper organizes the Macdonald term and the reducible terms by
addition of effective cycles.  Nekrasov and Okounkov introduced a symmetric
algebra over the Chow variety by pushing external powers along cycle-addition
maps and taking symmetric-group invariants
\cite[Section~2.3.5]{NekrasovOkounkov}; see also
\cite[Section~5.3.3]{OkounkovTakagi}.  We use its derived constructible form,
which retains cohomological shifts, Koszul signs, and the permutation local
systems of equal factors.  Write
$B_d=\operatorname{Chow}_{dH}(\Pp^2)$.  Multiplication of homogeneous
polynomials gives finite addition maps $B_a\times B_b\to B_{a+b}$, so
$\coprod_{d\ge0}B_d$ is a graded Chow monoid.  Constructible complexes on
its components carry a convolution product, and $\ExpChow$ denotes the
completed derived symmetric algebra for this product.  The same construction
applies to an additive semigroup of effective curve classes satisfying the
connectedness and density hypotheses of Section~\ref{sec:hilbert}.

Let $\mathcal A_d(q)$ be the connected Macdonald complex in degree $d$, and
let $\mathcal Z_{\PT}(q,Q)$ be the normalized semisimplified PT
direct-image series: its $(\chi,\beta)$-coefficient is
$q^\chi(R\Pi_{\chi,\beta,*}\Phi_{\chi,\beta}^{\PT})^{\mathrm{ss}}
[-\chi]Q^\beta$.  The central conjecture is the coefficientwise identity
\begin{equation}\label{eq:intro-exponential}
  \mathcal Z_{\PT}(q,Q)
  \cong
  \ExpChow\!\left(\bigoplus_{d\ge1}\mathcal A_d(q)Q^d\right).
\end{equation}
The degree-$d$ coefficient on the right is a direct sum indexed by
partitions of $d$.  Repeated parts are governed by derived symmetric powers,
not by division of an ordered product by the order of a symmetric group.
All coefficients and isomorphisms in this statement are taken in the
completed semisimple graded constructible category.  Equation
\eqref{eq:intro-exponential} is not deduced in this paper from the support
conditions below.

The left side of \eqref{eq:intro-exponential} is defined after choosing
compatible orientations for which the stable-pair vanishing-cycle sheaves
are globally defined.  Put $B_\beta=|\Ocal_S(\beta)|$ and
$g_\beta=p_a(\beta)$.  Write
$j_\beta:B_\beta^{\red}\hookrightarrow B_\beta$ for the reduced-cycle open
immersion and put $D_\beta^{\nr}=B_\beta\setminus B_\beta^{\red}$.  For
$\chi=1-g_\beta+n$, the condition
$(\mathrm{Hilb})_{\beta,n}$ means that the stable-pair vanishing-cycle sheaf
over reduced curves agrees with the relative-Hilbert intersection complex:
\begin{equation}\label{eq:intro-hilb-condition}
 \alpha_{\beta,n}^*
 \bigl(\Phi_{\beta,n}^{\PT}|_{P_{\beta,n}^{\red}}\bigr)
 \cong\IC_{\mathcal C_{\beta,\red}^{[n]}}.
\end{equation}
Here $P_{\beta,n}^{\red}$ is the inverse image of $B_\beta^{\red}$,
$\mathcal C_{\beta,\red}^{[n]}$ is the relative Hilbert scheme over that
open set, and $\alpha_{\beta,n}$ is the surface stable-pair isomorphism.  This notation
keeps the condition distinct from the intersection-complex symbol $\IC$.

Put $K_{\chi,\beta}=R\Pi_{\chi,\beta,*}\Phi_{\chi,\beta}^{\PT}$.  The
condition $(\mathrm{NR})_{\chi,\beta}$ means that no simple constituent of
any ${}^pH^iK_{\chi,\beta}$ has support contained in $D_\beta^{\nr}$.
Equivalently,
\begin{equation}\label{eq:intro-nr-condition}
 K_{\chi,\beta}^{\mathrm{ss}}
 \cong{}^pj_{\beta,!*}
       (j_\beta^*K_{\chi,\beta})^{\mathrm{ss}}.
\end{equation}
Thus $(\mathrm{NR})$ is a statement about strict supports after proper
pushforward.  It does not assert that the PT sheaf has zero stalks over
nonreduced curves.  The precise definitions are
Definitions~\ref{def:reduced-hilbert-comparison}
and~\ref{def:nonreduced-support}.

For local $\Pp^2$, put $g_d=(d-1)(d-2)/2$ and abbreviate
$(\mathrm{Hilb})_{d,n}=(\mathrm{Hilb})_{dH,n}$ and
$(\mathrm{NR})_{d,n}=(\mathrm{NR})_{1-g_d+n,dH}$.  Thus
$(\mathrm{NR})_{2,4}$ concerns the direct image from
$P_5(\Tot(K_{\Pp^2}),2H)$ to $B_2$.

\begin{theorem}[Nonreduced support and PT/GV]
\label{thm:intro-chow}
Let $\Gamma$ be an effective-class semigroup as in
Section~\ref{sec:chow-exponential}, and suppose that every factorization source has
a dense open subset whose cycle sum is reduced.

\begin{enumerate}[label=\textnormal{(\roman*)},leftmargin=9mm]
\item Every coefficient $\mathcal R_{\chi,\beta}$ of the Chow exponential is
semisimple and is recovered from the reduced Chow locus:
\[
 \mathcal R_{\chi,\beta}
 \cong{}^pj_{\beta,!*}j_\beta^*\mathcal R_{\chi,\beta}.
\]

\item Suppose compatible oriented PT vanishing-cycle sheaves have been
chosen and \eqref{eq:intro-exponential} is known after restriction to every
$B_\beta^{\red}$.  Then, in the completed semisimple graded category, the
identity holds on the full Chow varieties precisely when
$(\mathrm{NR})_{\chi,\beta}$ holds for every coefficient.

\item Suppose \eqref{eq:intro-exponential} holds.  The polarization of the
Jacobian local system and an ample class on each projective Chow variety
give two commuting Lefschetz actions.  Their gradings give the refined PT/GV product
\eqref{eq:ckk-product}.  Applying hypercohomological Euler characteristic
gives the numerical PT/GV correspondence.
\end{enumerate}
\end{theorem}

\begin{proof}
The three statements are Proposition~\ref{prop:exponential-minimal-extension},
Corollary~\ref{cor:reduced-exponential}, and
Proposition~\ref{prop:refined-realization}, respectively.
\end{proof}

Consequently, once the PT/GV identity is established on the reduced Chow
loci, it extends to every Chow variety in the completed semisimple graded
category precisely when all conditions $(\mathrm{NR})_{\chi,\beta}$ hold.
Its Euler characteristic gives the numerical PT/GV formula, and its two
Lefschetz gradings give the refined formula.

For local $\Pp^2$, this also removes the bound in the KKV recursion.  If the
reduced identity and all conditions $(\mathrm{NR})_{d,n}$ hold, the KKV
formula of Corollary~\ref{cor:kkv-all-degree-nr} holds for every $n$.  The
bound $n\le d+1$ is needed only to replace the PT vanishing-cycle coefficient
by the elementary Euler characteristic of the smooth incidence space.

The reduced identity in part \textnormal{(ii)} contains more information
than $(\mathrm{Hilb})_{\beta,n}$.  The reduced Hilbert comparison replaces
the PT sheaf by
the relative-Hilbert intersection complex; the local
Migliorini--Shende--Viviani calculation then determines its generic
reducible summands.  These two steps are kept separate throughout the
paper.

For the complete plane linear system, the sheaf conditions and the first
reducible coefficient can be proved directly.

\begin{theorem}[Plane stable-pair coefficients]
\label{thm:intro-plane}
Let $X=\Tot(K_{\Pp^2})$ and $g_d=(d-1)(d-2)/2$.  Let
$\mathcal C_d^{[n]}\to B_d$ be the relative Hilbert scheme,
$\Pi_{d,n}$ the PT support map, and $\Phi_{d,n}^{\PT}$ the vanishing-cycle
sheaf for the coefficientwise orientation chosen in part~\textnormal{(i)}.

\begin{enumerate}[label=\textnormal{(\roman*)},leftmargin=9mm]
\item For $d\ge1$ and $0\le n\le d+1$,
\[
 P_{1-g_d+n}(X,d)\cong\mathcal C_d^{[n]},
 \qquad
 \Phi_{d,n}^{\PT}\cong\IC_{\mathcal C_d^{[n]}}.
\]
The relative Hilbert scheme is smooth and the second isomorphism uses its
canonical orientation.  Thus $(\mathrm{Hilb})_{d,n}$ holds.

\item For the same orientation, if $d\ge3$ and $0\le n\le d+1$, then
$(\mathrm{NR})_{d,n}$ holds.

\item Put $K_{d,n}=R\Pi_{d,n,*}\Phi_{d,n}^{\PT}$ and let
$\mathcal S_{d,n}$ be its full-support Macdonald summand.  For $d\ge3$,
\[
 K_{d,n}\cong\mathcal S_{d,n}\quad(0\le n\le d-2),
 \qquad
 K_{d,d-1}\cong
 \mathcal S_{d,d-1}\oplus\IC_{Z_{d-1,1}}.
\]
Here $Z_{d-1,1}$ is the closure of the locus of a degree-$(d-1)$ curve
together with a line.
\end{enumerate}
\end{theorem}

\begin{proof}
Parts \textnormal{(i)} and \textnormal{(ii)} are
Theorem~\ref{thm:plane-smooth-range} and
Corollary~\ref{cor:plane-nonreduced}.  Part \textnormal{(iii)} is
Corollary~\ref{cor:before-first-boundary} together with
Theorem~\ref{thm:plane-first-window}.
\end{proof}

Taking Euler characteristics gives the unitriangular Macdonald recursion of
Katz--Klemm--Vafa and its reducible-curve correction.  In particular, the
quartic row and the quintic invariants in genera two through six follow from
the proved support range.  The genus-one and genus-zero coefficients are
expressed in terms of the residual complexes; they take the displayed KKV
values when the corresponding coefficients of
\eqref{eq:intro-exponential} hold.

The zero-section argument is not restricted to the plane.

\begin{theorem}[Del Pezzo incidence calculations]
\label{thm:intro-del-pezzo}
Let $S$ be a del Pezzo surface and let $\tau_S(\beta)$ be the threshold in
\eqref{eq:del-pezzo-threshold}.
\begin{enumerate}[label=\textnormal{(\roman*)},leftmargin=9mm]
\item Suppose $n<\tau_S(\beta)$ and every length-$n$ subscheme of $S$
imposes independent conditions on $|\beta|$.  Then
$P_{1-g_\beta+n}(\Tot(K_S),\beta)$ is the smooth incidence projective
bundle over $\Hilb^n(S)$, and its coefficientwise vanishing-cycle sheaf is
the intersection complex.

\item For local $\Pp^1\times\Pp^1$ in class $(2,3)$, the full-support
Macdonald coefficients are
\[
 (n_{(2,3)}^0,n_{(2,3)}^1,n_{(2,3)}^2)=(-110,68,-12).
\]

\item For $S_r=\operatorname{Bl}_{p_1,\ldots,p_r}\Pp^2$, with general
points and $5\le r\le8$, the aggregated full-support Macdonald
coefficients in anticanonical degree $9-r$ are
\[
\begin{array}{c|rrrr}
S_r&E_5&E_6&E_7&E_8\\ \hline
N_{9-r}^1&5&-4&3&-2\\
N_{9-r}^0&-192&243&-272&252.
\end{array}
\]
\end{enumerate}
\end{theorem}

\begin{proof}
Part \textnormal{(i)} is Theorem~\ref{thm:del-pezzo-incidence-range}.
Parts \textnormal{(ii)} and \textnormal{(iii)} are proved in
Section~\ref{sec:examples}; the former includes the complete codimension-two
reducible correction, and the latter follows from the anticanonical
incidence spaces and the enumeration of rational classes.
\end{proof}

The first calculation beyond the smooth-incidence range is degree two with
four points.  It separates the shifted constant sheaf used by an incidence
Euler-characteristic calculation, the intersection complex of the incidence
space, and the PT vanishing-cycle sheaf.

\begin{theorem}[Degree two with four points]
\label{thm:intro-degree-two}
Let
\[
 M=P_5(\Tot(K_{\Pp^2}),2),\qquad Y=\mathcal C_2^{[4]},
\]
and let $N_2\subset B_2$ be the locus of double lines.
Then $M=Y\cup_{\mathfrak B}R$ scheme-theoretically, where $R$ is the smooth
compactification of the Ferrand-ribbon locus.  Along
$\mathfrak S\cong\mathfrak B$, the incidence component has a transverse
threefold-node singularity, and at every point of $\mathfrak B$ the
completed classical germ of $M$ is
\[
 \mathbb A^6\times
 \operatorname{Crit}\bigl(c(x_1x_2+x_3x_4)\bigr)
\]
as a formal scheme.  In $\operatorname{Perv}(Y)$ there is a short exact sequence
\begin{equation}\label{eq:intro-degree-two-correction}
 0\longrightarrow\Q_{\mathfrak S}[6]
 \longrightarrow\Q_Y[9]
 \longrightarrow\IC_Y\longrightarrow0.
\end{equation}

The canonical PT $d$-critical germ along $\mathfrak B$ is represented by
$cq$, up to the smooth factor and quadratic stabilization.  Write $\Phi_M$
for its vanishing-cycle sheaf with the orientation of
Proposition~\ref{prop:degree-two-attachment}.  The support maps
$\Pi:M\to B_2$ and $\pi:Y\to B_2$ fit into a triangle
\[
 R\pi_*\Q_Y[9]\longrightarrow R\Pi_*\Phi_M
 \longrightarrow\Q_{N_2}[7]
 \xrightarrow{\delta}(R\pi_*\Q_Y[9])[1].
\]
The $N_2$-supported components of the attaching morphism $\delta$, together
with the perverse decomposition of $R\pi_*\Q_Y[9]$, determine
$(\mathrm{NR})_{2,4}$.
\end{theorem}

\begin{proof}
The scheme structure and completed local germ are
Corollary~\ref{cor:degree-two-union} and
Theorem~\ref{thm:degree-two-classical-germ}.  The short exact sequence is
Proposition~\ref{prop:degree-two-vc}.  The PT $d$-critical statement and
the triangle are
Propositions~\ref{prop:degree-two-pt-critical-model}
and~\ref{prop:degree-two-attachment}.
\end{proof}

Sequence \eqref{eq:intro-degree-two-correction} is the degree-two incidence
correction in perverse-sheaf form.  Since
$\mathfrak S\to N_2$ is a $\Pp^4$-bundle, its proper direct image contains
the Lefschetz string
\[
 \bigoplus_{i=0}^4\IC_{N_2}[4-2i].
\]
The extreme term $\IC_{N_2}[4]$ is the primitive nonreduced contribution
relevant to the Katz--Klemm--Vafa inversion.  The PT attachment triangle
supplies the second primitive map.  Conjecture~\ref{conj:degree-two-kkv}
asserts the decomposition of $A^{\mathrm{ss}}$, the nonvanishing of both
primitive maps, and the existence of the comparison morphism $\kappa_2$
with the stated kernel.  Together these remove the primitive term rather
than making the PT sheaf vanish over double lines.

More precisely, writing $H=R\pi_*\IC_Y$ and $T=R\Pi_*\Phi_M$, the proposed
degree-two KKV correction is the short exact sequence
\[
 0\longrightarrow\IC_{N_2}\longrightarrow{}^pH^{-3}(H^{\mathrm{ss}})
 \xrightarrow{\ \kappa_2\ }{}^pH^{-3}(T^{\mathrm{ss}})
 \longrightarrow0.
\]
Neither the local critical chart nor the Euler characteristics determine
the existence of $\kappa_2$ or the required primitive-map nonvanishing.
The final section formulates the analogous comparison at $2L+D$ as a
conjecture for every degree, using the finite addition of a double line and
a residual degree-$(d-2)$ curve.

Sections~\ref{sec:pt}--\ref{sec:chow-exponential} establish the
stable-pair/relative-Hilbert comparison, the local MSV support terms, and
the Chow exponential.  Sections~\ref{sec:plane}--\ref{sec:examples}
derive the plane support range, the KKV recursion, and the del Pezzo
calculations.  Section~\ref{sec:degree-two-correction} studies degree
$(2,4)$, including the incidence and ribbon components, the canonical PT
critical germ, the attachment triangle, and the higher-degree double-line
conjecture.

\section{Stable pairs over the reduced Chow base}\label{sec:pt}

For a vector bundle $E$, we use the convention that $\Pp(E)$ parametrizes
one-dimensional subspaces of its fibres.

Let $S$ be a smooth del Pezzo surface, let $p:X=\Tot(K_S)\to S$ be the
projection, and fix a nonzero effective class $\beta$.  Put

\[
  B_\beta=|\Ocal_S(\beta)|,
  \qquad
  g_\beta=p_a(\beta)
  =1+\frac{\beta\cdot(\beta+K_S)}2.
\]

For $n\ge0$, write $P_{\beta,n}=P_{1-g_\beta+n}(X,\beta)$.

A point of $P_{\beta,n}$ is a stable pair in the sense of \cite{PT},
$I^\bullet=[\Ocal_X\to F]$, where $F$ is pure of dimension one, has
compact support and class $\beta$, and satisfies
$\chi(F)=1-g_\beta+n$.  The Fitting support of $p_*F$ defines
$\Pi_{\beta,n}:P_{\beta,n}\to B_\beta$.

This morphism is projective.  The fibre coordinate restricts to a section
of negative degree on every reduced component of a compact curve, so every
support is set-theoretically contained in the zero section.  Fix a projective
compactification of $X$.  The projective stable-pair construction
\cite[Section~1]{PT} makes the family with the given Hilbert polynomial
bounded.  Since every member is set-theoretically supported on $S$,
Noetherianity gives a uniform power of the ideal of $S$ which annihilates
every pair sheaf.  Hence $P_{\beta,n}$ is a
closed stable-pair moduli space on a fixed projective infinitesimal
thickening of $S$, and its Fitting-support morphism is proper.

Let $j_\beta:B_\beta^{\red}\hookrightarrow B_\beta$ be the open locus of
reduced divisors, and set
$P_{\beta,n}^{\red}=P_{\beta,n}\times_{B_\beta}B_\beta^{\red}$.

\subsection{Stable pairs on the zero section}

The affine morphism $p$ identifies a compactly supported coherent sheaf
$F$ on $X$ with a coherent sheaf $G=p_*F$ on $S$, together with an
action of $K_S^{-1}$.  By adjunction, this action is a Higgs field
$\phi_F:G\to G\otimes K_S$.
Since $p|_{\Supp F}$ is affine and proper, it is finite; purity of $F$
therefore implies purity of $G$.

The sheaf $F$ is scheme-theoretically supported on the zero section if and
only if $\phi_F=0$.

\begin{proposition}\label{prop:reduced-zero-section}
If $\Fitt_0(G)$ is a reduced divisor, then $\phi_F=0$.  Hence every
stable pair in $P_{\beta,n}^{\red}$ is scheme-theoretically supported on
the zero section.
\end{proposition}

\begin{proof}
Let $C=\Fitt_0(G)=C_1\cup\cdots\cup C_r$.  At the generic point of
$C_i$, choose a transverse parameter $t$ in the regular local ring of $S$.
The order of the zeroth Fitting ideal along $C_i$ is the length of $G$ over
the discrete valuation ring generated by $t$.  Since $C$ is reduced, this
order is one.  It follows that $tG$ vanishes at the generic point and that
$G$ has rank one over $\Ocal_{C_i}$.  Consequently $I_CG$ is supported at
finitely many closed points.  Purity of $G$ excludes zero-dimensional
subsheaves, hence $I_CG=0$.  Thus $G$ is an $\Ocal_C$-module which is
torsion-free of rank one on every component.

Let $L_i$ be the torsion-free quotient of the pullback of $G$ to the
normalization $\nu_i:\widetilde C_i\to C_i$.  It is a line bundle on the
smooth curve $\widetilde C_i$.  Pulling back the Higgs field and passing to
torsion-free quotients gives a homomorphism
$L_i\to L_i\otimes\nu_i^*K_S$.  After tensoring by $L_i^{-1}$, this is a
section of $\nu_i^*K_S$.  Its degree is $K_S\cdot C_i<0$, because $S$ is del
Pezzo, so the section is zero.  The original Higgs field is therefore zero
at the generic point of every component of $C$.  Its image is
zero-dimensional, while $G\otimes K_S$ is pure of dimension one.  Hence the
image vanishes and $\phi_F=0$.
\end{proof}

Let $\mathcal C_\beta\subset S\times B_\beta$ be the universal divisor.  Set
$\mathcal C_\beta^{[n]}=\Hilb^n(\mathcal C_\beta/B_\beta)$ and denote its
projection by $\pi_\beta^{[n]}:\mathcal C_\beta^{[n]}\to B_\beta$.  Write
$\mathcal C_{\beta,\red}^{[n]}$ for its restriction to
$B_\beta^{\red}$.

\begin{proposition}\label{prop:hilbert-pt}
There is a canonical isomorphism over $B_\beta^{\red}$
\[
  \alpha_{\beta,n}:\mathcal C_{\beta,\red}^{[n]}
  \xrightarrow{\sim}P_{\beta,n}^{\red}.
\]
On closed points it sends $(C,Z)$ to
$[\Ocal_X\to\iota_{C,*}(I_{Z/C}^{\vee})]$.
\end{proposition}

\begin{proof}
The standard surface stable-pair correspondence identifies the stable-pair
space on $S$ with the relative Hilbert scheme
\cite[Proposition~1]{CKK}.  We recall the construction.
Since $C$ is a Cartier divisor on a smooth surface, it is Gorenstein.  A
length-$n$ subscheme $Z\subset C$ gives the rank-one torsion-free sheaf
$I_{Z/C}^{\vee}$ and an exact sequence
\[
  0\longrightarrow\Ocal_C\longrightarrow I_{Z/C}^{\vee}
  \longrightarrow Q\longrightarrow0,
  \qquad \operatorname{length}(Q)=n.
\]
After pushforward to $X$, this is a stable pair of Euler characteristic
$1-g_\beta+n$.

Conversely, Proposition~\ref{prop:reduced-zero-section} writes every pair
over $B_\beta^{\red}$ as
$0\to\Ocal_C\to G\to Q\to0$.  Duality on the Gorenstein curve gives
\[
  0\longrightarrow G^\vee\longrightarrow\Ocal_C
  \longrightarrow\mathcal Ext_C^1(Q,\Ocal_C)\longrightarrow0.
\]
Indeed, $\mathcal{H}\!om_C(Q,\Ocal_C)=0$ because $\Ocal_C$ is pure, and
$\mathcal Ext_C^1(G,\Ocal_C)=0$ because a torsion-free sheaf on a Gorenstein
curve is maximal Cohen--Macaulay.  Local duality identifies
$\mathcal Ext_C^1(Q,\Ocal_C)$ with the dual of the zero-dimensional sheaf
$Q$, so it has the same length $n$.  The kernel $G^\vee\subset\Ocal_C$ is
therefore the ideal $I_{Z/C}$ of a unique length-$n$ subscheme $Z\subset C$.
Gorenstein duality also gives $G\cong G^{\vee\vee}$, and hence
$G\cong I_{Z/C}^\vee$.

For a family of Gorenstein curves, the family statement in
\cite[Proposition~1]{CKK} proves that these duality constructions commute with
base change and give inverse transformations of the two moduli functors.
Thus the closed-point correspondence above is induced by a canonical
isomorphism of the relative Hilbert scheme with the stable-pair moduli space.
\end{proof}

The range in which every stable pair remains on the zero section admits a
useful numerical form for any del Pezzo surface.  If \(\beta\) has no
decomposition into two nonzero effective classes, set
\(\tau_S(\beta)=\infty\).  Otherwise define
\begin{equation}\label{eq:del-pezzo-threshold}
 \tau_S(\beta)=
 \min_{\substack{\ell\ge2,\ \gamma_i>0\\
                  \gamma_0+\cdots+\gamma_{\ell-1}=\beta}}
 \left\{
   \sum_{i<j}\gamma_i\cdot\gamma_j
   +\sum_{i=1}^{\ell-1}i(-K_S\cdot\gamma_i)
 \right\}.
\end{equation}
The order of the classes in this minimum matters.

\begin{lemma}[The del Pezzo zero-section threshold]
\label{lem:del-pezzo-zero-section-threshold}
Let \(C\subset\Tot(K_S)\) be a compact Cohen--Macaulay curve of class
\(\beta\) which is not scheme-theoretically contained in the zero section.
Then
\[
 \chi(\Ocal_C)\ge 1-g_\beta+\tau_S(\beta).
\]
Consequently every stable pair in \(P_{\beta,n}\) is supported on the zero
section when \(n<\tau_S(\beta)\).
\end{lemma}

\begin{proof}
Let \(I\) be the ideal of the zero section and \(J\) the ideal of \(C\).
The reduced support of a compact curve lies in the zero section, so
\(I^{r+1}\subset J\) for some \(r\).  Put
\[
 G_k=(J+I^k)/(J+I^{k+1}),\qquad 0\le k\le r.
\]
Since \(I^k/I^{k+1}\cong K_S^{-k}\), the sheaf \(G_k\) is a quotient of
\(K_S^{-k}\).  The reflexive hull of its kernel is
\(K_S^{-k}(-D_k)\) for an effective divisor \(D_k\).  If
\(\gamma_k=[D_k]\), comparison with the reflexive hull gives an exact
sequence
\[
 0\longrightarrow T_k\longrightarrow G_k
 \longrightarrow\Ocal_{D_k}(-kK_S)\longrightarrow0,
\]
where \(T_k\) has finite length.  Hence
\[
 \chi(G_k)\ge
 \chi\bigl(\Ocal_{D_k}(-kK_S)\bigr)
 =1-g_{\gamma_k}+k(-K_S\cdot\gamma_k)
\]
whenever \(D_k\ne0\).  Layers with \(D_k=0\) have finite length and may
be omitted from the lower bound.

If $C$ is not contained in the zero section, choose the largest $k>0$ for
which $I^k\Ocal_C\ne0$.  Then $I^k\Ocal_C=G_k$.  If this sheaf had finite
length, it would be a nonzero zero-dimensional subsheaf of the
Cohen--Macaulay sheaf $\Ocal_C$, contrary to purity.  Thus $D_k\ne0$ for
some $k>0$.

The positive-dimensional layers have total class \(\beta\).  Write their
indices as \(0=k_0<k_1<\cdots<k_{\ell-1}\) and their classes as
\(\gamma_0,\ldots,\gamma_{\ell-1}\).  The first index is zero because
\(G_0\) contains the reduced support, and \(k_i\ge i\).  Adjunction gives
\[
 \sum_i(1-g_{\gamma_i})
 =1-g_\beta+\sum_{i<j}\gamma_i\cdot\gamma_j.
\]
As \(-K_S\) is ample, replacing \(k_i\) by \(i\) can only decrease the
right-hand side of the preceding lower bound.  A curve not contained in
the zero section has at least one positive layer with \(k_i>0\), so
\(\ell\ge2\).  Formula \eqref{eq:del-pezzo-threshold} now proves the
inequality.  Finally, the section of a stable pair embeds \(\Ocal_C\) in
its pair sheaf with zero-dimensional cokernel.  Thus
\(\chi(F)\ge\chi(\Ocal_C)\), which proves the last assertion.
\end{proof}

\subsection{Vanishing cycles and the relative Hilbert scheme}

We assume that the PT deformation theory supplies an oriented $d$-critical
structure on $P_{\beta,n}$ and that its local vanishing-cycle sheaves glue
globally.  Write
$\Phi_{\beta,n}^{\PT}\in\operatorname{Perv}(P_{\beta,n})$ for the resulting
perverse sheaf.  No general existence theorem for this orientation is used.

The local-to-global construction for an oriented $d$-critical locus is
\cite[Theorem~6.9]{BBDJS}; related orientation choices for compactly
supported objects on Calabi--Yau threefolds are discussed in
\cite{JoyceUpmeier}.

\begin{definition}[Reduced Hilbert comparison]
\label{def:reduced-hilbert-comparison}
For the chosen global PT orientation, the condition
$(\mathrm{Hilb})_{\beta,n}$ is the assertion that, under the canonical
isomorphism of Proposition~\ref{prop:hilbert-pt},
\[
  \alpha_{\beta,n}^*
  \bigl(\Phi_{\beta,n}^{\PT}|_{P_{\beta,n}^{\red}}\bigr)
  \cong
  \IC_{\mathcal C_{\beta,\red}^{[n]}}.
\]
For $S=\Pp^2$ and $\beta=dH$, we write $(\mathrm{Hilb})_{d,n}$.
\end{definition}

\begin{remark}
Definition~\ref{def:reduced-hilbert-comparison} is a statement over the open locus
of reduced curves.  When the relative Hilbert scheme is singular, it must be
established separately.  Pullback to a closed reducible stratum does not
suffice, since closed pullback is not perverse $t$-exact and need not
preserve intersection complexes.
\end{remark}

In the plane range used below,
Theorem~\ref{thm:plane-smooth-range} proves this comparison directly.
Write
\[
  D_d=\dim|\Ocal_{\Pp^2}(d)|=\frac{d(d+3)}2,
  \qquad
  g_d=\frac{(d-1)(d-2)}2.
\]

\begin{lemma}[The zero-section threshold]\label{lem:zero-section-threshold}
Let $d\ge2$, and let $C\subset\Tot(K_{\Pp^2})$ be a compact
Cohen--Macaulay curve of degree $d$ which is not scheme-theoretically
contained in the zero section.
Then
\[
  \chi(\Ocal_C)\ge 1-g_d+d+2.
\]
The bound is sharp.  Consequently, a stable pair of degree $d$ and Euler
characteristic $1-g_d+n$ is supported on the zero section whenever
$n\le d+1$.
\end{lemma}

\begin{proof}
Let $I$ be the ideal of the zero section and $J$ the ideal of $C$.
Set-theoretically, every compact curve in $\Tot(K_{\Pp^2})$ lies in the
zero section: on the normalization of each reduced component, a displacement
from the zero section would be a section of a line bundle of negative degree.
Thus $I^{r+1}\subset J$ for some $r$.

The filtration by $J+I^k$ gives
\[
 \chi(\Ocal_C)=\sum_{k=0}^r\chi(G_k),\qquad
 G_k:=\frac{J+I^k}{J+I^{k+1}}.
\]
Since $I^k/I^{k+1}\cong\Ocal_{\Pp^2}(3k)$, there is a surjection
$\Ocal_{\Pp^2}(3k)\twoheadrightarrow G_k$.
After tensoring its kernel by
$\Ocal(-3k)$, let $K_k\subset\Ocal_{\Pp^2}$ be the resulting rank-one
torsion-free ideal.  Its reflexive hull has the form
$K_k^{\vee\vee}=\Ocal_{\Pp^2}(-D_k)$ for an effective divisor $D_k$; write
$d_k=\deg D_k$.  Comparing
$K_k\subset K_k^{\vee\vee}\subset\Ocal_{\Pp^2}$ gives
\[
 0\longrightarrow (K_k^{\vee\vee}/K_k)(3k)
 \longrightarrow G_k\longrightarrow\Ocal_{D_k}(3k)\longrightarrow0.
\]
The first term has finite length, so
\[
 \chi(G_k)\ge\chi(\Ocal_{D_k}(3k))
 =3kd_k+1-\frac{(d_k-1)(d_k-2)}2.
\]
The fundamental cycle of $C$ is the sum of those of the layers; hence
$\sum_kd_k=d$.  This proves the filtration inequality used in
\cite[proof of Proposition~2]{CKK}, with the $n=d+2$ arithmetic corrected.

If some $d_k$ with $k>0$ is nonzero, the minimum is attained for
$d_0=d-1$ and $d_1=1$.  To verify this, write
$f_k(m)=3km+1-(m-1)(m-2)/2$.  Replacing layers
$d_1=a\ge0$ and $d_k=b>0$, $k>1$, by the single layer $d_1=a+b$ changes the
sum by
$f_1(a)+f_k(b)-f_1(a+b)=ab+3(k-1)b\ge0$; iterating leaves only
$d_0=d-m$ and $d_1=m$.  For $1\le m\le d$, the difference from $m=1$
is $(m-1)(d-m+2)\ge0$.  The minimum therefore equals
\[
 1-\frac{(d-2)(d-3)}2+4=1-g_d+d+2.
\]

The bound is attained globally by
$J=I^2+I I_L+I_L^{d-1}$, where $I_L$ is the pullback of the ideal of a
line $L$.  Locally this is $(p^2,p\ell,\ell^{d-1})$, and the quotient is
free over a coordinate along $L$, with basis
$1,\ell,\ldots,\ell^{d-2},p$; it is therefore Cohen--Macaulay.  Moreover,
\[
 0\longrightarrow\Ocal_L(3)\longrightarrow\Ocal_C
 \longrightarrow\Ocal_{(d-1)L}\longrightarrow0,
\]
so $[C]=dH$ and
$\chi(\Ocal_C)=4+1-(d-2)(d-3)/2=1-g_d+d+2$.

Finally, the section of a stable pair has zero-dimensional cokernel, so
$\chi(F)\ge\chi(\Ocal_C)$.  This excludes non-zero-section support for
$n\le d+1$.
\end{proof}

\begin{lemma}[Normal deformations of a surface pair]
\label{lem:normal-deformations-surface-pair}
Let $S$ be a smooth surface, let $X=\Tot(K_S)$, and let
$s:\Ocal_S\to F$ be a stable pair on $S$.  Denote by $\mathfrak P_S$ and
$\mathfrak P_X$ the deformation functors of this pair on $S$ and on $X$,
respectively.  The natural map $\mathfrak P_S\to\mathfrak P_X$ identifies
$\mathfrak P_S$ with the closed zero-Higgs subfunctor, and there is an exact
sequence of Zariski tangent spaces
\[
 0\longrightarrow T_s\mathfrak P_S\longrightarrow T_s\mathfrak P_X
 \longrightarrow\Hom_S(F,F\otimes K_S)\longrightarrow0.
\]
Every $\theta\in\Hom_S(F,F\otimes K_S)$ integrates to a one-parameter
family of $\Ocal_X$-module structures on the fixed sheaf $F$, with Higgs
field $t\theta$.  The element $s(1)$ defines a stable-pair section for every
$t$.
\end{lemma}

\begin{proof}
Since $X\to S$ is affine, the same correspondence works in families: a
$T$-flat coherent sheaf on $X\times T$ whose support is finite over
$S\times T$ is a $T$-flat coherent sheaf $F_T$ on $S\times T$ together
with an action
\[
 K_S^{-1}\otimes F_T\longrightarrow F_T,
\]
or, by adjunction, a Higgs field $F_T\to F_T\otimes K_S$.  Vanishing of
this morphism is a closed condition and represents the surface-pair functor
as the closed zero-Higgs subfunctor.

Over the dual numbers, the Higgs field of a deformation reducing to the
zero action has the form $\epsilon\theta$.  Dividing by $\epsilon$ and
restricting to the central fibre defines the tangent map to
$\Hom_S(F,F\otimes K_S)$.  Its kernel consists exactly of first-order
surface-pair deformations.  This proves exactness at the first two terms of
the displayed sequence.

For a fixed $\theta$, contraction with $t\theta$ defines the action of the
degree-one generator $K_S^{-1}$ on $F\otimes\C[t]$.  The symmetric algebra
is freely generated by this line bundle, so iteration defines an
$\operatorname{Sym}(K_S^{-1})$-module with no additional integrability
equation.  Let $s_t:\Ocal_X\to F_t$ be the $\Ocal_X$-linear map determined
by $s_t(1)=s(1)$.  Its image contains the original $\Ocal_S$-submodule
generated by $s(1)$; hence $\operatorname{coker}(s_t)$ is a quotient of the
zero-dimensional sheaf $\operatorname{coker}(s)$.  The underlying
$\Ocal_S$-module is unchanged, so purity is also unchanged.  Thus $s_t$ is
a stable pair for every $t$.  Applying this construction over the dual
numbers proves that every normal tangent vector is attained, and hence
proves surjectivity of the last arrow.
\end{proof}

\begin{theorem}[Plane Hilbert comparison]\label{thm:plane-smooth-range}
Let $d\ge1$ and $X=\Tot(K_{\Pp^2})$.  For $0\le n\le d+1$, the stable-pair moduli
space is a smooth relative Hilbert scheme,
\[
  P_{1-g_d+n}(X,d)\cong\mathcal C_d^{[n]}.
\]
Its smooth $d$-critical structure has virtual canonical bundle
$K_{\mathcal C_d^{[n]}}^{\otimes2}$.  Choose the canonical orientation
$K_{\mathcal C_d^{[n]}}$, and let $\Phi_{d,n}^{\PT}$ denote the
vanishing-cycle sheaf for this choice.  Then
\[
  \Phi_{d,n}^{\PT}
  \cong\Q_{\mathcal C_d^{[n]}}[\dim\mathcal C_d^{[n]}]
  =\IC_{\mathcal C_d^{[n]}}.
\]
In particular, $(\mathrm{Hilb})_{d,n}$ holds.
\end{theorem}

\begin{proof}
Proposition~1 of \cite{CKK} identifies stable pairs on the surface with the
relative Hilbert scheme, as a moduli functor.  Lemma
\ref{lem:zero-section-threshold} shows that they have the same geometric
points in the stated range when $d\ge2$; for $d=1$, the same conclusion is
immediate because the degree-one Fitting divisor is a reduced line, so
Proposition~\ref{prop:reduced-zero-section} forces the Higgs field to
vanish.  By
Lemma~\ref{lem:normal-deformations-surface-pair}, the surface-pair functor is
a closed subfunctor, and a nonzero normal tangent vector
$\theta:F\to F\otimes K_{\Pp^2}$ integrates to stable pairs with Higgs field
$t\theta$.  For $t\ne0$ these pairs are not scheme-theoretically supported
on the zero section.  This contradicts
Lemma~\ref{lem:zero-section-threshold} when $d\ge2$ and
Proposition~\ref{prop:reduced-zero-section} when $d=1$.
Thus the relative tangent space is zero.

The closed immersion of the surface pair functor is consequently bijective
on geometric points and an isomorphism on Zariski tangent spaces.  At a
point, write the corresponding quotient of local rings as $A\to A/I$.
Both moduli spaces are of finite type over $\C$, and bijectivity on
geometric points makes this closed immersion a homeomorphism.  Thus
$V(I)=\operatorname{Spec}A$; since $A$ is Noetherian, $I$ is nilpotent.
Nilpotence gives $\dim A=\dim(A/I)$, and equality of tangent spaces gives
$\operatorname{embdim}A=\operatorname{embdim}(A/I)$.  The surface moduli is
smooth by the calculation below, so
$\operatorname{embdim}(A/I)=\dim(A/I)$.  Hence $A$ is regular and reduced;
the nilpotent ideal $I$ is zero.  The closed immersion is an isomorphism.

It remains to analyze the incidence scheme.  Let
$\mathcal Z\subset\Pp^2\times\Hilb^n(\Pp^2)$ be the universal
subscheme and let $p_1,p_2$ be the two projections.  The ideal of every
length-$n$ subscheme $Z\subset\Pp^2$ is $n$-regular.  Since $n\le d+1$, one
has $d\ge n-1$, and therefore
$H^1(\Pp^2,I_Z(d))=0$.  The exact sequence
$0\to I_Z(d)\to\Ocal_{\Pp^2}(d)\to\Ocal_Z\to0$ then shows that
$h^0(I_Z(d))=D_d+1-n$ for every $Z$.

Cohomology and base change now gives a locally free sheaf
$E_{d,n}=p_{2,*}(I_{\mathcal Z}\otimes p_1^*\Ocal_{\Pp^2}(d))$ of rank
$D_d+1-n$, with vanishing first higher direct image.  A point of
$\Pp(E_{d,n})$ is a pair consisting of $Z$ and a degree-$d$ equation
vanishing on $Z$.  Hence
\[
  \mathcal C_d^{[n]}\cong\Pp(E_{d,n})
  \longrightarrow\Hilb^n(\Pp^2)
\]
is a projective bundle of relative dimension $D_d-n$.  The base
$\Hilb^n(\Pp^2)$ is smooth of dimension $2n$, so the total space is smooth
of dimension $D_d+n$.

For a smooth truncation $Y$ of a $d$-critical locus, the virtual canonical
bundle is $K_Y^{\otimes2}$.  Locally, a critical chart for $Y$ is a
stabilization of the zero function on $Y$ by a nondegenerate quadratic form
in the normal directions.  The stabilization theorem for vanishing cycles
therefore identifies the perverse sheaf of an oriented chart with
$\Q_Y[\dim Y]$ tensored by the rank-one local system determined by the
orientation \cite[Theorem~6.9]{BBDJS}.  The square root $K_Y$ gives the
trivial orientation local system.  Taking
$Y=\mathcal C_d^{[n]}$ proves
$\Phi_{d,n}^{\PT}\cong\Q_Y[D_d+n]=\IC_Y$.  A different square root would
twist this sheaf by its associated rank-one sign local system.
\end{proof}

\begin{theorem}[Del Pezzo incidence range]
\label{thm:del-pezzo-incidence-range}
Let \(S\) be a del Pezzo surface, let \(\beta\ne0\) be effective, and put
\(D_\beta=\dim|\beta|\).  Fix \(n<\tau_S(\beta)\), and suppose that every
length-\(n\) subscheme \(Z\subset S\) imposes independent conditions on
\(H^0(S,\Ocal_S(\beta))\).  Then
\[
 P_{1-g_\beta+n}(\Tot(K_S),\beta)
 \cong\mathcal C_\beta^{[n]}
 \cong\Pp(E_{\beta,n})\longrightarrow\Hilb^n(S),
\]
where \(E_{\beta,n}\) has rank \(D_\beta+1-n\).  In particular the
moduli space is smooth of dimension \(D_\beta+n\).  For the orientation
\(K_{\mathcal C_\beta^{[n]}}\),
\[
 \Phi_{\beta,n}^{\PT}=\IC_{\mathcal C_\beta^{[n]}},\qquad
 \chi(B_\beta,R\Pi_{\beta,n,*}\Phi_{\beta,n}^{\PT})
 =(-1)^{D_\beta+n}(D_\beta+1-n)\chi(\Hilb^n(S)).
\]
\end{theorem}

\begin{proof}
Lemma~\ref{lem:del-pezzo-zero-section-threshold} shows that every geometric
point is a surface pair.  The surface-pair functor is the closed zero-Higgs
subfunctor by Lemma~\ref{lem:normal-deformations-surface-pair}.  If a
normal tangent vector were nonzero, the family with Higgs field
\(t\theta\) constructed in that lemma would give non-zero-section pairs
with the same numerical invariants.  Thus the closed immersion of the
surface-pair moduli space is bijective on geometric points and an
isomorphism on tangent spaces.

The surface stable-pair correspondence
\cite[Proposition~1]{CKK} identifies the surface-pair space with the
relative Hilbert scheme.  The surjective universal evaluation map has a
vector-bundle kernel of rank \(D_\beta+1-n\), so this relative Hilbert
scheme is the indicated projective bundle.  It is smooth because \(\Hilb^n(S)\) is
smooth.  At a point the closed immersion is a
quotient \(A\to A/I\) with nilpotent kernel.  Equality of tangent spaces
gives equal embedding dimensions, while \(A/I\) is regular and has the
same dimension as \(A\).  Hence \(A\) is regular and reduced, so \(I=0\).
This proves the scheme-theoretic isomorphism.  The orientation and
vanishing-cycle statement are the smooth critical-locus calculation used
in Theorem~\ref{thm:plane-smooth-range}.  Finally, a projective space of
dimension \(D_\beta-n\) has Euler characteristic \(D_\beta+1-n\), which
gives the last formula.
\end{proof}

\begin{remark}
The orientation in Theorem~\ref{thm:plane-smooth-range} is a
coefficientwise choice.  Compatibility of these choices with addition of
curve classes and derived symmetric powers is not a consequence of the
theorem; it is part of the orientation compatibility required in
Conjecture~\ref{conj:cohomological-exponential}.
\end{remark}

\begin{proposition}\label{prop:open-base-change}
If $(\mathrm{Hilb})_{\beta,n}$ holds, then
\[
  j_\beta^*R\Pi_{\beta,n,*}\Phi_{\beta,n}^{\PT}
  \cong
  R\pi_{\beta,\red,*}^{[n]}
  \IC_{\mathcal C_{\beta,\red}^{[n]}}.
\]
\end{proposition}

\begin{proof}
Let $i:P_{\beta,n}^{\red}\hookrightarrow P_{\beta,n}$ be the inverse image
of $j_\beta$, and write $\Pi_{\beta,n}^{\red}$ for the restricted support
map.  The square formed by $i$, $j_\beta$, and the two support maps is
Cartesian.  Since $j_\beta$ is open, proper base change for
$\Pi_{\beta,n}$ gives
\[
 j_\beta^*R\Pi_{\beta,n,*}\Phi_{\beta,n}^{\PT}
 \cong
 R(\Pi_{\beta,n}^{\red})_*i^*\Phi_{\beta,n}^{\PT}.
\]
Proposition~\ref{prop:reduced-zero-section} identifies the reduced PT space
with stable pairs on $S$, and Proposition~\ref{prop:hilbert-pt} identifies
this space, together with its support morphism, with
$\pi_{\beta,\red}^{[n]}:\mathcal C_{\beta,\red}^{[n]}\to
B_\beta^{\red}$.  Under this isomorphism the reduced Hilbert comparison replaces
$i^*\Phi_{\beta,n}^{\PT}$ by
$\IC_{\mathcal C_{\beta,\red}^{[n]}}$, which gives the asserted formula.
\end{proof}

\subsection{The nonreduced-support condition}

\begin{definition}[Nonreduced-support condition]
\label{def:nonreduced-support}
Put
\[
 D_\beta^{\nr}=B_\beta\setminus B_\beta^{\red}.
\]
For the chosen global PT orientation, the condition
$(\mathrm{NR})_{\beta,n}$ holds if, for every $i$, no simple
Jordan--H\"older constituent of
\[
 {}^pH^iR\Pi_{\beta,n,*}\Phi_{\beta,n}^{\PT}
\]
has support contained in $D_\beta^{\nr}$.  For Euler-characteristic
indexing, let $P_\chi(X,\beta)$ be the stable-pair moduli space with
$\chi(F)=\chi$, let $\Phi_{\chi,\beta}^{\PT}$ be its vanishing-cycle sheaf
for the chosen orientation, and let
$\Pi_{\chi,\beta}:P_\chi(X,\beta)\to B_\beta$ be the support map.  The
notation $(\mathrm{NR})_{\chi,\beta}$ denotes the same condition for
$R\Pi_{\chi,\beta,*}\Phi_{\chi,\beta}^{\PT}$.  For
$S=\Pp^2$ and $\beta=dH$, we write $(\mathrm{NR})_{d,n}$.
\end{definition}

Definition~\ref{def:nonreduced-support} concerns strict supports after proper
pushforward.  It does not say that $\Phi_{\beta,n}^{\PT}$ has zero stalks
over nonreduced curves, or that its restriction to
$\Pi_{\beta,n}^{-1}(D_\beta^{\nr})$ vanishes.  Full-support intersection
complexes usually have nonzero stalks on $D_\beta^{\nr}$.

For a constructible complex $K$, write
$K^{\mathrm{ss}}=\bigoplus_i({}^pH^iK)^{\mathrm{ss}}[-i]$.

Let $j:U\hookrightarrow B$ be a dense open immersion.  For a complex $L$ on
$U$ with semisimple perverse cohomology, define its \emph{perverse minimal
extension} by
\begin{equation}\label{eq:perverse-minimal-extension}
  {}^pj_{!*}L:=\bigoplus_i j_{!*}({}^pH^iL)[-i].
\end{equation}
Thus middle extension is applied coefficientwise to the perverse cohomology
sheaves, not directly to an arbitrary derived complex.

\begin{proposition}[Minimal-extension criterion]
\label{prop:minimal-extension-criterion}
For $K\in D_c^b(B)$, the following are equivalent:

\begin{enumerate}[label=$\roman*$,leftmargin=8mm]
\item no simple constituent of ${}^pH^iK$ is supported in $B\setminus U$;
\item $K^{\mathrm{ss}}\cong{}^pj_{!*}(j^*K)^{\mathrm{ss}}$.
\end{enumerate}

The right side is therefore the unique semisimple extension of
$(j^*K)^{\mathrm{ss}}$ having no strict support in the complement.
\end{proposition}

\begin{proof}
For every $i$, write the semisimplification of the perverse cohomology as
\[
 ({}^pH^iK)^{\mathrm{ss}}
 =\bigoplus_{\alpha}\IC_{Z_\alpha}(M_\alpha),
\]
where $Z_\alpha\subset B$ is irreducible and $M_\alpha$ is an irreducible
local system on a smooth dense open subset of $Z_\alpha$.  Since $U$ is
open, either $Z_\alpha\subset B\setminus U$ or $Z_\alpha\cap U$ is a dense
open subset of $Z_\alpha$.  In the first case the restriction of the
summand to $U$ is zero.  In the second case
$\IC_{Z_\alpha}(M_\alpha)$ is, by the defining absence of subobjects and
quotients on the boundary, the middle extension of its restriction to
$Z_\alpha\cap U$.

It follows that ${}^pj_{!*}(j^*K)^{\mathrm{ss}}$ is the direct sum of
exactly those simple constituents for which $Z_\alpha\cap U\ne\varnothing$.
It equals $K^{\mathrm{ss}}$ precisely when there is no constituent with
$Z_\alpha\subset B\setminus U$.  The same description also proves
uniqueness: any semisimple extension without a boundary constituent is the
direct sum of these prescribed middle extensions.
\end{proof}

If $(\mathrm{Hilb})_{\beta,n}$ holds, semisimplifying the perverse cohomology
gives a decomposition
\begin{equation}\label{eq:pt-minimal-extension}
 \bigl(R\Pi_{\beta,n,*}\Phi_{\beta,n}^{\PT}\bigr)^{\mathrm{ss}}
 \cong{}^pj_{\beta,!*}
 \left(R\pi_{\beta,\red,*}^{[n]}
 \IC_{\mathcal C_{\beta,\red}^{[n]}}\right)^{\mathrm{ss}}
 \oplus\mathcal N_{\beta,n}^{\PT},
\end{equation}
where every simple constituent of $\mathcal N_{\beta,n}^{\PT}$ is supported
in $D_\beta^{\nr}$.  The condition $(\mathrm{NR})_{\beta,n}$ is precisely
$\mathcal N_{\beta,n}^{\PT}=0$.  The open set $B_\beta^{\red}$ contains both
irreducible and reducible reduced curves, so all partition supports are
already present before taking the minimal extension.  An isomorphism before
semisimplification requires semisimplicity, for example from a suitable
purity statement.

\subsection{Support exclusion for plane curves}

We will use the following standard support estimate.

\begin{proposition}\label{prop:support-bound}
Let $f:Y\to B$ be projective, with $Y$ smooth and $B$ irreducible, and
$\dim Y=\dim B+r$.  Suppose that every fiber has dimension at most $r$.
If $Z\subset B$ is the support of a simple perverse constituent of
$Rf_*\Q_Y[\dim Y]$, then
\[
  \codim_B Z\le r.
\]
\end{proposition}

\begin{proof}
By the decomposition theorem \cite{BBD}, write the corresponding summand as
$\IC_Z(L)[-k]$.  Relative hard Lefschetz pairs the summands in degrees $k$
and $-k$, so we may choose an occurrence with $k\ge0$.  Let $z$ be a point
in the smooth locus of $Z$ at which $L$ is a local system, and put
$c=\codim_BZ$.  At $z$, the complex $\IC_Z(L)$ is $L[\dim Z]$.
Consequently the chosen summand gives a nonzero class in
\[
 H^{\dim Y-\dim Z+k}(Y_z,\Q)
 =H^{r+c+k}(Y_z,\Q).
\]
Every fiber has complex dimension at most $r$, hence its cohomology vanishes
in degrees greater than $2r$.  Thus $r+c+k\le2r$, or
$c\le r-k\le r$, as required.
\end{proof}

\begin{lemma}\label{lem:plane-nonreduced-codim}
For $d\ge2$, the nonreduced locus in $B_d$ has codimension $2d-1$.
\end{lemma}

\begin{proof}
Every nonreduced divisor can be written as $mE+D$, where $m\ge2$, $E$ is an
irreducible curve of degree $e$, and $me\le d$.  The locus of such divisors
is contained in the image of $B_e\times B_{d-me}$ and therefore has
codimension at least
\[
 D_d-D_e-D_{d-me}
 =\frac e2\bigl(2dm-(m^2+1)e+3(m-1)\bigr).
\]
For fixed $e$, increasing $m$ by one changes the expression by
$e(2d-(2m+1)e+3)/2$, which is positive whenever $(m+1)e\le d$.
It is therefore enough to set $m=2$.  The difference between the resulting
codimension and $2d-1$ is
$(e-1)(4d-5e-2)/2$, which is nonnegative for $1\le e\le d/2$.
The minimum is attained for $(m,e)=(2,1)$.  The image of
$B_1\times B_{d-2}\to B_d$, $(L,D)\mapsto2L+D$, has dimension
$D_1+D_{d-2}$ and hence codimension $2d-1$, proving equality.
\end{proof}

\begin{corollary}\label{cor:plane-nonreduced}
For $d\ge3$ and $0\le n\le d+1$, no simple constituent of
$R\pi_{d,*}^{[n]}\IC_{\mathcal C_d^{[n]}}$
is supported inside the nonreduced locus.  Thus
$(\mathrm{NR})_{d,n}$ holds for the smooth stable-pair model in this range.
\end{corollary}

\begin{proof}
Theorem~\ref{thm:plane-smooth-range} gives a smooth source of dimension
$D_d+n$, so the relative dimension is $n$.  Fix a curve $C\in B_d$ and
consider the Hilbert--Chow morphism
$\Hilb^n(C)\to\operatorname{Sym}^n(C)$.  The stratum of cycles
$\sum_{i=1}^s m_ip_i$ with $s$ distinct support points has dimension at most
$s$.  Over such a cycle, the punctual factors embed in
$\prod_i\Hilb^{m_i}(\mathbb A^2,0)$.  Brian\c{c}on's dimension theorem gives
dimension at most $\sum_i(m_i-1)=n-s$ for this product
\cite{Briancon}.  Hence the inverse image of each stratum has dimension at
most $s+(n-s)=n$, and therefore every fiber of $\pi_d^{[n]}$ has dimension
at most $n$.

Proposition~\ref{prop:support-bound} now bounds the codimension of any strict
support by $n$.  Lemma~\ref{lem:plane-nonreduced-codim} gives codimension
$2d-1$ for the nonreduced locus.  Since $2d-1>d+1\ge n$, no strict support
can be contained in that locus.
\end{proof}

\section{Support decomposition for relative Hilbert schemes}
\label{sec:hilbert}

We work over the reduced locus and assume throughout this section that
every curve in $B_\beta^{\red}$ is geometrically connected and that the locus
$B_\beta^{\sm}$ of geometrically connected smooth curves is dense.  These
conditions hold for plane curves.  Without connectedness, the Macdonald
denominator and the treatment of partitions with $N(\boldsymbol\beta)=0$
must be modified.  Closures of partition strata are taken in $B_\beta$;
their intersections with $B_\beta^{\red}$ are the supports relevant to the
reduced direct image.  The main input is the support decomposition of
Migliorini--Shende--Viviani \cite[Theorem~5.10]{MSV}.  We use their theorem
only after translating from their convention, in which intersection complexes
begin in cohomological degree zero, to the usual perverse normalization.
Although the theorem is written in $\ell$-adic notation, the authors allow
$\Q$-coefficients in characteristic zero.  We use its rational Betti
realization over $\C$.

\subsection{The full-support term}

Let $B_\beta^{\sm}\subset B_\beta^{\red}$ be the smooth-curve locus,
let $\rho_\beta:\mathcal C_\beta^{\sm}\to B_\beta^{\sm}$ be the universal
smooth curve, and put $V_\beta=R^1\rho_{\beta,*}\Q$.
If $M$ is a local system on a dense smooth open subset of $B_\beta$, then
$\IC_{B_\beta}(M)$ denotes its middle-extension intersection complex on
$B_\beta$.

For a local system $V$ and $0\le k\le2n$, define its finite Macdonald
coefficient by

\begin{equation}\label{eq:finite-macdonald}
  \Mac_n^k(V)
  =[x^ny^k]\,
  \frac{\sum_{a\ge0}\bigwedge^aV\,x^ay^a}
       {(1-x)(1-xy^2)}.
\end{equation}

Macdonald's formula \cite{Macdonald} gives
$R^k\rho_{\beta,*}^{[n]}\Q\cong\Mac_n^k(V_\beta)$, where
$\rho_\beta^{[n]}$ is the relative $n$-th symmetric product.
We therefore set

\begin{equation}\label{eq:smooth-support}
  \mathcal S_{\beta,n}
  =\bigoplus_{k=0}^{2n}
  \IC_{B_\beta}
  \bigl(\Mac_n^k(V_\beta)\bigr)[n-k].
\end{equation}

Its restriction to $B_\beta^{\red}$ is the corresponding middle extension
inside the reduced locus.  If $\mathcal C_{\beta,\red}^{[n]}$ is smooth,
this restriction is the sum of all full-$B_\beta^{\red}$-support
constituents of
$R\pi_{\beta,\red,*}^{[n]}\Q[\dim B_\beta+n]$.

Indeed, the restriction of the direct image to $B_\beta^{\sm}$ is
$
\bigoplus_k\Mac_n^k(V_\beta)[\dim B_\beta+n-k]
$, and a full-support simple constituent is the middle extension of its
restriction to this open set.

\subsection{Reducible strata and partial normalization}

Let $\boldsymbol\beta=(\beta_1,\ldots,\beta_\ell)$ be a decomposition
$\beta_1+\cdots+\beta_\ell=\beta$ into nonzero effective classes.  Let
$\Sigma_{\boldsymbol\beta}^{\circ}\subset B_\beta^{\red}$ be the locus of
curves $C=C_1\cup\cdots\cup C_\ell$ such that $C_i$ is a smooth connected
curve of class $\beta_i$, all
intersections are transverse, and no three components meet at one point.  We
assume that this locus is nonempty.  Its component-labelled cover is an open
subset $\widetilde\Sigma_{\boldsymbol\beta}^{\circ}
\subset\prod_{i=1}^{\ell}B_{\beta_i}^{\sm}$.

Decompositions are unordered; whenever a direct sum is indexed by
$\boldsymbol\beta$, one representative is chosen for each unordered
decomposition.

The group which permutes equal classes acts on this cover.  Set
$N(\boldsymbol\beta)=\sum_{i<j}\beta_i\cdot\beta_j$.

On the labelled cover, normalize all nodes joining different components.
The resulting family is the disjoint union
$\widetilde{\mathcal C}_{\boldsymbol\beta}^{\nu}
=\coprod_{i=1}^{\ell}\widetilde{\mathcal C}_i$ over
$\widetilde\Sigma_{\boldsymbol\beta}^{\circ}$.

For $s\ge0$, its relative Hilbert scheme is

\begin{equation}\label{eq:partial-normalization-hilbert}
  \bigl(\widetilde{\mathcal C}_{\boldsymbol\beta}^{\nu}\bigr)^{[s]}
  =
  \coprod_{r_1+\cdots+r_\ell=s}
  \prod_{i=1}^{\ell}\widetilde{\mathcal C}_i^{[r_i]}.
\end{equation}

Let $\widetilde\rho_{\boldsymbol\beta,s}$ denote its projection to the
labelled stratum.  The equivariant support decomposition gives descent isomorphisms
on $R^k\widetilde\rho_{\boldsymbol\beta,s,*}\Q$.

We denote the descended local system on
$\Sigma_{\boldsymbol\beta}^{\circ}$ by
$\mathcal H_{\boldsymbol\beta,s}^k$.  When no two $\beta_i$ are equal,
there is no component-permutation quotient and its pullback is simply

\begin{equation}\label{eq:partial-normalization-local-system}
 \bigoplus_{\substack{r_1+\cdots+r_\ell=s\\
                       k_1+\cdots+k_\ell=k}}
 \bigboxtimes_{i=1}^{\ell}
 R^{k_i}\rho_{\beta_i,*}^{[r_i]}\Q.
\end{equation}

For repeated classes, these descent isomorphisms are part of the equivariant MSV
decomposition; it must not be replaced by taking invariants without keeping
the action on the normal directions.  This distinction produces, for example,
the sign local system on the reducible conic locus.

Put
$Z_{\boldsymbol\beta}=\overline{\Sigma_{\boldsymbol\beta}^{\circ}}
\subset B_\beta$.

For $s\ge0$, define the standard-perverse-normalized partial-normalization
term

\begin{equation}\label{eq:boundary-term}
  \mathcal T_{\boldsymbol\beta,s}
  =\bigoplus_{k=0}^{2s}
  \IC_{Z_{\boldsymbol\beta}}
  \bigl(\mathcal H_{\boldsymbol\beta,s}^k\bigr)[s-k].
\end{equation}

If mixed Hodge modules are retained, this term carries the Tate twist
$N(\boldsymbol\beta)$ in the convention of \cite{MSV}.  We work with
rational constructible complexes, where this twist is invisible.

\begin{lemma}[Equivariant descent on a partition stratum]
\label{lem:partition-descent}
After replacing the component-labelled cover by a finite \`etale Galois
closure, the local MSV strict-support decomposition is equivariant under the
deck group in each perverse cohomology sheaf.  Consequently the
partial-normalization local systems in
\eqref{eq:partial-normalization-local-system}, summed over the deck orbit of
an unordered partition, descend to the local systems
$\mathcal H_{\boldsymbol\beta,s}^k$ on
$\Sigma_{\boldsymbol\beta}^{\circ}$.
\end{lemma}

\begin{proof}
Let $G$ be the deck group of the chosen Galois closure.  Pullback by an
element of $G$ permutes the labelled components and the nodes joining them,
and therefore preserves the local MSV decomposition while permuting its
simple constituents.  In a fixed perverse cohomology sheaf, take the direct
sum of all constituents whose strict supports belong to the $G$-orbit of the
chosen labelled partition.  Semisimplicity makes this the unique maximal
subobject having precisely these strict supports.  It is consequently
$G$-stable, and the natural pullback isomorphisms satisfy the cocycle
condition inherited from the original direct image.

This maximal-strict-support characterization makes the summand independent
of the splitting used in the decomposition theorem.  Finite \`etale descent
produces a complex on the unlabelled stratum.  On
cohomology in degree $k$ its pullback is the K\"unneth local system in
\eqref{eq:partial-normalization-local-system}; this defines
$\mathcal H_{\boldsymbol\beta,s}^k$.  When classes repeat, the $G$-action
also permutes the smoothing coordinates normal to the stratum.  Retaining
this action is what determines the descended local system.
\end{proof}

\subsection{The local MSV calculation}

The global formula of Migliorini--Shende--Viviani extends to an
independently broken H-smooth family
\cite[Definition~5.14 and Corollary~5.16]{MSV}.  The definition requires
nonsingular total spaces for all relative Hilbert schemes and the existence
of a relative fine compactified Jacobian, the condition on the component
space in
\cite[Proposition~5.4]{MSV}, and a dense nodal open set in each cogenus
stratum.  The complete plane family fails the first requirement: the
relative four-point Hilbert scheme of conics is singular by
Proposition~\ref{prop:conic-correction-node}.  Component permutations are
handled in \cite{MSV} after a finite \`etale cover and are not themselves an
obstruction to H-smoothness.  We therefore do not apply the global formula
to the complete plane family.

The generic complex on a fixed reducible support is determined by the versal
local calculation of \cite[Theorem~5.10]{MSV}.

\begin{lemma}[Generic transverse nodal pullback]
\label{lem:transverse-nodal-base-change}
Let $(\mathcal C_V,C)\to(V,0)$ be a labelled versal deformation of a
reduced nodal curve, and let $t_1,\ldots,t_N$ be the smoothing coordinates
of its nodes.  Let $g:(B,b)\to(V,0)$ classify a deformation of $C$, with
$B$ smooth, and suppose that
\[
 d(t_1\circ g,\ldots,t_N\circ g)_b:T_bB\longrightarrow\C^N
\]
is surjective.  After shrinking the germs, normalized pullback carries each
MSV summand supported on a node-persistence stratum as follows: at a general
point of an irreducible component of its inverse image, and after shrinking
about that point, the pullback is the intersection complex of that component
with the pulled-back partial-normalization local system.  Moreover,
\[
 g^*R\pi_{V,*}^{[n]}\Q[\dim V+n][\dim B-\dim V]
 \cong R\pi_{B,*}^{[n]}\Q[\dim B+n].
\]
\end{lemma}

\begin{proof}
For $I\subset\{1,\ldots,N\}$, write
\[
 Z_I^\circ=\{t_i=0\ (i\in I),\ t_j\ne0\ (j\notin I)\}.
\]
Near a general point of $Z_I^\circ$, the remaining components of the curve
are smooth and the partial-normalization cohomology is a local system.
Shrink to this open set.  The versal base has product coordinates
$(t_1,\ldots,t_N,y_1,\ldots,y_m)$, and the node-persistence loci are the
corresponding coordinate subspaces.

After shrinking $B$, the map
$(t_1\circ g,\ldots,t_N\circ g)$ is a submersion.  The constant-rank theorem
makes these functions part of a coordinate system on $B$.  Hence the inverse
image of each node-persistence locus is smooth of the expected codimension
at the general point of every relevant component, and $g$ is transverse
there.  If $W$ is such a component, transverse pullback gives
\[
 g^*\IC_Z(L)[\dim B-\dim V]
 \cong\IC_W(g^*L)
\]
for the corresponding summand.  Removing the proper closed set where the
partial normalization ceases to be smooth gives the stated generic
assertion.  Finally, the relative
Hilbert scheme of the pulled-back curve is the base change of
$\mathcal C_V^{[n]}$.  Proper base change gives the second display, and the
shift $\dim B-\dim V$ changes the source normalization from
$[\dim V+n]$ to $[\dim B+n]$.
\end{proof}

\begin{proposition}[Generic partition summand]
\label{prop:generic-partition}
Fix $n\ge0$ and a decomposition $\boldsymbol\beta$ with
$N(\boldsymbol\beta)\le n$.  Suppose that
$\mathcal C_{\beta,\red}^{[n]}$ is smooth and that, at the generic point of
$\Sigma_{\boldsymbol\beta}^{\circ}$, the Kodaira--Spencer map surjects onto
the smoothing spaces of the nodes joining distinct components.  Then the
strict-support summand of
$R\pi_{\beta,\red,*}^{[n]}
\IC_{\mathcal C_{\beta,\red}^{[n]}}$ with support
$Z_{\boldsymbol\beta}\cap B_\beta^{\red}$ is
$j_\beta^*\mathcal T_{\boldsymbol\beta,n-N(\boldsymbol\beta)}$.
\end{proposition}

\begin{proof}
Put $N=N(\boldsymbol\beta)$.  At a general point of the stratum, the
components are smooth and their only singularities are the $N$ transverse
nodes joining distinct components.  After an analytic shrinking, choose a
labelled versal deformation of this projective curve and a classifying map
$g$ from the germ of $B_\beta$ to its versal base.  The node-smoothing
coordinates on that base pull back to the local smoothing parameters of the
plane family.  The assumed Kodaira--Spencer surjectivity therefore makes
$g$ transverse to every stratum on which a prescribed subset of the nodes
persists.

Theorem~5.10 of \cite{MSV} gives the strict-support decomposition on the
versal base.  Proper base change pulls back the full direct image, while
Lemma~\ref{lem:transverse-nodal-base-change} identifies, near a general
point of each inverse-image component, the corresponding
partial-normalization local system.
This does not make the complete plane family versal or globally
independently broken H-smooth.  Only transversality to the nodal strata is
used; the directions in which the components themselves vary remain in the
base and are recorded by the pulled-back cohomology local systems.

The local support formula attaches to the codimension-$N$ stratum in the
$n$-th relative Hilbert scheme the cohomology of the length-$s$ Hilbert
scheme of the partial normalization, where $s=n-N$.  On the
component-labelled cover the partial normalization is
$\coprod_i\widetilde{\mathcal C}_i$, so its length-$s$ Hilbert scheme
decomposes as in \eqref{eq:partial-normalization-hilbert}.  Applying K\"unneth
to each component gives, in cohomological degree $k$, the local system
displayed in \eqref{eq:partial-normalization-local-system}.  The perverse
normalization in the local theorem is the shift $[s-k]$; its Tate twist by
$N$ disappears after passing to rational constructible complexes.

The components become labelled after a finite \`etale base change; see
\cite[Proposition~5.4]{MSV}.  Lemma~\ref{lem:partition-descent} then
descends the local systems to $\mathcal H_{\boldsymbol\beta,s}^k$, including
the permutation action for repeated classes.  Finally, the source is smooth and
the relative Hilbert morphism is projective, so the decomposition theorem
applies.  A summand with strict support
$Z_{\boldsymbol\beta}\cap B_\beta^{\red}$ is uniquely the middle extension
of its restriction to the dense smooth stratum.  Summing over $k$ with
shifts $[s-k]$ gives $j_\beta^*\mathcal T_{\boldsymbol\beta,s}$.  The
argument determines this strict support and does not classify summands
supported on proper closed subsets of it.
\end{proof}

If $(\mathrm{Hilb})_{\beta,n}$ holds, Proposition
\ref{prop:open-base-change} transfers this statement to the PT direct image
over the reduced locus.  If $(\mathrm{NR})_{\beta,n}$ also holds,
Proposition~\ref{prop:minimal-extension-criterion} gives its unique middle extension to
$B_\beta$, after semisimplification.  These two assumptions still do not
exclude additional strict supports on special reduced loci; that stronger
claim is part of the exponential conjecture below.

\section{Symmetric products on the Chow monoid}
\label{sec:chow-exponential}

The symmetric algebra over the Chow variety was introduced by
Nekrasov--Okounkov \cite[Section~2.3.5]{NekrasovOkounkov}; see also
\cite[Section~5.3.3]{OkounkovTakagi}.  The derived constructible form used
here is as follows.

We retain the connectedness and density hypotheses of
Section~\ref{sec:hilbert} and fix an additive subsemigroup
$\Gamma\subset\Eff(S)$ on which they hold; write
$\Gamma_+=\Gamma\setminus\{0\}$.  For local $\Pp^2$ one takes
$\Gamma=\{dH:d\ge0\}$.  Put $B_0=\operatorname{Spec}\C$.  Addition of effective
cycles gives finite, hence proper, maps
$\mu_{\alpha,\beta}:B_\alpha\times B_\beta\to B_{\alpha+\beta}$,
since a fixed effective cycle has only finitely many allocations of its
irreducible components between two subcycles.  Thus
$\coprod_{\beta\in\Gamma}B_\beta$ is a graded commutative monoid.

For each $\beta$, let $\mathsf C_\beta$ be the additive idempotent-complete
subcategory of $D_c^b(B_\beta,\Q)$ whose objects are finite direct sums of
shifts of semisimple perverse sheaves.  We use the formal sums
$\widehat{\mathsf C}_\Gamma=
\prod_{\beta\in\Gamma}\mathsf C_\beta((q))Q^\beta$, with the $q$-powers
bounded below for each fixed $\beta$.  The effective-class grading is
locally finite: only finitely many multiplicity functions occur in a fixed
coefficient.  For an arbitrary constructible complex, we first take the
coefficientwise semisimplification defined in Section~\ref{sec:pt}.  All
isomorphisms in this section are read in $\widehat{\mathsf C}_\Gamma$.

For
series $A=\bigoplus A_\beta Q^\beta$ and
$A'=\bigoplus A'_\beta Q^\beta$, define convolution by
\begin{equation}\label{eq:chow-convolution}
 (A\star A')_\beta
 =\bigoplus_{\alpha+\gamma=\beta}
 R\mu_{\alpha,\gamma,*}(A_\alpha\boxtimes A'_\gamma).
\end{equation}
For a del Pezzo surface the required local finiteness follows because its
effective cone is rational polyhedral and $-K_S$ is positive on every
nonzero effective class.

For $F\in D_c^b(B_\beta,\Q)((q))$, let
$\sigma_m:B_\beta^m\to\operatorname{Sym}^mB_\beta$ be the quotient and set
$\operatorname{Sym}^m(F)=(R\sigma_{m,*}F^{\boxtimes m})^{\mathfrak S_m}$.
The action is the canonical Koszul action on complexes.  Since the
coefficients are rational, taking invariants is a direct summand.  For a
$q$-series this definition is taken in the completed category of graded
complexes, with every cohomological shift retained inside the symmetric
power \cite{MaximSaitoSchurmann}.

For local $\Pp^2$, a partition
$\lambda=(1^{m_1}2^{m_2}\cdots)\vdash d$ determines the multiplication map
\begin{equation}\label{eq:chow-partition-map}
 \mu_\lambda:
 \prod_{e\ge1}\operatorname{Sym}^{m_e}B_e\longrightarrow B_d,
 \qquad ([f_{e,j}])\longmapsto\left[\prod_{e,j}f_{e,j}\right].
\end{equation}
Thus a series $A=\bigoplus_{e\ge1}A_eQ^e$ has the Chow exponential
\begin{equation}\label{eq:chow-exponential-definition}
 \left.\ExpChow(A)\right|_{B_d}
 =\bigoplus_{\lambda=(1^{m_1}2^{m_2}\cdots)\vdash d}
 R\mu_{\lambda,*}
 \left(\bigboxtimes_{e\ge1}\operatorname{Sym}^{m_e}(A_e)\right),
\end{equation}
with degree-zero term $\Q_{B_0}$.  This is the symmetric product of a
formal series whose degree-$e$ coefficient lives on the degree-$e$ Chow
variety.  It includes both the diagonals and the permutation monodromy of
equal factors.  For a general effective-class semigroup, one replaces
$\lambda$ by a finite multiplicity function $(m_\gamma)$ satisfying
$\sum_\gamma m_\gamma\gamma=\beta$ and uses the corresponding cycle-addition
map, denoted $\mu_{(m_\gamma)}$; formula
\eqref{eq:chow-exponential-definition} is unchanged.

\begin{proposition}[Exponential law]\label{prop:chow-exponential-law}
There is a natural isomorphism
$\ExpChow(A\oplus A')\cong\ExpChow(A)\star\ExpChow(A')$.
\end{proposition}

\begin{proof}
Fix $e$ and $m$.  Expanding the external tensor power gives a direct sum
indexed by subsets $I\subset\{1,\ldots,m\}$: the positions in $I$ carry
$A_e$ and the remaining positions carry $A'_e$.  The symmetric group acts
transitively on subsets of cardinality $i$, with stabilizer
$\mathfrak S_i\times\mathfrak S_{m-i}$.  Taking invariants, with the Koszul
action on the tensor factors, therefore gives
\[
 \operatorname{Sym}^m(A_e\oplus A'_e)
 \cong\bigoplus_{i+j=m}
 R\nu_{i,j,*}\bigl(\operatorname{Sym}^i(A_e)
 \boxtimes\operatorname{Sym}^j(A'_e)\bigr),
\]
where $\nu_{i,j}:\operatorname{Sym}^iB_e\times
\operatorname{Sym}^jB_e\to\operatorname{Sym}^mB_e$ adds the two cycles.
This is the usual induction--invariants identity over $\Q$.

Apply this identity independently in every degree $e$.  A choice of all
pairs $(i_e,j_e)$ separates a partition coefficient into a class
$\alpha=\sum_e e i_e$ coming from $A$ and a class
$\gamma=\sum_e e j_e$ coming from $A'$.  The composite of the maps
$\nu_{i_e,j_e}$ with cycle addition is
$\mu_{\alpha,\gamma}$.  Summing over $\alpha+\gamma=\beta$ gives precisely
the convolution formula \eqref{eq:chow-convolution}.
\end{proof}

Put $g_\beta=p_a(\beta)$ and define
the full-support GV complex and the universal point factor by
\begin{equation}\label{eq:gv-complex}
 \mathcal G_\beta(q)
 =\bigoplus_{a=0}^{2g_\beta}q^{a+1-g_\beta}
   \IC_{B_\beta}(\wedge^aV_\beta)[g_\beta-1-a],
 \qquad
 D(q)=\frac1{(1-q)(1-q[-2])}.
\end{equation}
Here $q$ is an even grading variable and $q^rK[s]$ retains the displayed
cohomological shift as part of its coefficient.  The terminology is
motivated by the vanishing-cycle definition of Gopakumar--Vafa invariants in
\cite{MaulikToda}.  For toric del Pezzo surfaces and ample $\beta$,
Maulik--Shen identify the underlying full-support summands with the direct
image from the moduli of one-dimensional semistable sheaves
\cite[Theorem~0.4]{MaulikShen}; $\mathcal G_\beta(q)$ is their PT-normalized,
$q$-graded reorganization.  Outside this setting we use
\eqref{eq:gv-complex} as the definition and do not assert a separate
identification with a sheaf-moduli vanishing-cycle complex.

Macdonald's formula \eqref{eq:finite-macdonald} is equivalent to
\begin{equation}\label{eq:normalized-macdonald-complex}
 \mathcal A_\beta(q):=D(q)\otimes\mathcal G_\beta(q)
 =\bigoplus_{n\ge0}q^{n+1-g_\beta}
   \mathcal S_{\beta,n}[g_\beta-1-n].
\end{equation}
The displayed shift is $-\chi(F)$ for a stable pair with
$\chi(F)=1-g_\beta+n$ and is retained inside every symmetric power.

For arbitrary Euler characteristic, denote the oriented vanishing-cycle
sheaf on $P_\chi(X,\beta)$ by $\Phi_{\chi,\beta}^{\PT}$.  Define the
normalized PT series by
\begin{equation}\label{eq:normalized-pt-series}
 \mathcal Z_{\PT}(q,Q)
 =1+\bigoplus_{\beta\in\Gamma_+}\ \bigoplus_{\chi\in\mathbb Z}
 q^\chi
 \bigl(R\Pi_{\chi,\beta,*}\Phi_{\chi,\beta}^{\PT}\bigr)^{\mathrm{ss}}
 [-\chi]Q^\beta.
\end{equation}
For $\chi=1-g_\beta+n$, this is the normalization in
\eqref{eq:normalized-macdonald-complex}.  If $n<0$, the restriction to
$B_\beta^{\red}$ is empty because the surface stable-pair correspondence
would identify $n$ with the length of the cokernel; under
$(\mathrm{NR})_{\chi,\beta}$ its semisimplified direct image is therefore
zero.  Thus the series reduces to $n\ge0$ when these support conditions hold.
If purity makes the direct images semisimple, the superscript ${\mathrm{ss}}$
can be omitted.

The conjecture requires a compatible choice of orientations for all
$(\chi,\beta)$.  On the open set of pairwise disjoint reduced supports, this
means that the factorization isomorphisms lift to the oriented $d$-critical
charts, satisfy the associativity and permutation cocycles, and induce on
vanishing cycles the Koszul signs used in
\eqref{eq:chow-exponential-definition}.  Orientations for one coefficient,
such as those used in Sections~\ref{sec:plane}
and~\ref{sec:degree-two-correction}, do
not by themselves give this compatibility between coefficients.

\begin{conjecture}[Cohomological GV/PT exponential]
\label{conj:cohomological-exponential}
There is a coefficientwise isomorphism in the completed semisimple graded
category
\begin{equation}\label{eq:cohomological-exponential}
 \mathcal Z_{\PT}(q,Q)
 \cong
 \ExpChow\!\left(
   \bigoplus_{\beta\in\Gamma_+}\mathcal A_\beta(q)Q^\beta\right).
\end{equation}
\end{conjecture}

For an independently broken H-smooth system, the local MSV support formulas
globalize by \cite[Corollary~5.16]{MSV}.
Conjecture~\ref{conj:cohomological-exponential} is their Chow-valued
factorization form for
the complete local-surface family; it is not a consequence of the two PT
sheaf conditions in Section~\ref{sec:pt}.  Those conditions compare the PT
coefficient with the relative Hilbert scheme on the reduced locus and control
extension across the nonreduced complement.  They do not exclude special
strict supports inside the reduced locus.

We impose the following moving-factor condition on $\Gamma$:
\emph{(MF) for every finite multiplicity function $(m_\gamma)$, the source
$\prod_\gamma\operatorname{Sym}^{m_\gamma}B_\gamma$ has a dense open subset
whose cycle sum is reduced.}
It holds for the complete linear systems $B_d=|\Ocal_{\Pp^2}(d)|$, since a
general tuple consists of distinct smooth curves with no common component.

\begin{proposition}[Minimal extension of the Chow exponential]
\label{prop:exponential-minimal-extension}
Assume \emph{(MF)}.
Let $\mathcal R_{\chi,\beta}$ be a coefficient of
$\ExpChow(\bigoplus_{\gamma\in\Gamma_+}
\mathcal A_\gamma(q)Q^\gamma)$.  This coefficient is semisimple, and
\begin{equation}\label{eq:exponential-is-minimal}
 \mathcal R_{\chi,\beta}
 \cong{}^pj_{\beta,!*}j_\beta^*\mathcal R_{\chi,\beta}.
\end{equation}
In particular, no simple constituent of $\mathcal R_{\chi,\beta}$ is
supported entirely on $D_\beta^{\nr}$.
\end{proposition}

\begin{proof}
Fix a coefficient $(\chi,\beta)$ and one multiplicity function
$(m_\gamma)$ contributing to it.  Each coefficient of
$\mathcal A_\gamma(q)$ is a finite direct sum of shifts of
$\IC_{B_\gamma}(M)$, with $M$ a polarizable local system on the smooth-curve
locus.  These local systems are exterior powers of a polarized variation of
Hodge structure, so their IC extensions underlie polarizable pure Hodge
modules.  External products of these summands are intersection complexes on
the corresponding products.  The quotient
$B_\gamma^{m_\gamma}\to\operatorname{Sym}^{m_\gamma}B_\gamma$ is finite,
and the derived symmetric power is the direct summand selected by the
Koszul $\mathfrak S_{m_\gamma}$-action
\cite{MaximSaitoSchurmann}.  Purity and the decomposition theorem therefore
make every such symmetric-power coefficient semisimple.

The product of the quotient maps and the factorization map
$\mu_{(m_\gamma)}$ are finite.  A finite map is small onto its image, so the
pushforward of each nonzero simple summand is the IC extension of its
generic finite-pushforward local system.  In particular, its strict support
is the image of the relevant irreducible component of the factorization
source, rather than a proper closed subset of that image.

By the moving-factor hypothesis, this source contains a dense open subset
whose total divisor is reduced.  Its image is dense in the strict support,
so the support meets $B_\beta^{\red}$ densely.  The same holds for every
multiplicity function and every simple constituent.  Proposition
\ref{prop:minimal-extension-criterion} now gives
\eqref{eq:exponential-is-minimal} coefficientwise.
\end{proof}

\begin{corollary}[Extension from the reduced locus]
\label{cor:reduced-exponential}
Assume \emph{(MF)}, and suppose the restriction of
\eqref{eq:cohomological-exponential} to every
$B_\beta^{\red}$ is known.  Then the full Chow-valued identity holds if and
only if $(\mathrm{NR})_{\chi,\beta}$ holds for every coefficient, in the
completed semisimple graded category.  Equivalently, defining the
semisimplified PT pushforward to be the perverse minimal extension of its
reduced restriction recovers Conjecture~\ref{conj:cohomological-exponential}.
\end{corollary}

\begin{proof}
Fix $(\chi,\beta)$.  By the assumed identity on $B_\beta^{\red}$, the
restrictions of the PT coefficient and of the Chow-exponential coefficient
are isomorphic.  Proposition~\ref{prop:exponential-minimal-extension} says
that the latter is the middle extension of this common restriction.  On the
PT side, Proposition~\ref{prop:minimal-extension-criterion} decomposes the
semisimplification as this same middle extension plus the sum of all simple
constituents supported in $D_\beta^{\nr}$.  The two global coefficients are
therefore isomorphic exactly when that latter sum is zero, which is
$(\mathrm{NR})_{\chi,\beta}$.  Applying the argument to every coefficient
proves the assertion in the completed category.
\end{proof}

\subsection{The refined PT/GV realization}

For a projective variety $B$ and $K\in D_c^b(B,\Q)$, put
$P_B(K;r)=\sum_k\dim\mathbb H^k(B,K)r^k$ and
$\chi_j(z)=\sum_{m=-j}^{j}z^{2m}$.  Extend $P_B$ coefficientwise to
Chow-valued series and denote the resulting realization simply by $P$.
The symbol $\mathbb L^{1/2}$ below is a formal grading variable, chosen to
match the motivic notation of \cite{CKK}.  For
$-g_\beta\le i\le g_\beta$, set
$\mathcal P_{\beta,i}=\IC_{B_\beta}(\wedge^{g_\beta+i}V_\beta)$.
Relative hard Lefschetz and hard Lefschetz on the projective Chow variety
give a unique character decomposition
\begin{equation}\label{eq:refined-gv-character}
 \sum_{i,k}\dim\mathbb H^k(B_\beta,\mathcal P_{\beta,i})x^iy^k
 =\sum_{j_L,j_R}N^\beta_{j_L,j_R}
   \chi_{j_L}(x)\chi_{j_R}(y).
\end{equation}
Thus the perverse degree records the $\mathrm{SU}(2)_L$ weight and the
hypercohomological degree records the $\mathrm{SU}(2)_R$ weight; here
$j_L,j_R\in\frac12\mathbb Z_{\ge0}$, and all weight products below run in
steps of one.

Define the refined realization of the PT direct images by
\begin{equation}\label{eq:refined-pt-realization}
 Z_{\PT}^{\mathrm{ref}}(\mathfrak q,Q;\mathbb L^{1/2})
 =1+\sum_{\beta\in\Gamma_+}\sum_{\chi\in\mathbb Z}
 \mathfrak q^\chi Q^\beta
 P_{B_\beta}\!\left(
   (R\Pi_{\chi,\beta,*}\Phi_{\chi,\beta}^{\PT})^{\mathrm{ss}};
   -\mathbb L^{1/2}\right).
\end{equation}

\begin{proposition}[Refined realization of the Chow exponential]
\label{prop:refined-realization}
Assume Conjecture~\ref{conj:cohomological-exponential}.  The local system
$V_\beta$ is a polarized variation of Hodge structure, so its exterior
powers and their IC extensions underlie polarizable pure Hodge modules
\cite{SaitoHodge}.  Its polarization supplies the first Lefschetz action,
and an ample class on the projective variety $B_\beta$ supplies the second.
The two actions commute.  Then
\begin{align}\label{eq:ckk-product}
 Z_{\PT}^{\mathrm{ref}}
 ={}&\prod_{\beta\in\Gamma_+}\prod_{j_L,j_R}
 \prod_{m_L=-j_L}^{j_L}\prod_{m_R=-j_R}^{j_R}
 \prod_{m=1}^{\infty}\prod_{a=0}^{m-1} \notag\\[-2mm]
 &\left(1-\mathbb L^{-m/2+1/2+a-m_R}
 (-\mathfrak q)^{m-2m_L}Q^\beta\right)^{
 (-1)^{2(j_L+j_R)}N^\beta_{j_L,j_R}}.
\end{align}
This is the refined PT/GV product of Choi--Katz--Klemm
\cite[Equation~(8.1)]{CKK}.
\end{proposition}

\begin{proof}
We first record the realization of a derived symmetric algebra.  Proper
pushforward and K\"unneth identify the hypercohomology of $\ExpChow(A)$ with
the super-symmetric algebra on
$\bigoplus_{\beta,\chi}\mathbb H^*(B_\beta,A_{\beta,\chi})q^\chi Q^\beta$.
A one-dimensional even class of monomial weight $z$ contributes
$(1-z)^{-1}$, whereas an odd class contributes $1+z$.  Hence
\begin{equation}\label{eq:cohomological-product-realization}
 P(\ExpChow(A);r)
 =\prod_{\beta,\chi,k}
 \left(1-(-1)^kq^\chi Q^\beta r^k\right)^{
 (-1)^{k+1}\dim\mathbb H^k(B_\beta,A_{\beta,\chi})}.
\end{equation}

Rename the exterior-power index in \eqref{eq:gv-complex} as $b$ and put
$i=b-g_\beta$.  Then
\begin{equation}\label{eq:connected-refined-expansion}
 \mathcal G_\beta(q)
 =\sum_{i=-g_\beta}^{g_\beta}
 q^{i+1}\mathcal P_{\beta,i}[-i-1],
 \qquad
 D(q)=\sum_{m\ge1}\sum_{a=0}^{m-1}q^{m-1}[-2a].
\end{equation}
The second identity follows by expanding the two geometric series in
$((1-q)(1-q[-2]))^{-1}$ and grouping terms with fixed total exponent
$m-1$.  Consequently, the $(m,a,i)$ term of $\mathcal A_\beta(q)$ is
$q^{m+i}\mathcal P_{\beta,i}[-i-1-2a]$.  A class in
$\mathbb H^k(B_\beta,\mathcal P_{\beta,i})$ therefore has total
cohomological degree
\[
 t=i+k+1+2a
\]
in the connected complex.  Its contribution to
\eqref{eq:cohomological-product-realization} is
\begin{equation}\label{eq:single-refined-factor}
 \left(1-(-1)^tq^{m+i}Q^\beta r^t\right)^{
 (-1)^{t+1}\dim\mathbb H^k(B_\beta,\mathcal P_{\beta,i})}.
\end{equation}

The normalization on the PT side is compatible with the substitution
$\mathfrak q=qr$, because
\[
 q^\chi P_B(K[-\chi];r)
 =q^\chi r^\chi P_B(K;r)
 =\mathfrak q^\chi P_B(K;r).
\]
We may therefore
apply \eqref{eq:single-refined-factor} to the normalized series.
Use the two Lefschetz symmetries to replace $(i,k)$ by $(-i,-k)$.
After writing $q=\mathfrak q/r$, the expression inside the parentheses
becomes
\[
 1-(-1)^{-i-k+1}\mathfrak q^{m-i}Q^\beta
 r^{-k+1+2a-m},
\]
and its outer parity is $(-1)^{i+k}$.  Set
$r=-\mathbb L^{1/2}$.  The exponent of $\mathbb L$ is
$-m/2+1/2+a-k/2$, while the total sign inside the factor is
$(-1)^{m-i}$; the latter combines with $\mathfrak q^{m-i}$ to give
$(-\mathfrak q)^{m-i}$.

Finally write $i=2m_L$ and $k=2m_R$.  Within the character
$\chi_{j_L}(x)\chi_{j_R}(y)$, the parity $i+k$ is constant and equals
$2(j_L+j_R)$ modulo two.  Thus each weight contributes the factor in
\eqref{eq:ckk-product} with exponent $(-1)^{2(j_L+j_R)}$.  Taking all
weights and then the multiplicity $N^\beta_{j_L,j_R}$ proves the product
formula.
\end{proof}

\begin{corollary}[Refined realization of the minimal extension]
\label{cor:minimal-refined-product}
Assume \emph{(MF)} and that \eqref{eq:cohomological-exponential} is known after restriction
to every reduced locus, as in Corollary~\ref{cor:reduced-exponential}.
The series obtained by replacing every PT pushforward with the perverse
minimal extension of its reduced restriction satisfies
\eqref{eq:ckk-product}.  If all the conditions
$(\mathrm{NR})_{\chi,\beta}$ hold, this series is the refined realization of
the semisimplified PT vanishing-cycle pushforwards.
\end{corollary}

\begin{proof}
Corollary~\ref{cor:reduced-exponential} identifies the minimal-extension
series with the Chow exponential.  Apply
Proposition~\ref{prop:refined-realization}; under $(\mathrm{NR})$ the
minimal-extension and PT coefficients agree.
\end{proof}

\begin{corollary}[PT/GV correspondence]
\label{cor:nr-ptgv}
Assume \emph{(MF)} and suppose that
\eqref{eq:cohomological-exponential} is known on every reduced Chow locus.
If $(\mathrm{NR})_{\chi,\beta}$ holds for all $(\chi,\beta)$, then the
cohomological PT/GV exponential holds on the full Chow varieties.  Its
hypercohomological Euler characteristic is the numerical PT/GV
correspondence.  Its two Lefschetz gradings give the refined product
\eqref{eq:ckk-product}.
\end{corollary}

\begin{proof}
Corollary~\ref{cor:reduced-exponential} gives the cohomological identity.
Taking Euler characteristics gives the numerical identity, and
Proposition~\ref{prop:refined-realization} gives the refined product.
\end{proof}

At the Betti level, Proposition~\ref{prop:refined-realization} is a statement
about the doubly graded Poincar\'e realization.  A compatible lift to pure
monodromic mixed Hodge modules gives the corresponding Hodge or motivic
realization; a rational perverse sheaf alone does not define an identity of
motives.

The shifts are additive.  If
$\boldsymbol\beta=(\beta_1,\ldots,\beta_\ell)$ and
$r_1+\cdots+r_\ell=n-N(\boldsymbol\beta)$, then
\begin{equation}\label{eq:genus-addition}
 g_\beta=\sum_i g_{\beta_i}+N(\boldsymbol\beta)-\ell+1,
 \qquad
 1-g_\beta+n=\sum_i(1-g_{\beta_i}+r_i).
\end{equation}
The negative of the second equality gives additivity of the cohomological
shift in \eqref{eq:normalized-pt-series}.  Hence no extra sign or shift is
placed outside the symmetric power.

For a decomposition $\boldsymbol\beta$ with multiplicities $m_\gamma$ and
$s\ge0$, put $n=N(\boldsymbol\beta)+s$.  Denote by
$\mathcal U_{\boldsymbol\beta,s}$ the Chow-factorization coefficient defined
by
\begin{equation}\label{eq:chow-factorization-term}
 [q^{n+1-g_\beta}]\,
 R\mu_{(m_\gamma),*}
 \left(\bigboxtimes_\gamma
 \operatorname{Sym}^{m_\gamma}\mathcal A_\gamma(q)\right)
 =\mathcal U_{\boldsymbol\beta,s}[g_\beta-1-n].
\end{equation}
For repeated classes the derived symmetric powers in this definition retain
the Koszul descent isomorphisms.

\begin{proposition}\label{prop:exponential-coefficient}
Assume Conjecture~\ref{conj:cohomological-exponential} in the coefficient
$q^{n+1-g_\beta}Q^\beta$.  Then
\begin{equation}\label{eq:exponential-coefficient}
 \bigl(R\Pi_{\beta,n,*}\Phi_{\beta,n}^{\PT}\bigr)^{\mathrm{ss}}
 \cong \mathcal S_{\beta,n}\oplus
 \bigoplus_{\substack{\boldsymbol\beta\ne(\beta)\\
                       N(\boldsymbol\beta)\le n}}
 \mathcal U_{\boldsymbol\beta,n-N(\boldsymbol\beta)}.
\end{equation}
\end{proposition}

\begin{proof}
Write the connected series in class $\gamma$ as
\[
 \mathcal A_\gamma(q)
 =\sum_{r\ge0}q^{1-g_\gamma+r}
  \mathcal S_{\gamma,r}[g_\gamma-1-r].
\]
Fix a multiplicity function $(m_\gamma)$ and temporarily label its
$\ell=\sum_\gamma m_\gamma$ factors by
$\gamma_1,\ldots,\gamma_\ell$.  Choosing the term of length $r_i$ from the
$i$-th factor gives $q$-degree
$\sum_i(1-g_{\gamma_i}+r_i)$.  If
$\sum_i\gamma_i=\beta$, the genus identity
\eqref{eq:genus-addition} rewrites this degree as
\[
 1-g_\beta+N(\boldsymbol\gamma)+\sum_i r_i.
\]
It equals $1-g_\beta+n$ precisely when
$\sum_i r_i=n-N(\boldsymbol\gamma)$.  For the same choice, the sum of the
cohomological shifts is
\[
 \sum_i(g_{\gamma_i}-1-r_i)=g_\beta-1-n.
\]

If a class occurs repeatedly, the symmetric group permutes the labelled
choices of the $r_i$ with the Koszul signs determined by these shifts.
Passing to its invariant direct summand is exactly the derived symmetric
power in \eqref{eq:chow-factorization-term}.  Thus the contribution of the
factorization type $\boldsymbol\gamma$ is
$\mathcal U_{\boldsymbol\gamma,n-N(\boldsymbol\gamma)}[g_\beta-1-n]$.
The one-part type gives
$\mathcal S_{\beta,n}[g_\beta-1-n]$ by
\eqref{eq:normalized-macdonald-complex}.  Extracting the common shift from
the coefficient of Conjecture~\ref{conj:cohomological-exponential} proves
\eqref{eq:exponential-coefficient}.
\end{proof}

For distinct classes, Proposition~\ref{prop:distinct-product} below proves
$\mathcal U_{\boldsymbol\beta,s}\cong
\mathcal T_{\boldsymbol\beta,s}$.  For repeated classes their identification
is not a formal consequence of Lemma~\ref{lem:partition-descent}.  In the
smooth plane range it follows coefficientwise from the exponential conjecture
and the uniqueness argument in Proposition~\ref{prop:plane-supports}; the
first case is verified directly for two lines in
Proposition~\ref{prop:conic-sign}.  The diagonal of
$\operatorname{Sym}^{m}B_e$ is already part of this middle extension, even
though it maps to nonreduced cycles; see
Proposition~\ref{prop:exponential-minimal-extension}.

\section{The first reducible range for plane curves}\label{sec:plane}

Let $S=\Pp^2$, let $H$ be the class of a line, and write

\[
  B_d=|\Ocal_{\Pp^2}(d)|\cong\Pp^{D_d},
  \qquad
  D_d=\frac{d(d+3)}2,
  \qquad
  g_d=\frac{(d-1)(d-2)}2.
\]

For $0\le n\le d+1$, retain the notation $K_{d,n}$ from
Theorem~\ref{thm:intro-plane}.  By
Theorem~\ref{thm:plane-smooth-range},
$K_{d,n}\cong
R\pi_{d,*}^{[n]}\Q_{\mathcal C_d^{[n]}}[D_d+n]$
for the orientation chosen there.  If $d\ge3$,
Corollary~\ref{cor:plane-nonreduced} shows that it has no simple constituent
supported inside the nonreduced locus.  The notation
$\mathcal T_{\lambda,s}$ is the global middle extension defined in
\eqref{eq:boundary-term}.

\subsection{Codimension and node deformations}

Let $\lambda=(d_1,\ldots,d_\ell)$ be a partition of $d$.  Denote by
$\Sigma_\lambda^\circ$ the dense locus of reduced curves
$C=C_1\cup\cdots\cup C_\ell$, where $C_i$ is smooth of degree $d_i$
and all intersections are transverse with no triple points.  The addition map
$B_{d_1}\times\cdots\times B_{d_\ell}\to B_d$ is generically finite on
this locus.  Hence the closure $Z_\lambda$ has
codimension

\begin{equation}\label{eq:plane-partition-codim}
  \codim_{B_d}Z_\lambda
  =D_d-\sum_iD_{d_i}
  =\sum_{i<j}d_i d_j
  =N(\lambda).
\end{equation}

The same number is the number of nodes joining distinct components of a curve
in $\Sigma_\lambda^\circ$.

\begin{lemma}\label{lem:plane-node-versality}
Let $C\in B_d$ be a reduced nodal plane curve and let
$Z\subset\Pp^2$ be any collection of $N\le d+1$ of its nodes.  The
tangent space of $B_d$ maps surjectively to the direct sum of the smoothing
parameters of the nodes in $Z$.
\end{lemma}

\begin{proof}
At a node $p$, a first-order perturbation $f+\varepsilon g$ maps to the
class of $g(p)$ in the one-dimensional $T^1$-space of the node, after a
choice of local coordinates.  Thus the required map is the evaluation map
\[
  H^0(\Pp^2,\Ocal(d))/\langle f\rangle
  \longrightarrow H^0(Z,\Ocal_Z).
\]
The case $N=0$ is trivial.  Otherwise the class of $f$ maps to zero.  For each
node $p_i$, choose one line through each $p_j$ with $j\ne i$ and avoiding
$p_i$.  Their product has degree
$N-1$, vanishes at every $p_j$ with $j\ne i$, and not at $p_i$; multiply by
a general form of degree $d-N+1$.  These forms give a basis of values on
$Z$, so evaluation is surjective.
\end{proof}

If $N(\lambda)\le n\le d+1$, the generic curve in
$\Sigma_\lambda^\circ$ has exactly $N(\lambda)$ nodes, and
Lemma~\ref{lem:plane-node-versality} makes the classifying map transverse
to all relevant node-persistence strata there.
Proposition~\ref{prop:generic-partition} therefore identifies the strict
support term over $B_d^{\red}$ with
$j_d^*\mathcal T_{\lambda,n-N(\lambda)}$, including the finite descent under
component monodromy.  Corollary~\ref{cor:plane-nonreduced} gives its global
middle extension.  This determines the complete summand with strict support
$Z_\lambda$, but does not classify summands supported on proper closed
subvarieties of $Z_\lambda$.

For every nontrivial partition $\lambda$ of $d$, one has
$N(\lambda)\ge d-1$, with equality only for $(d-1,1)$.  Indeed,
$2N(\lambda)=d^2-\sum_i d_i^2$, and among nontrivial partitions the sum of
the squares is largest for $(d-1,1)$.

\begin{proposition}\label{prop:plane-supports}
Let $d\ge3$ and $0\le n\le d+1$.  There is a decomposition
\[
  K_{d,n}
  \cong
  \mathcal S_{d,n}
  \oplus
  \bigoplus_{\substack{\lambda\vdash d,\ \ell(\lambda)>1\\
                        N(\lambda)\le n}}
  \mathcal T_{\lambda,n-N(\lambda)}
  \oplus\mathcal E_{d,n},
\]
where no simple constituent of $\mathcal E_{d,n}$ has full support or support
equal to one of the displayed $Z_\lambda$.  One has
$\mathcal E_{d,n}=0$ for $0\le n\le d-1$.  If
Conjecture~\ref{conj:cohomological-exponential} holds in the coefficient
$q^{n+1-g_d}Q^d$, then $\mathcal E_{d,n}=0$ as well.
\end{proposition}

\begin{proof}
Theorem~\ref{thm:plane-smooth-range} makes
$\mathcal C_d^{[n]}$ smooth, and $\pi_d^{[n]}$ is projective.  The
decomposition theorem therefore makes $K_{d,n}$ semisimple.  Its simple
constituents with support $B_d$ form $\mathcal S_{d,n}$.  For each partition
$\lambda$ with $N(\lambda)\le n$, Lemma~\ref{lem:plane-node-versality} and
Proposition~\ref{prop:generic-partition} determine the complete summand with
strict support $Z_\lambda$ as
$\mathcal T_{\lambda,n-N(\lambda)}$.  Define $\mathcal E_{d,n}$ to be the
direct sum of the remaining simple constituents.  This gives the asserted
decomposition and the stated exclusion of supports.

We prove first that $\mathcal E_{d,n}=0$ for $n\le d-1$.
Corollary~\ref{cor:plane-nonreduced} excludes constituents supported in the
nonreduced locus.  Over the open locus of integral curves, the restricted
family consists of integral locally planar curves and its relative Hilbert
scheme is smooth by Theorem~\ref{thm:plane-smooth-range}.  The support theorem
of \cite[Theorem~1]{MiglioriniShende} therefore applies and excludes proper
supports.  Since the full-support
constituents have already been removed, every support of
$\mathcal E_{d,n}$ must lie in the reducible locus.

Every irreducible component of the reducible locus has codimension at least
$d-1$, and $Z_{d-1,1}$ is the unique component of codimension $d-1$.
Proposition~\ref{prop:support-bound} bounds the codimension of a support by
$n$.  Hence no residual support exists for $n\le d-2$.  When $n=d-1$, a
residual support would have to equal $Z_{d-1,1}$: every proper closed subset
of this irreducible variety has codimension at least $d$ in $B_d$.  But the
complete summand with this support is already
$\mathcal T_{(d-1,1),0}$, so no residual constituent remains.

If the stated exponential coefficient holds,
Proposition~\ref{prop:exponential-coefficient} expresses $K_{d,n}$ as the
one-part term and the factorization terms $\mathcal U_{\lambda,s}$.  The
finite quotient and addition argument in
Proposition~\ref{prop:exponential-minimal-extension}, whose moving-factor
hypothesis holds for plane curves, shows that every nonzero simple
constituent of $\mathcal U_{\lambda,s}$ has strict support exactly
$Z_\lambda$.  Proposition~\ref{prop:generic-partition} identifies its
restriction on the dense reduced stratum with that of
$\mathcal T_{\lambda,s}$.  A semisimple middle extension is uniquely
determined by this generic local system, so
$\mathcal U_{\lambda,s}\cong\mathcal T_{\lambda,s}$.  The exponential side
has no further supports and $\mathcal E_{d,n}=0$.
\end{proof}

For $d\ge4$, the second smallest partition codimension is
$N(d-2,2)=2d-4$.  This follows by transferring degree from a smaller part to
a larger one; any partition with at least three parts has larger value.

\begin{corollary}\label{cor:before-first-boundary}
If $d\ge3$ and $0\le n\le d-2$, then
$K_{d,n}\cong\mathcal S_{d,n}$.
\end{corollary}

\begin{proof}
Every nontrivial partition satisfies $N(\lambda)\ge d-1$, so the direct sum
of partition terms in Proposition~\ref{prop:plane-supports} is empty when
$n\le d-2$.  The same proposition gives $\mathcal E_{d,n}=0$ in this range,
leaving only the full-support summand $\mathcal S_{d,n}$.
\end{proof}

\subsection{The first reducible coefficients}

For $a+b=d$, $a\ne b$, and $m\ge0$, abbreviate
$\mathcal B_{a,b}^{(m)}=\mathcal T_{(a,b),m}$.

On the dense stratum $\Sigma_{a,b}^\circ$, its local systems are the
cohomology sheaves of

\begin{equation}\label{eq:two-component-local-complex}
  \bigoplus_{r=0}^{m}
  R\rho_{a,*}^{[r]}\Q
  \boxtimes
  R\rho_{b,*}^{[m-r]}\Q,
\end{equation}

and $\mathcal B_{a,b}^{(m)}$ is their middle extension with the shifts
$[m-k]$ of \eqref{eq:boundary-term}.  In particular,
$\mathcal B_{a,b}^{(0)}=\IC_{Z_{a,b}}$ after forgetting the Tate twist.

\begin{theorem}\label{thm:plane-first-window}
For $d\ge3$, the first reducible coefficient is
\begin{equation}\label{eq:first-plane-boundary}
 K_{d,d-1}
 \cong
 \mathcal S_{d,d-1}\oplus\mathcal B_{d-1,1}^{(0)}.
\end{equation}
\end{theorem}

\begin{proof}
The partition $(d-1,1)$ is the unique reducible support of codimension at
most $d-1$.  Its excess length is zero, and
Proposition~\ref{prop:plane-supports} has no residual term in this
coefficient.
\end{proof}

\begin{corollary}[First coefficients of the Chow exponential]
\label{cor:plane-exponential-window}
Assume Conjecture~\ref{conj:cohomological-exponential} in the coefficients
under consideration.  For $d\ge5$,
\begin{equation}\label{eq:second-plane-boundary}
 K_{d,d}
 \cong
 \mathcal S_{d,d}\oplus\mathcal B_{d-1,1}^{(1)}.
\end{equation}
For $d\ge6$,
\begin{equation}\label{eq:third-plane-boundary}
 K_{d,d+1}
 \cong
 \mathcal S_{d,d+1}\oplus\mathcal B_{d-1,1}^{(2)}.
\end{equation}
For $(d,n)=(5,6)$,
\begin{equation}\label{eq:quintic-six-boundary}
  K_{5,6}
  \cong
  \mathcal S_{5,6}
  \oplus\mathcal B_{4,1}^{(2)}
  \oplus\mathcal B_{3,2}^{(0)}.
\end{equation}
\end{corollary}

\begin{proof}
The next partition codimension is $N(d-2,2)=2d-4$.  It exceeds $d$ for
$d\ge5$ and exceeds $d+1$ for $d\ge6$, so only $(d-1,1)$ occurs in the
first two formulas.  Its excess lengths are one and two.  In degree five,
$N(3,2)=6$, so the $(3,2)$ term enters with excess zero.  Now apply
Proposition~\ref{prop:plane-supports}.
\end{proof}

\begin{remark}
The exclusions $d\ge5$ and $d\ge6$ avoid additional low-degree
partitions.  For example, $(2,2)$ already contributes to $K_{4,4}$.  Such
cases have the partition terms, together with the possible residual
complex, described by Proposition~\ref{prop:plane-supports}.  The displayed
formulas isolate the range in which only the most unbalanced partition is
expected.
\end{remark}

\section{The Katz--Klemm--Vafa recursion}
\label{sec:kkv-interpretation}

Taking hypercohomological Euler characteristics turns the support terms into
the Katz--Klemm--Vafa recursion.  Write $\chi(B,K)$ for the Euler
characteristic of a constructible complex and set
$p_{d,n}:=\chi(B_d,K_{d,n})=(-1)^{D_d+n}\chi(\mathcal C_d^{[n]})$.
For $0\le n\le d+1$, define the partition contribution by
\begin{equation}\label{eq:kkv-rdn}
  r_{d,n}:={}
  \sum_{\substack{\lambda\vdash d,\ \ell(\lambda)>1\\
                   N(\lambda)\le n}}
  \chi\bigl(B_d,\mathcal T_{\lambda,n-N(\lambda)}\bigr).
\end{equation}

Let $\varepsilon_{d,n}:=\chi(B_d,\mathcal E_{d,n})$ be the Euler
characteristic of the residual complex in
Proposition~\ref{prop:plane-supports}, and put
$a_{d,n}:=\chi(B_d,\mathcal S_{d,n})$.  In the smooth incidence range,
Proposition~\ref{prop:plane-supports} says
$a_{d,n}=p_{d,n}-r_{d,n}-\varepsilon_{d,n}$.

Proposition~\ref{prop:plane-supports} gives $\mathcal E_{d,n}=0$
for $n\le d-1$ by the support argument.  It also vanishes whenever the corresponding
coefficient of Conjecture~\ref{conj:cohomological-exponential} holds.  The
term $\varepsilon_{d,n}$ records what is not determined by the local MSV
calculation without the global exponential conjecture.

The incidence side of this formula is elementary in the smooth range.  Set
$e_n=\chi(\Hilb^n(\Pp^2))$.  G\"ottsche's formula \cite{Goettsche} gives

\begin{equation}\label{eq:goettsche-p2}
  \sum_{n\ge0}e_nq^n=\prod_{m\ge1}(1-q^m)^{-3},
  \qquad
  (e_0,\ldots,e_8)=(1,3,9,22,51,108,221,429,810).
\end{equation}

Since $\mathcal C_d^{[n]}$ is a $\Pp^{D_d-n}$-bundle over
$\Hilb^n(\Pp^2)$,

\begin{equation}\label{eq:incidence-general}
  p_{d,n}=(-1)^{D_d+n}(D_d-n+1)e_n,
  \qquad 0\le n\le d+1.
\end{equation}

Define the numerical full-support GV coefficients $n_d^h$ by the
Macdonald-basis expansion

\begin{equation}\label{eq:kkv-basis-series}
  \sum_{n\ge0}a_{d,n}q^{1-g_d+n}
  =
  \sum_{h=0}^{g_d}n_d^h q^{1-h}(1+q)^{2h-2}.
\end{equation}

Indeed, put
$c_a=\chi(B_d,\IC_{B_d}(\wedge^aV_d))$.  Formula
\eqref{eq:finite-macdonald} gives
\[
 \sum_{n\ge0}a_{d,n}q^{1-g_d+n}
 =\frac{q^{1-g_d}\sum_{a=0}^{2g_d}c_aq^a}{(1+q)^2}.
\]
The symplectic pairing on $V_d$ gives $c_a=c_{2g_d-a}$, and the polynomials
$q^{g_d-h}(1+q)^{2h}$, $0\le h\le g_d$, form a unitriangular basis of the
palindromic polynomials of degree $2g_d$.  This proves existence, uniqueness,
and integrality of the coefficients.  We use generalized binomial
coefficients, with $\binom rm=0$ for $m<0$.

\begin{proposition}[The sheaf-theoretic KKV recursion]
\label{prop:kkv-recursion}
For $d\ge3$ and $0\le n\le d+1$,

\begin{equation}\label{eq:kkv-master}
  (-1)^{D_d+n}\chi(\mathcal C_d^{[n]})
  =
  \sum_{h=0}^{g_d}
  \binom{2h-2}{h-g_d+n}n_d^h
  +r_{d,n}+\varepsilon_{d,n}.
\end{equation}

Equivalently, for
$0\le\delta\le\min\{g_d,d+1\}$, the invariants are recovered
unitriangularly by

\begin{equation}\label{eq:kkv-unitriangular}
  n_d^{g_d-\delta}
  =a_{d,\delta}
   -\sum_{j=0}^{\delta-1}
    \binom{2g_d-2j-2}{\delta-j}n_d^{g_d-j}.
\end{equation}
\end{proposition}

\begin{proof}
For a fixed genus $h$, coefficient extraction gives
\[
 [q^{1-g_d+n}]\,q^{1-h}(1+q)^{2h-2}
 =\binom{2h-2}{h-g_d+n}.
\]
Taking the coefficient of $q^{1-g_d+n}$ in
\eqref{eq:kkv-basis-series} therefore expresses $a_{d,n}$ as the sum in
\eqref{eq:kkv-master}.  Proposition~\ref{prop:plane-supports} gives
$p_{d,n}=a_{d,n}+r_{d,n}+\varepsilon_{d,n}$, while
$p_{d,n}=(-1)^{D_d+n}\chi(\mathcal C_d^{[n]})$ by definition.  Substitution
proves \eqref{eq:kkv-master}.

Set $n=\delta$ and write $h=g_d-j$.  The lower binomial index becomes
$\delta-j$.  It is negative for $j>\delta$, so only
$0\le j\le\delta$ occurs, and for $j=\delta$ the coefficient is
$\binom{2g_d-2\delta-2}{0}=1$.  Solving the resulting triangular equation
for $n_d^{g_d-\delta}$ yields \eqref{eq:kkv-unitriangular}.
\end{proof}

For arbitrary $n\ge0$, set
\[
 p_{d,n}^{\PT}
 =\chi\bigl(B_d,R\Pi_{d,n,*}\Phi_{d,n}^{\PT}\bigr)
\]
and define the factorization contribution
\begin{equation}\label{eq:kkv-all-degree-factorization}
 r_{d,n}^{\ExpChow}
 =\sum_{\substack{\lambda\vdash d,\ \ell(\lambda)>1\\
                   N(\lambda)\le n}}
   \chi\bigl(B_d,\mathcal U_{\lambda,n-N(\lambda)}\bigr).
\end{equation}

\begin{corollary}[KKV recursion]
\label{cor:kkv-all-degree-nr}
Assume \emph{(MF)}, suppose the cohomological PT/GV identity is known on
the reduced Chow loci, and assume $(\mathrm{NR})_{\chi,dH}$ for every
$\chi$ and $d$.  Then, for every $d\ge1$ and $n\ge0$,
\begin{equation}\label{eq:kkv-all-degree}
 p_{d,n}^{\PT}
 =\sum_{h=0}^{g_d}
   \binom{2h-2}{h-g_d+n}n_d^h
  +r_{d,n}^{\ExpChow}.
\end{equation}
For $d\ge3$ and $0\le n\le d+1$, this is
\eqref{eq:kkv-master} with $\varepsilon_{d,n}=0$.  For $n>d+1$,
\eqref{eq:kkv-all-degree} remains valid with the PT vanishing-cycle
coefficient $p_{d,n}^{\PT}$ in place of the shifted-constant Euler
characteristic of the relative Hilbert scheme.
\end{corollary}

\begin{proof}
Corollary~\ref{cor:nr-ptgv} gives the cohomological exponential on the full
Chow varieties.  Proposition~\ref{prop:exponential-coefficient} decomposes
its degree-$(d,n)$ coefficient into the full-support term
$\mathcal S_{d,n}$ and the factorization terms in
\eqref{eq:kkv-all-degree-factorization}.  Taking Euler characteristics and
using $\chi(K)=\chi(K^{\mathrm{ss}})$ and
\eqref{eq:kkv-basis-series} gives
\eqref{eq:kkv-all-degree}.  In the smooth range,
Theorem~\ref{thm:plane-smooth-range} identifies $p_{d,n}^{\PT}$ with
$(-1)^{D_d+n}\chi(\mathcal C_d^{[n]})$, while
Proposition~\ref{prop:plane-supports} identifies the factorization terms
with the summands used in $r_{d,n}$ and removes the residual complex.
\end{proof}

\begin{remark}
The bound $n\le d+1$ belongs to the elementary incidence description, not
to the KKV recursion itself.  Beyond that bound the stable-pair space can
contain non-zero-section components.  If $(\mathrm{Hilb})_{d,n}$ holds,
the reduced restriction of its vanishing-cycle sheaf is computed by the
relative-Hilbert intersection complex; $(\mathrm{NR})_{d,n}$ then makes the
semisimplified proper PT direct image its perverse minimal extension.  The ordinary Euler
characteristic of the singular incidence space is not a substitute for
$p_{d,n}^{\PT}$.
\end{remark}

For example, the first three steps are

\begin{align*}
  n_d^{g_d}&=a_{d,0},\\
  n_d^{g_d-1}&=a_{d,1}-(2g_d-2)a_{d,0},\\
  n_d^{g_d-2}
  &=a_{d,2}-(2g_d-4)a_{d,1}
    +\frac{(2g_d-2)(2g_d-5)}2a_{d,0}.
\end{align*}

These are the equations underlying
\cite[Equations~(4.13)--(4.15)]{KKV}.  Equation~(4.15) there prints
$b(g,0)=0$, whereas Equation~(4.13) forces $b(g,0)=1$; the unitriangular
formula fixes the normalization and the signs.

\subsection{Reducible curves and the product correction}

The subtraction $r_{d,n}$ is the sheaf-theoretic form of the reducible-curve
correction in KKV.  Write
$\mathsf A_d(q):=\sum_{n\ge0}a_{d,n}q^{1-g_d+n}
=\sum_{h=0}^{g_d}n_d^h q^{1-h}(1+q)^{2h-2}$.

Suppose that $\lambda=(d_1,\ldots,d_\ell)$ has pairwise distinct parts, and
let $\mu_\lambda:\prod_iB_{d_i}\to Z_\lambda\hookrightarrow B_d$ be the
factorization morphism.  Multiplication of forms is finite, the source is a
product of irreducible projective spaces and hence normal, and unique
factorization of a general reduced plane curve makes the map generically
one-to-one.  Thus $\mu_\lambda$ is the normalization of its image
$Z_\lambda$.

\begin{proposition}\label{prop:distinct-product}
For every $s\ge0$,

\begin{equation}\label{eq:boundary-complex-product}
  \mathcal T_{\lambda,s}
  \cong
  \bigoplus_{r_1+\cdots+r_\ell=s}
  \mu_{\lambda,*}
  \bigl(
    \mathcal S_{d_1,r_1}\boxtimes\cdots\boxtimes
    \mathcal S_{d_\ell,r_\ell}
  \bigr).
\end{equation}

Consequently,

\begin{equation}\label{eq:kkv-reducible-product}
  \sum_{s\ge0}
  \chi\bigl(B_d,\mathcal T_{\lambda,s}\bigr)
  q^{1-g_d+N(\lambda)+s}
  =\prod_{i=1}^{\ell}\mathsf A_{d_i}(q).
\end{equation}
\end{proposition}

\begin{proof}
Let $U\subset\prod_iB_{d_i}^{\sm}$ be the open set on which the component
curves meet transversely and no three meet at a point.  Because the degrees
$d_i$ are pairwise distinct, the factorization map identifies $U$ with the
component-labelled cover of $\Sigma_\lambda^\circ$; no nontrivial
permutation preserves the degrees.

Over $U$, partial normalization at the nodes joining different components
is the disjoint union of the universal curves of degrees $d_i$.  Its
length-$s$ relative Hilbert scheme is therefore
\[
 \coprod_{r_1+\cdots+r_\ell=s}
 \prod_i\mathcal C_{d_i}^{[r_i]}.
\]
K\"unneth identifies its degree-$k$ direct-image local system with the sum
over $k_1+\cdots+k_\ell=k$ of
$\boxtimes_iR^{k_i}\rho_{d_i,*}^{[r_i]}\Q$.  The shifts in the external
product of the full-support complexes add to
$\sum_i(r_i-k_i)=s-k$.  Hence the restriction of the right side of
\eqref{eq:boundary-complex-product} to the dense stratum is exactly the
local system and shift defining $\mathcal T_{\lambda,s}$.

The map $\mu_\lambda$ is finite and generically one-to-one, and its source
is normal; it is therefore the normalization of $Z_\lambda$.  A finite map
is perverse $t$-exact and small onto its image.  Its pushforward of an IC
complex is consequently the middle extension of the finite-pushforward
local system on the dense stratum.  The preceding generic identification
thus extends uniquely to all of $Z_\lambda$, proving
\eqref{eq:boundary-complex-product}.

Properness of $\mu_\lambda$ and K\"unneth now give
\[
 \chi(B_d,\mathcal T_{\lambda,s})
 =\sum_{r_1+\cdots+r_\ell=s}
   \prod_i\chi(B_{d_i},\mathcal S_{d_i,r_i}).
\]
Finally, \eqref{eq:genus-addition} yields
$1-g_d+N(\lambda)+s=\sum_i(1-g_{d_i}+r_i)$.  Multiplying by the
corresponding power of $q$ and summing over $s$ gives
\eqref{eq:kkv-reducible-product}.
\end{proof}

Here $g_d=\sum_i g_{d_i}+N(\lambda)-\ell+1$.  If
$h_\lambda=\sum_i h_i-\ell+1$, then
$\prod_iq^{1-h_i}(1+q)^{2h_i-2}
=q^{1-h_\lambda}(1+q)^{2h_\lambda-2}$.

Thus, when $h_\lambda\ge0$, the same unitriangular change of basis shows that
the support of type $\lambda$ contributes $\prod_i n_{d_i}^{h_i}$ in genus
$h_\lambda=\sum_i h_i-\ell+1$.  Negative $h_\lambda$ belongs to the
disconnected stable-pair expansion and is not called a connected GV genus.
Since
$p_{d,n}=a_{d,n}+r_{d,n}+\varepsilon_{d,n}$, this product is subtracted from
the value obtained by applying the Macdonald inversion to the total incidence
coefficients.  When the corresponding coefficient of the cohomological
exponential holds, this is exactly the geometric content of the correction
in \cite[Equation~(5.6)]{KKV}.  In general $\varepsilon_{d,n}$ records the
unresolved strict-support terms.

For repeated parts, \eqref{eq:kkv-reducible-product} must be replaced by the
graded equivariant symmetric power.  Dividing an ordinary product by
$|\operatorname{Aut}(\lambda)|$ loses the monodromy in the normal directions.
The sign local system for two lines in
Proposition~\ref{prop:conic-sign} is the first example.  A normalized
excess-zero degree-$e$ factor has total cohomological parity
$D_e+g_e-1\equiv e^2\pmod2$, so Koszul transposition acts by
$(-1)^{e^2}$.  Thus a line is odd whereas a conic is even.
The excess-zero invariant IC descent for two conics has Euler
characteristic $\binom{\chi(\Pp^5)+1}{2}=21$: it is the invariant summand of
the finite pushforward from $\Pp^5\times\Pp^5$.  It is therefore not obtained
by copying the line sign or by dividing an ordinary product by two.

\subsection{Stable-pair and BPS coefficients}

The stable-pair and BPS coefficient bases are related by the triangular
Macdonald transformation.  For the quintic
$(4,1)$-support, the BPS-basis corrections in genera $2,1,0$ are
$(45,-306,693)$; in the stable-pair coefficient basis they become
$(45,-216,738)$, since $-306+2\cdot45=-216$ and $693+45=738$.

The $(3,2)$-support adds
$n_3^1n_2^0=(-10)(-6)=60$ in genus zero.  Hence the proper-support
coefficient displayed by the partition terms at $n=6$ is
$738+60=798$, whereas the connected genus-zero KKV
correction is $693+60=753$.  The difference is exactly the triangular
Macdonald change of basis.

The coefficient $n_d^{g_d-\delta}$ is defined here from the full-support IC
summand; it is not identified with the ordinary Euler characteristic of the
closure of the $\delta$-nodal locus.  At the coefficient $n=d+2$, the
incidence space is singular along curves $2L+C_{d-2}$ and non-zero-section
stable pairs are no longer excluded.  This is the first coefficient beyond
the smooth-incidence theorem and is treated separately in the final
section.

\subsection{Comparison of the direct images}

The three direct images can be compared without assuming that they coincide.
Suppose the reduced Hilbert comparison holds.  Denote the perverse minimal
extension of the reduced direct image by
$K_{d,n}^{\min}:={}^pj_{d,!*}(R\pi_{d,\red,*}^{[n]}
\IC_{\mathcal C_{d,\red}^{[n]}})^{\mathrm{ss}}$.
Write $K_{d,n}^{\Hilb}:=(R\pi_{d,*}^{[n]}
\IC_{\mathcal C_d^{[n]}})^{\mathrm{ss}}$ for the Hilbert-side direct image,
and write $K_{d,n}^{\PT}:=(R\Pi_{d,n,*}
\Phi_{d,n}^{\PT})^{\mathrm{ss}}$ for the stable-pair direct image.

\begin{proposition}\label{prop:hilbert-pt-defects}
There are noncanonical decompositions
\begin{equation}\label{eq:hilbert-pt-defects}
 K_{d,n}^{\Hilb}\cong K_{d,n}^{\min}\oplus
 \mathcal N_{d,n}^{\Hilb},
 \qquad
 K_{d,n}^{\PT}\cong K_{d,n}^{\min}\oplus
 \mathcal N_{d,n}^{\PT},
\end{equation}
where every simple constituent of either $\mathcal N$ is supported in the
nonreduced locus.  Moreover, $(\mathrm{NR})_{d,n}$ holds if and only if
$\mathcal N_{d,n}^{\PT}=0$.
\end{proposition}

\begin{proof}
Open base change and the restriction property of intersection complexes give
\[
 j_d^*K_{d,n}^{\Hilb}
 \cong
 \bigl(R\pi_{d,\red,*}^{[n]}
 \IC_{\mathcal C_{d,\red}^{[n]}}\bigr)^{\mathrm{ss}}.
\]
The decomposition theorem makes the Hilbert-side direct image semisimple.
Split its simple constituents into those whose supports meet
$B_d^{\red}$ and those supported in the complement.  By
Proposition~\ref{prop:minimal-extension-criterion}, the sum of the first
group is the perverse middle extension of the displayed restriction, hence
is $K_{d,n}^{\min}$; the second group is
$\mathcal N_{d,n}^{\Hilb}$.

For the PT complex, the reduced Hilbert comparison and
Proposition~\ref{prop:open-base-change} identify its restriction with the
same reduced direct image.  Apply the identical separation of simple
constituents to its semisimplification.  This gives the second decomposition
and identifies the complementary summand as $\mathcal N_{d,n}^{\PT}$.
Definition~\ref{def:nonreduced-support} says precisely that this summand is
zero.
\end{proof}

The possible distinction is therefore between a nonzero
$\mathcal N_{d,n}^{\Hilb}$ and the conjectural vanishing of
$\mathcal N_{d,n}^{\PT}$.

Katz--Klemm--Vafa take the ordinary Euler characteristic of the singular
incidence space.  After the common perverse normalization, the corresponding
object is the shifted constant sheaf rather than its intersection complex.
If
$K_{d,n}^{\mathrm{const}}:=R\pi_{d,*}^{[n]}
\Q_{\mathcal C_d^{[n]}}[D_d+n]$, then in the
Grothendieck group its deviation from the minimal PT candidate splits as
\begin{equation}\label{eq:constant-ic-pt-defect}
 [K_{d,n}^{\mathrm{const}}]-[K_{d,n}^{\min}]
 =\bigl([K_{d,n}^{\mathrm{const}}]-[K_{d,n}^{\Hilb}]\bigr)
  +\bigl([K_{d,n}^{\Hilb}]-[K_{d,n}^{\min}]\bigr).
\end{equation}
The first term measures the singularity of the incidence space; the second
measures strict supports of its IC direct image on the nonreduced locus.
Numerical Euler characteristics cannot separate them.  A local IC and
vanishing-cycle calculation is therefore still required to prove that the
entire observed correction is the second term.

\section{Calculations on del Pezzo surfaces}
\label{sec:examples}

For an effective class $\beta$ with a smooth connected member, write
\[
 K_{\beta,n}=R\Pi_{\beta,n,*}\Phi_{\beta,n}^{\PT},\qquad
 p_{\beta,n}=\chi(B_\beta,K_{\beta,n}),\qquad
 a_{\beta,n}=\chi(B_\beta,\mathcal S_{\beta,n}).
\]
The full-support Macdonald coefficients $n_\beta^h$ are defined by
\begin{equation}\label{eq:del-pezzo-macdonald-basis}
 \sum_{n\ge0}a_{\beta,n}q^{1-g_\beta+n}
 =\sum_{h=0}^{g_\beta}n_\beta^h q^{1-h}(1+q)^{2h-2}.
\end{equation}

\subsection{\texorpdfstring{Bidegree $(2,3)$ on
$\Pp^1\times\Pp^1$}{Bidegree (2,3) on P1 x P1}}

Let \(S=\Pp^1\times\Pp^1\), and write a divisor class as \((a,b)\),
with \(1\le a\le b\).  Then
\[
 D_{a,b}=ab+a+b,\qquad g_{a,b}=(a-1)(b-1).
\]
The threshold of \eqref{eq:del-pezzo-threshold} is
\begin{equation}\label{eq:p1p1-threshold}
 \tau_S(a,b)=a+2.
\end{equation}
Indeed, for an ordered decomposition
\((a,b)=\sum_i(a_i,b_i)\), the expression minimized in
\eqref{eq:del-pezzo-threshold} is
\[
 \sum_{i<j}(a_i b_j+a_j b_i)
 +2\sum_{i\ge1}i(a_i+b_i).
\]
Combining all positive-index layers into the first layer can only decrease
this expression: their mutual intersections are removed and their weights
decrease.  Put that layer equal to \((A,B)\).  The resulting value is
\[
 E(A,B)=aB+bA-2AB+2A+2B.
\]

If \(A=0\), this is \(B(a+2)\); if \(B=0\), it is \(A(b+2)\).
When \(A,B>0\), regard it as a linear function of \(A\).  If its
coefficient \(b-2B+2\) is nonnegative, its minimum occurs at \(A=1\)
and is at least \(b+2+a\).  If the coefficient is negative, its minimum
occurs at \(A=a\) and is
\(a(b-B)+2a+2B\ge a+2\).  Equality is attained by the decomposition
\((a,b)=(a,b-1)+(0,1)\), proving \eqref{eq:p1p1-threshold}.

We also need the following elementary jet-separation fact.  If
\(a,b\ge n-1\), then
\begin{equation}\label{eq:p1p1-jet-separation}
 H^1(S,I_Z(a,b))=0
 \quad\text{for every length-$n$ subscheme }Z\subset S.
\end{equation}

For $n=0$ the assertion is immediate.  For $n>0$, choose a fibre
\(F\in|\Ocal(1,0)|\) through a point of \(Z\), put
\(m=\operatorname{length}(Z\cap F)\), and use the residual sequence
\[
 0\to I_{\operatorname{Res}_F Z}(a-1,b)
 \to I_Z(a,b)
 \to I_{Z\cap F,F}(b)\to0.
\]
Here $m\ge1$, and the display is the standard residual sequence for the
Cartier divisor $F$.  The residual scheme has length \(n-m\).  Induction
on \(n\) applies to
the first term because \(a-1,b\ge n-m-1\), while the last term has no
first cohomology because \(b-m\ge-1\).  This proves
\eqref{eq:p1p1-jet-separation}.  Theorem
\ref{thm:del-pezzo-incidence-range} therefore identifies the PT space with
the smooth incidence space for \(0\le n\le a+1\).

Take \(\beta=(2,3)\).  Here \(D_\beta=11\) and \(g_\beta=2\).
G\"ottsche's formula gives Euler characteristics \(1,4,14\) for the first
three Hilbert schemes of points on \(S\).  Hence the incidence coefficients
are
\begin{equation}\label{eq:p1p1-23-incidence}
 (p_{\beta,0},p_{\beta,1},p_{\beta,2})=(-12,44,-140).
\end{equation}

For a decomposition \((a,b)=(a_1,b_1)+(a_2,b_2)\), the addition locus has
codimension \(a_1b_2+a_2b_1\), the number of intersection nodes of the
generic components.  In class \((2,3)\), the unique codimension-two locus
is
\[
 Z_{(0,1),(2,2)}:\quad (2,3)=(0,1)+(2,2).
\]
All other reducible loci have codimension at least three.  The possible
doubled classes are \((0,1)\), \((1,0)\), and \((1,1)\); their addition
loci have codimensions \(5,7,7\), respectively.  Thus no nonreduced locus
can support a summand in relative dimension two.

At the generic point of \(Z_{(0,1),(2,2)}\), the two components meet in
two transverse nodes.  Surjectivity onto their smoothing parameters is the
length-two case of \eqref{eq:p1p1-jet-separation}.  Proposition
\ref{prop:generic-partition} gives the complete proper-support summand.
The integral-curve support theorem and Proposition~\ref{prop:support-bound}
exclude every other support, so
\[
 K_{\beta,2}\cong\mathcal S_{\beta,2}
 \oplus\IC_{Z_{(0,1),(2,2)}}.
\]

The addition map from
\(B_{0,1}\times B_{2,2}\cong\Pp^1\times\Pp^8\) is finite and generically
one-to-one, hence is the normalization; being finite, it is small.  With the
nine-dimensional IC shift its Euler characteristic is \(-18\).  The
full-support coefficients are consequently \((-12,44,-122)\).  Writing
them in the genus-two Macdonald basis gives
\[
 -12=n_{(2,3)}^2,\qquad
 44=n_{(2,3)}^1+2n_{(2,3)}^2,\qquad
 -122=n_{(2,3)}^0+n_{(2,3)}^2.
\]
Therefore the entirely determined row is
\begin{equation}\label{eq:p1p1-23-row}
 \boxed{(n_{(2,3)}^0,n_{(2,3)}^1,n_{(2,3)}^2)=(-110,68,-12).}
\end{equation}
It agrees with the local \(\Pp^1\times\Pp^1\) calculation of
\cite{CKK}.

\subsection{Full-support coefficients in anticanonical degree}

Let
\(S_r=\operatorname{Bl}_{p_1,\ldots,p_r}\Pp^2\), where the points are
general, \(5\le r\le8\), and put \(\kappa=K_{S_r}^2=9-r\).  We use
\(E_r\) as a shorthand for this del Pezzo surface.  For
\(\beta=aH-\sum_i b_iE_i\),
\[
 -K_{S_r}\cdot\beta=3a-\sum_i b_i,\qquad
 g_\beta=1+\frac{a(a-3)-\sum_i b_i(b_i-1)}2.
\]
Define the anticanonical-degree aggregate of the full-support Macdonald
coefficients by
\[
 N_d^h(S_r)=
 \sum_{\substack{-K_{S_r}\cdot\beta=d\\
                  |\beta|\text{ has a smooth connected member}}}
 n_\beta^h.
\]

Suppose \(-K_{S_r}\cdot\beta=\kappa\).  The Hodge index theorem gives
\(\beta^2\le\kappa\), hence \(g_\beta\le1\); equality forces
\(\beta=-K_{S_r}\).  A rational class has \(\beta^2=\kappa-2\), so its
plane coefficients satisfy
\begin{equation}\label{eq:exceptional-diophantine}
 \sum_i b_i=3a-\kappa,\qquad
 \sum_i b_i^2=a^2+2-\kappa.
\end{equation}
Solving these equations, and including the exceptional curves when
\(r=8\), gives the following list.  The bracketed number is the number of
permutations of the multiplicities; these counts also appear in the del
Pezzo tables of \cite{KKV}.

\begingroup
\footnotesize
\begin{equation}\label{eq:exceptional-rational-classes}
\begin{array}{c|l|r}
r&\text{types }(a;b_1,\ldots,b_r)&R_\kappa\\ \hline
5&(2;1^2,0^3)[10],\ (3;2,1^3,0)[20],\
  (4;2^3,1^2)[10]&40\\
6&(1;0^6)[1],\ (2;1^3,0^3)[20],\ (3;2,1^4,0)[30],\
  (4;2^3,1^3)[20],\ (5;2^6)[1]&72\\
7&(1;1,0^6)[7],\ (2;1^4,0^3)[35],\ (3;2,1^5,0)[42],\
  (4;2^3,1^4)[35],\ (5;2^6,1)[7]&126\\
8&E_i[8],\ (1;1^2,0^6)[28],\ (2;1^5,0^3)[56],\
  (3;2,1^6,0)[56],\ (4;2^3,1^5)[56],\
  (5;2^6,1^2)[28],\ (6;3,2^7)[8]&240.
\end{array}
\end{equation}
\endgroup

For completeness, the two equalities in
\eqref{eq:exceptional-diophantine} and Cauchy--Schwarz give
\[
 (3a-\kappa)^2\le r(a^2+2-\kappa).
\]
For $r=5,6,7,8$, respectively, this leaves
$2\le a\le4$, $1\le a\le5$, $1\le a\le5$, and $0\le a\le7$.
Substitution in the two equalities gives exactly the displayed multiplicity
patterns; for $r=8$, the case $a=0$ gives the exceptional curves.  The
entries are effective for general points, for example by
the standard quadratic Cremona reductions and class tables in
\cite[Appendix~A]{KKV}.

We first check that the coefficients used below lie on the zero section.
Let $\gamma$ be a proper nonzero effective class with
$d_\gamma=-K_{S_r}\cdot\gamma<\kappa$.  Hodge index gives
$\gamma^2\le d_\gamma^2/\kappa<d_\gamma$; adjunction therefore gives
$p_a(\gamma)\le0$.  For an ordered layer decomposition
$-K_{S_r}=\gamma_0+\cdots+\gamma_{\ell-1}$, the genus identity in the proof
of Lemma~\ref{lem:del-pezzo-zero-section-threshold} gives
\[
 \sum_{i<j}\gamma_i\cdot\gamma_j
 =\ell-\sum_i p_a(\gamma_i)\ge\ell,
 \qquad
 \sum_{i=1}^{\ell-1}i d_{\gamma_i}
 \ge\frac{\ell(\ell-1)}2.
\]
Thus $\tau_{S_r}(-K_{S_r})\ge3$ for $r=5,6,7$.  When $r=8$,
ampleness of $-K_{S_8}$ excludes a decomposition because the positive
integral anticanonical degrees of two summands cannot add to one; hence the
threshold is infinite.

For \(r=5,6,7\), the anticanonical system is the basepoint-free space
\(|-K_{S_r}|\cong\Pp^\kappa\).  Its one-point incidence variety is a
\(\Pp^{\kappa-1}\)-bundle over \(S_r\), and
\(\chi(S_r)=3+r=12-\kappa\).  Thus
\[
 p_{-K,0}=(-1)^\kappa(\kappa+1),\qquad
 p_{-K,1}=(-1)^{\kappa+1}\kappa(12-\kappa).
\]

Every anticanonical divisor on these general blow-ups is reduced.  In the
plane model it is the proper transform of a cubic through at least five
general points.  A nonreduced cubic is $3L$ or $2L+M$ and cannot contain
five general points.  If an exceptional curve $E_i$ occurred with
multiplicity at least two, the residual plane cubic would have multiplicity
at least three at $p_i$ and pass through the other $r-1$ points.  These
conditions have expected projective dimension
$9-6-(r-1)=4-r<0$, and they are independent for general points.  Thus this
case is also impossible.

At a general point of a reducible anticanonical stratum, write the divisor as
$A_1+\cdots+A_\ell$, with the $A_i$ smooth irreducible components.  The
preceding Hodge-index argument gives $p_a(A_i)\le0$, while smooth
connectedness gives $p_a(A_i)\ge0$; hence every $A_i$ is rational.
Adjunction gives
$A_i^2=(-K\cdot A_i)-2$, and the sequence
$0\to\Ocal_{S_r}\to\Ocal_{S_r}(A_i)\to\Ocal_{A_i}(A_i)\to0$
then gives $\dim|A_i|=(-K\cdot A_i)-1$.  Since the anticanonical degrees
sum to $\kappa$, the corresponding addition stratum has codimension
\[
 \kappa-\sum_i\bigl((-K\cdot A_i)-1\bigr)=\ell.
\]

The nonintegral locus therefore has codimension at least two.  The support
bound for the one-point family permits only codimension one, and the
integral-curve support theorem excludes a proper support on the remaining
open set.  Hence the two displayed numbers are the genus-one and genus-zero
full-support coefficients of $-K$.

For a rational class $\beta$ in
\eqref{eq:exceptional-rational-classes}, the same genus identity applied to
an ordered layer decomposition gives
\[
 \sum_{i<j}\gamma_i\cdot\gamma_j\ge\ell-1,
 \qquad
 \sum_{i=1}^{\ell-1}i(-K\cdot\gamma_i)
 \ge\frac{\ell(\ell-1)}2.
\]
Thus $\tau_{S_r}(\beta)\ge2$.  For $r\le7$ the listed classes are nef, so
Kodaira vanishing and Riemann--Roch give
$|\beta|\cong\Pp^{\kappa-1}$.  Their genus-zero coefficient is therefore
$(-1)^{\kappa-1}\kappa$.  For $r=8$ they are the 240 $(-1)$-curves and
their linear systems are points.

It remains to treat the base-pointed anticanonical pencil on $S_8$.  The
incidence space $\mathcal C_{-K}^{[1]}$ is the blow-up of $S_8$ at the
unique base point and is smooth, with Euler characteristic $12$.  Since
$\tau_{S_8}(-K)=\infty$, every geometric PT pair lies on the zero section.
The closed zero-Higgs, tangent-space, and local-ring argument of
Theorem~\ref{thm:del-pezzo-incidence-range} then proves
\[
 P_1(\Tot(K_{S_8}),-K_{S_8})\cong\mathcal C_{-K}^{[1]}.
\]

For its smooth orientation the vanishing-cycle sheaf is the intersection
complex.  Every anticanonical member is integral, so the integral-curve
support theorem shows that $p_{-K,1}=12$ is a full-support coefficient;
also $p_{-K,0}=-2$.  Combining the elliptic and rational classes gives
\begin{equation}\label{eq:exceptional-aggregate-table}
\boxed{
\begin{array}{c|c|rr}
S&d&N_d^1&N_d^0\\ \hline
E_5&4&5&-192\\
E_6&3&-4&243\\
E_7&2&3&-272\\
E_8&1&-2&252.
\end{array}}
\end{equation}

These full-support Macdonald coefficients agree with the corresponding KKV
aggregate coefficients.  Their
identification with a sheaf-moduli definition of GV invariants is not used
in the calculation.

\subsection{The sign local system for reducible conics}

Let $d=2$.  The dense reduced singular locus
$\Sigma_{1,1}^\circ\subset B_2$ parametrizes unions
$C=L_1\cup L_2$ of two distinct lines.  Its component-labelled cover is

\[
  \varpi:\widetilde\Sigma_{1,1}^\circ
  =\bigl((\Pp^2)^\vee\times(\Pp^2)^\vee\bigr)\setminus\Delta
  \longrightarrow \Sigma_{1,1}^\circ,
\]

with deck group $\mathfrak S_2$.  Denote by
$\mathcal L_{\mathrm{sgn}}$ the rank-one local system associated with the
sign representation.

\begin{proposition}\label{prop:conic-sign}
There is a decomposition
\[
 K_{2,1}\cong\mathcal S_{2,1}\oplus
 \IC_{\Sigma_{1,1}}(\mathcal L_{\mathrm{sgn}}).
\]
Thus the proper-support local system is the sign local system, rather than
the constant local system on the singular-conic divisor.
\end{proposition}

\begin{proof}
Theorem~\ref{thm:plane-smooth-range},
Lemma~\ref{lem:plane-node-versality}, and
Proposition~\ref{prop:generic-partition} produce the rank-one proper-support
local system on the dense nodal conic locus.  The nonreduced locus has
codimension three, greater than the relative dimension one, so
Proposition~\ref{prop:support-bound} excludes a further strict-support term
there and the local system middle-extends to the singular-conic divisor.
Over the labelled cover,
$H^2(L_1\cup L_2,\Q)=\Q[L_1]\oplus\Q[L_2]$.
In a one-parameter smoothing of $L_1\cup L_2$, the fundamental class of the
smooth conic specializes to $[L_1]+[L_2]$.  This invariant line is therefore
the specialization of the full-support local system.  The complementary
line, which is the local MSV proper-support term, is generated by
$[L_1]-[L_2]$.  The deck involution exchanges $L_1$ and $L_2$, so it acts by
$-1$ on this line.  The proper-support local system consequently descends as
$\mathcal L_{\mathrm{sgn}}$.

The integral-curve support theorem excludes
another proper support over the irreducible locus.  The support bound allows
only codimension at most one, while the double-line locus has codimension
three.  Hence the displayed middle extension is the only proper-support
summand.
\end{proof}

The same sign is forced by the Chow-symmetric description.  The degree-one GV
complex is placed in an odd cohomological shift.  The Koszul sign in its
symmetric square makes the transposition act by the sign character; finite
descent alone would only give the decomposition
$\varpi_*\Q=\Q\oplus\mathcal L_{\mathrm{sgn}}$ and would not determine which
isotypical term occurs.

\subsection{The quartic row}

For $d=4$, one has $D_4=14$ and $g_4=3$.  The Euler characteristics
of the first four Hilbert schemes of points on $\Pp^2$ are
$1,3,9,22$.  Formula \eqref{eq:incidence-general} therefore gives

\[
\begin{array}{c|rrrr}
n&0&1&2&3\\ \hline
p_{4,n}&15&-42&117&-264.
\end{array}
\]

Every nontrivial partition of four has codimension at least three, and the
unique codimension-three support is $Z_{3,1}$.  It first occurs at $n=3$,
with excess length zero.  Theorem~\ref{thm:plane-first-window} gives its
complete summand, whose Euler characteristic is
\[
 r_{4,3}=p_{3,0}p_{1,0}=(-10)\cdot3=-30.
\]
There is no residual support in this range by
Proposition~\ref{prop:plane-supports}.  Hence the full-support coefficients
are $(15,-42,117,-234)$.

In genus three, the four Macdonald equations are
\[
 \begin{aligned}
 15&=n_4^3,\\
 -42&=n_4^2+4n_4^3,\\
 117&=n_4^1+2n_4^2+6n_4^3,\\
 -234&=n_4^0+n_4^2+4n_4^3.
 \end{aligned}
\]
Solving them in this order yields

\begin{equation}\label{eq:quartic-row}
\boxed{
\begin{array}{c|rrrr}
h&0&1&2&3\\ \hline
n_4^h&-192&231&-102&15.
\end{array}}
\end{equation}

This is the complete quartic row of \cite[Table~4]{KKV}.  The partition
$(2,2)$ first appears at $n=4$, so it does not enter the reconstruction of
this row.

\subsection{The quintic row}

For $d=5$, one has $B_5\cong\Pp^{20}$ and $g_5=6$.  For
$0\le n\le6$, the relative Hilbert scheme $\mathcal C_5^{[n]}$ is a
projective bundle of fiber dimension $20-n$ over $\Hilb^n(\Pp^2)$.  Thus
$\chi(\mathcal C_5^{[n]})=(21-n)\chi(\Hilb^n(\Pp^2))$.

The required Hilbert-scheme Euler characteristics are listed in
\eqref{eq:goettsche-p2}.

The smooth stable-pair space has dimension $20+n$, so its Behrend function
is $(-1)^{20+n}$.  In the numerical notation $p_{d,n}$ fixed above, we
obtain

\begin{equation}\label{eq:quintic-PT-table}
\begin{array}{c|rrrrrrr}
n&0&1&2&3&4&5&6\\ \hline
p_{5,n}&21&-60&171&-396&867&-1728&3315.
\end{array}
\end{equation}

Let $n_5^h$ be the unrefined degree-five genus-$h$ BPS number.  The Euler
characteristic of the full-support term is

\begin{equation}\label{eq:quintic-smooth-general}
  \chi(B_5,\mathcal S_{5,n})
  =\sum_{h=\max(0,6-n)}^6
   \binom{2h-2}{h-6+n}n_5^h.
\end{equation}

Here generalized binomial coefficients are defined by
$(1+x)^r=\sum_{k\ge0}\binom rkx^k$ for every integer $r$.  In particular,
$\binom{-2}{0}=1$, which supplies the genus-zero term when $n=6$.
Writing out \eqref{eq:quintic-smooth-general} gives

\begin{equation}\label{eq:quintic-smooth-table}
\begin{array}{c|l}
n&\chi(B_5,\mathcal S_{5,n})\\ \hline
0&n_5^6\\
1&n_5^5+10n_5^6\\
2&n_5^4+8n_5^5+45n_5^6\\
3&n_5^3+6n_5^4+28n_5^5+120n_5^6\\
4&n_5^2+4n_5^3+15n_5^4+56n_5^5+210n_5^6\\
5&n_5^1+2n_5^2+6n_5^3+20n_5^4+70n_5^5+252n_5^6\\
6&n_5^0+n_5^2+4n_5^3+15n_5^4+56n_5^5+210n_5^6.
\end{array}
\end{equation}

The partition codimensions in degree five are

\[
\begin{array}{c|rrrrrr}
\lambda&(4,1)&(3,2)&(3,1,1)&(2,2,1)&(2,1,1,1)&(1,1,1,1,1)\\ \hline
N(\lambda)&4&6&7&8&9&10.
\end{array}
\]

Thus Proposition~\ref{prop:plane-supports} gives no partition contribution
for $n\le3$, the $(4,1)$-term for $n=4,5$, and the $(4,1)$- and
$(3,2)$-terms for $n=6$.  Possible residual terms are kept separately in
$\varepsilon_{5,n}$ below.

The Euler characteristics of these terms can be obtained directly from the
lower-degree coefficients.  In every entry used below the length is below the
first reducible threshold of the component degree (or the component is a
line), so its total PT coefficient equals its full-support coefficient.  The
required values are

\[
\begin{array}{c|rrr}
\text{degree }1&p_{1,0}&p_{1,1}&p_{1,2}\\ \hline
&3&-6&9
\end{array}
\qquad
\begin{array}{c|rrr}
\text{degree }4&p_{4,0}&p_{4,1}&p_{4,2}\\ \hline
&15&-42&117.
\end{array}
\]

For the $(4,1)$-support, the partial-normalization convolution gives

\begin{align*}
 r_{5,4}^{(4,1)}&=15\cdot3=45,\\
 r_{5,5}^{(4,1)}&=15(-6)+(-42)3=-216,\\
 r_{5,6}^{(4,1)}&=15\cdot9+(-42)(-6)+117\cdot3=738.
\end{align*}

At $n=6$, the $(3,2)$-support contributes
$r_{5,6}^{(3,2)}=p_{3,0}p_{2,0}=(-10)(-6)=60$.  Consequently,

\begin{equation}\label{eq:quintic-boundary-table}
\begin{array}{c|rrrrrrr}
n&0&1&2&3&4&5&6\\ \hline
r_{5,n}&0&0&0&0&45&-216&798.
\end{array}
\end{equation}

For example, the first value is also immediate geometrically.  The
normalization of $Z_{4,1}$ is $B_4\times B_1$, so
$\chi(\IC_{Z_{4,1}})=\chi(B_4)\chi(B_1)=15\cdot3=45$.

Equations \eqref{eq:quintic-PT-table},
\eqref{eq:quintic-smooth-table}, and
\eqref{eq:quintic-boundary-table} give the triangular system

\begin{equation}\label{eq:quintic-triangular-system}
  p_{5,n}
  =\chi(B_5,\mathcal S_{5,n})+r_{5,n}+\varepsilon_{5,n},
  \qquad 0\le n\le6.
\end{equation}

Here the support theorem gives $\mathcal E_{5,n}=0$ for $n\le4$.  The first five
equations give $n_5^6=21$, $n_5^5=-270$, $n_5^4=1386$,
$n_5^3=-3672$, and $n_5^2=5430$.

The last two equations give
$n_5^1=-4452-\varepsilon_{5,5}$ and
$n_5^0=1695-\varepsilon_{5,6}$.

Under Conjecture~\ref{conj:cohomological-exponential} through the degree-five
coefficients $n=5,6$, the residual complexes vanish.  Therefore

\begin{equation}\label{eq:quintic-final-row}
\boxed{
\begin{array}{c|rrrrrrr}
h&0&1&2&3&4&5&6\\ \hline
n_5^h&1695&-4452&5430&-3672&1386&-270&21.
\end{array}}
\end{equation}

This is the degree-five local-$\Pp^2$ row of Katz--Klemm--Vafa
\cite{KKV}.  The entries from genus $2$ through genus $6$ follow from the
proved support range
in the present argument.  The genus-one and genus-zero entries use the stated
coefficients of the Chow exponential.  The local MSV calculation identifies
the relevant partition summands on their generic supports.

\section{Degree two with four points}
\label{sec:degree-two-correction}

Set $X=\Tot(K_{\Pp^2})$, $M=P_5(X,2)$, and
$Y=\mathcal C_2^{[4]}$.  Let
$N_2=\{2L:L\subset\Pp^2\text{ a line}\}\subset B_2$, the Veronese surface
of double lines.  Since $g_2=0$, the superscript $[4]$ records the length of
$Z$ on the zero-section component.  This is the first coefficient outside
Theorem~\ref{thm:plane-smooth-range}, and hence the first at which there are
two different moduli schemes to compare.

The zero-section estimate in \cite[Proposition~2]{CKK} was printed with the
range extending to $n=d+2$.  Lemma~\ref{lem:zero-section-threshold} gives
the corrected bound $n\le d+1$, and equality at $n=d+2$ is realized by
vertical ribbons.  The pair $(d,n)=(2,4)$ is the smallest case in which this
additional component appears.

The results of this section give the first scheme-theoretic
comparison beyond the zero-section range between a local-$\Pp^2$
stable-pair space and its surface Hilbert-flag model.  The latter is the
relative Hilbert, or incidence, scheme of pairs $(Z,C)$ with $Z\subset C$.
The calculation proves $M=Y\cup_{\mathfrak B}R$ scheme-theoretically: $Y$
is the zero-section incidence component, $R$ compactifies the non-planar
Ferrand ribbons, and their common flag space $\mathfrak B$ parametrizes four
points on a line.  Along $\mathfrak B$, the component $Y$ has a transverse
threefold node, whereas the completed classical germ of $M$ is
$\mathbb A^6\times\operatorname{Crit}(c(x_1x_2+x_3x_4))$.  Thus the
calculation determines the incidence of the two components and the formal
scheme at their
intersection, not only their closed points or Euler characteristics.  This
is why the relative Hilbert scheme alone does not recover the degree-$(2,4)$ PT
invariant.

Let $H=\Hilb^4(\Pp^2)$ and let $\mathfrak B$ be the space of pairs
$(L,Z)$ with $L\subset\Pp^2$ a line and $Z\subset L$ of length four.
The line containing $Z$ is unique, so $\mathfrak B\to(\Pp^2)^\vee$ is a
$\Pp^4$-bundle.  Its map to $H$ is a closed immersion: an infinitesimal
motion of $L$ that induces the trivial deformation of $Z$ is a section of
$N_{L/\Pp^2}=\Ocal_L(1)$ vanishing on the length-four divisor $Z$, hence is
zero; the proper, injective, unramified map is therefore a closed immersion.
Denote by
$\mathfrak S\subset Y$ the section
$(L,Z)\mapsto(Z,2L)$.

\begin{proposition}[The incidence singularity]\label{prop:conic-correction-node}
The singular locus of $Y$ is $\mathfrak S\cong\mathfrak B$.  For every
$b\in\mathfrak S$, the analytic germ of the pair $(Y,\mathfrak S)$ at $b$
is
\[
 (Y,\mathfrak S)\cong
 \bigl(\mathbb A^6\times
       \{x_1x_2+x_3x_4=0\},\ \mathbb A^6\times\{0\}\bigr).
\]
Thus, at $b\in\mathfrak S$, the stalk of $\IC_Y$ has one copy of $\Q$ in
degrees $-9$ and $-7$, whereas $\Q_Y[9]$ has only the copy in degree
$-9$.  Moreover,
\[
 \chi(Y)=117,\qquad
 \chi(Y,\Q_Y[9])=-117,\qquad
 \chi(Y,\IC_Y)=-132.
\]
\end{proposition}

\begin{proof}
Let $E$ be the rank-four tautological bundle on $H$.  The incidence space is
the projectivized kernel of the evaluation map
$H^0(\Pp^2,\Ocal(2))\otimes\Ocal_H\to E$.  This map has rank four off
$\mathfrak B$ and rank three on $\mathfrak B$.  At $Z\subset L$, its kernel
is $\ell H^0(\Pp^2,\Ocal(1))$, where $\ell$ is an equation of $L$, and its
cokernel is one-dimensional.  Indeed, a length-four scheme fails to impose
four independent conditions on conics precisely when its degree-two
Hilbert function is at most three.  Hilbert-function growth then forces its
degree-one Hilbert function to be at most two, which is equivalent to
containment in a line.  On a line the rank is exactly three.

The normal derivative of the evaluation map is the multiplication pairing
\[
 H^1(L,\Ocal_L(-3))\otimes H^0(L,\Ocal_L(1))
       \longrightarrow H^1(L,\Ocal_L(-2)).
\]
Indeed, normal deformations of $Z$ away from a fixed line form
$H^0(Z,\Ocal_Z(1))$, while motions of the containing line give the image of
$H^0(L,\Ocal_L(1))$.  The sequence
$0\to\Ocal_L(-3)\to\Ocal_L(1)\to\Ocal_Z(1)\to0$ therefore identifies their
quotient, and hence $N_{\mathfrak B/H,b}$, with
$H^1(L,\Ocal_L(-3))$.  The last factor in the display is the cokernel of
evaluation in degree two.  Serre duality
shows that the pairing is perfect.

The remaining kernel direction,
generated by $\ell^2$, has zero normal derivative.  After analytic changes
of bases, write the final row of the evaluation map as $(a,b,c)$.  The
perfect derivative makes $a,b$ regular parameters normal to
$\mathfrak B$.  The corank-one degeneracy ideal is $(a,b,c)$ and has
underlying locus $\mathfrak B$; since $\mathfrak B$ is smooth of
codimension two, $a,b$ generate its ideal scheme-theoretically and
$c\in(a,b)$.  A final column operation gives the exact row $(u,v,0)$.  On
the affine chart around
$[\ell^2]$ in the exceptional $\Pp^2$, the incidence equation is
$ux+vy=0$.  This is the threefold ordinary double point, with the six
coordinates along $\mathfrak B$ as smooth parameters.  The same normal
calculation applies to every length-four divisor on $L$.

For four distinct points this normal form can be seen directly.  Take
$L=(y=0)$ with affine coordinate $t$ and points $t_i$.  A nearby conic can
be written $q=(y+\ell(t))^2+r(t)$ with $\deg\ell\le1$ and $\deg r\le2$.
Putting $z_i=y_i+\ell(t_i)$ and eliminating the three coefficients of $r$
leaves the single equation
\[
 \sum_{i=1}^4\frac{z_i^2}{\prod_{j\ne i}(t_i-t_j)}=0.
\]
All four coefficients are nonzero, so this is a rank-four quadric.  The
Serre-duality argument above is its coordinate-free extension across
collisions of the four points.

A small resolution of the threefold node has exceptional fibre $\Pp^1$.
The small-map description of its intersection complex gives the two stated
stalk groups.  Finally, $Y\to H$ has fibre $\Pp^1$ off $\mathfrak B$ and
$\Pp^2$ on $\mathfrak B$.  Since
$\chi(H)=51$ and $\chi(\mathfrak B)=3\cdot5=15$,
\[
 \chi(Y)=2(51-15)+3\cdot15=117.
\]
The IC stalk Euler function is $-1$ off $\mathfrak S$ and $-2$ on
$\mathfrak S$.  Integrating this function gives
$-117-15=-132$.
\end{proof}

The fibre over a fixed double line also admits an elementary stratification.
Write
$\Ocal_{2L}=\Ocal_L\oplus\epsilon\Ocal_L(-1)$ with $\epsilon^2=0$.

\begin{lemma}[Four points on a double line]\label{lem:four-points-double-line}
The Hilbert scheme $\Hilb^4(2L)$ has three locally closed strata
$\mathcal H_{4,0}$, $\mathcal H_{3,1}$, and $\mathcal H_{2,2}$.  A point of
$\mathcal H_{a,b}$ is described by effective divisors $E\le D$ on $L$ with
$\deg D=a$, $\deg E=b$, $a+b=4$, together with a homomorphism
\[
 \Ocal_L(-D)\longrightarrow
 \Ocal_L(-1)/\Ocal_L(-1-E).
\]
Each stratum has dimension four.  The first stratum is
$\mathcal H_{4,0}=\Hilb^4(L)\cong\Pp^4$ and is exactly the singular locus
inside the fibre of $Y$.  A general point of $\mathcal H_{3,1}$ consists of
one transverse length-two point and two reduced points; a general point of
$\mathcal H_{2,2}$ consists of two transverse length-two points.  The total
incidence space is smooth along the last two strata.
\end{lemma}

\begin{proof}
For a colength-four ideal $I\subset\Ocal_{2L}$, let $\Ocal_L(-D)$ be its
image in $\Ocal_L$ and write
$I\cap\epsilon\Ocal_L(-1)=\epsilon\Ocal_L(-1-E)$.  Closure under
multiplication by $\epsilon$ is equivalent to $E\le D$, and the choice of
$I$ with these two associated graded pieces is precisely the displayed
homomorphism (equivalently, a lift after choosing the split ribbon
presentation).  The quotient has length $\deg D+\deg E$.  The space of
nested divisors has dimension $a$, and the lift space has dimension $b$, giving dimension
four.  Finally, a length-four scheme imposes only three conditions on
conics exactly when it is contained in a line.  Proposition
\ref{prop:conic-correction-node} then identifies $\mathcal H_{4,0}$ as the
singular stratum and proves smoothness on the other two.
\end{proof}

The incidence space does not exhaust $M$: there are also non-planar ribbon
pairs.  Let $p$ be a local equation of the zero section and let
$L=(\ell=0)$.  The ideal
$(p^2,\ell)$ defines a vertical ribbon $C$ with
$0\to\Ocal_L(3)\to\Ocal_C\to\Ocal_L\to0$; hence
$\chi(\Ocal_C)=5$.  The pair $\Ocal_X\to\Ocal_C$ is a point of $M$ which
does not belong to $Y$.  This is the degree-two equality case in
Lemma~\ref{lem:zero-section-threshold}.

\begin{proposition}[Ferrand ribbons]\label{prop:conic-correction-ribbons}
Fix a line $L$.  Every closed point of $M$ whose support is not
scheme-theoretically contained in the zero section has zero-dimensional
cokernel of length zero and is the structure sheaf of the Ferrand ribbon
determined by a subbundle $\Ocal_L(-3)\subset N_{L/X}$.  Such ribbons are
parametrized by pairs
\[
 (P,c)\in H^0(L,\Ocal_L(4))\oplus\C,\qquad c\ne0,
\]
up to a common scalar.  The fibrewise projective compactification is
$\Pp^5$.  Its hyperplane $c=0$ is naturally
$\operatorname{Sym}^4L\cong\Pp^4$: a nonzero quartic $P$ gives the planar
double line together with the zero scheme $Z=(P=0)$.

As $L$ varies, the non-plane locus is an $\mathbb A^5$-bundle over
$(\Pp^2)^\vee$.  Let
$\mathcal L\subset\Pp^2\times(\Pp^2)^\vee$ be the universal line and let
$p_2:\mathcal L\to(\Pp^2)^\vee$ be the projection.  Use bidegrees from the
two factors and set
\[
 \mathcal V=p_{2,*}\Ocal_{\mathcal L}(4,0),\qquad
 \mathcal W=\mathcal V\otimes\Ocal_{(\Pp^2)^\vee}(1).
\]
The rank-five bundle $\mathcal V$ has fibre $H^0(L,\Ocal_L(4))$.  The
projective bundle
\[
 R=\Pp(\mathcal W\oplus\Ocal)
 \longrightarrow(\Pp^2)^\vee
\]
carries a universal flat family and hence a proper morphism $r:R\to M$.
Its boundary $\Pp(\mathcal W)\cong\Pp(\mathcal V)$ maps isomorphically to
$\mathfrak B$.
\end{proposition}

\begin{proof}
For a non-plane stable pair, Lemma~\ref{lem:zero-section-threshold} gives
$5=\chi(F)\ge\chi(\Ocal_C)\ge5$.  Both inequalities are equalities, so the
pair has no zero-dimensional cokernel.  More explicitly, equality in the
proof of Lemma~\ref{lem:zero-section-threshold}, with $d=2$, forces exactly
two nonzero layers, both of degree one, and forces the finite-length error
terms in those layers to vanish.  The support of the nilpotent layer is
contained in the support of the degree-zero layer because it is a module
over that layer.  The two degree-one supports are therefore the same line
$L$, and
\[
 0\longrightarrow\Ocal_L(3)\longrightarrow\Ocal_C
 \longrightarrow\Ocal_L\longrightarrow0.
\]

The support of a pure one-dimensional sheaf has no embedded points; in a
smooth threefold such a one-dimensional unmixed scheme is
Cohen--Macaulay.  Finally,
$\chi(F)=\chi(\Ocal_C)+\operatorname{length}\operatorname{coker}(s)$, so equality makes
the cokernel zero and identifies $F$ with $\Ocal_C$.  Thus $C$ is a
Cohen--Macaulay double structure on $L$.  The Ferrand construction, reviewed in
\cite[Chapter~3]{Carlucci}, identifies such a
double structure with a line subbundle of
\[
 N_{L/X}=\Ocal_L(1)\oplus\Ocal_L(-3).
\]
If the subbundle is $\Ocal_L(a)$, then the nilradical in $\Ocal_C$ is
$\Ocal_L(-a)$ and $\chi(\Ocal_C)=2-a$.  The equality
$\chi(\Ocal_C)=5$ forces $a=-3$.  An injection
$\Ocal_L(-3)\to\Ocal_L(1)\oplus\Ocal_L(-3)$ is a pair $(P,c)$ as in the
statement, and it is a subbundle exactly when $c\ne0$.

When $c=0$, the section $P$ has a degree-four zero divisor.  Its saturation
is the planar summand $\Ocal_L(1)$, while the failure of saturation records
the length-four zero scheme of the limiting stable pair.  Take local
coordinates $z$ on $L$, $x$ normal to $L$ in $\Pp^2$, and $p$
normal to the zero section.  The flat family of pair sheaves has generators
$e_0,e_1$, section $1\mapsto e_0$, and relations
\[
 ce_1=pe_0,\qquad P(z)e_1=xe_0,\qquad xe_1=pe_1=0.
\]
It is free on $e_0,e_1$ over $\C[c,z]$, hence flat over the parameter $c$.
For $c\ne0$ it is cyclic and its annihilator is
$(x^2,xp,p^2,cx-Pp)$, the Ferrand ribbon ideal.  For $c=0$ it is supported
on the planar double line and the quotient by $\Ocal_{2L}e_0$ is
$\Ocal_Z$, where $Z=(P=0)$.  Equivalently, it is
$\mathcal{H}\!om(I_Z,\Ocal_{2L})$.  This identifies the boundary with
$\operatorname{Sym}^4L$ without introducing an embedded point into the
support curve.

The twist in the global parameter bundle comes from
\[
 N_{\mathcal L/(X\times(\Pp^2)^\vee)}
 \cong\Ocal_{\mathcal L}(1,1)\oplus\Ocal_{\mathcal L}(-3,0).
\]
Thus the homomorphisms from $\Ocal_{\mathcal L}(-3,0)$ into the two normal
summands form
$\mathcal W\oplus\Ocal$ on the dual plane.  Let
$\mathcal T=\Ocal_R(-1)\subset
\operatorname{pr}^*(\mathcal W\oplus\Ocal)$ be the tautological line.  On
a chart trivializing $\mathcal T$, its inclusion is the pair $(P,c)$ and
the preceding two-generator presentation applies.  On an overlap, changing
the trivialization multiplies $(P,c)$ by a unit and $e_1$ by the inverse
unit, so all four relations agree.  The chartwise modules therefore glue;
the local freeness used above proves flatness over $R$.  Together with the
section $1\mapsto e_0$, this universal module defines a morphism $R\to M$.
It is proper because $R$ is projective over the dual plane, and it is
injective on geometric points by the classification above.
\end{proof}

\begin{lemma}[Smooth ribbon locus]\label{lem:smooth-ribbon-locus}
The stable-pair space $M$ is smooth of dimension seven at every non-planar
ribbon pair, including the ribbon corresponding to $(P,c)=(0,1)$.
\end{lemma}

\begin{proof}
Normalize $c=1$.  The ribbon is the local complete intersection
\[
 C=(p^2,\ell-Pp),
 \]
so
$N_{C/X}\cong\Ocal_C(1)\oplus\Ocal_C(-6)$ and
$\Ocal_C\cong\Ocal_L\oplus\Ocal_L(3)$ as an $\Ocal_L$-module.  Hence
$h^0(N_{C/X})=2+5=7$ and $h^1(N_{C/X})=5+2=7$.  The stable-pair section is
surjective, so its tangent space is the Hilbert tangent
$H^0(N_{C/X})$.  Moving $L$ gives the two
$H^0(\Ocal_L(1))$ directions, and varying $P$ gives the five
$H^0(\Ocal_L(4))$ directions.  Their Kodaira--Spencer map is an
isomorphism.  The seven-dimensional ribbon family therefore has the full
Zariski tangent space, and the local ring is regular.
\end{proof}

Proposition~\ref{prop:conic-correction-ribbons} classifies the closed points
and supplies a flat family through the boundary.  The completed structure
along that boundary is determined below from the dual-obstruction-cone
description of local-surface pairs.

The corresponding calculation for the resolved conifold is known
scheme-theoretically.  In Carlucci's notation, $N_{2,4}$ has curve class two
and Euler characteristic four; the second index is not a Hilbert-scheme
length.  He proves that this stable-pair scheme has reduced space $\Pp^3$ and
\[
 N_{2,4}\cong
 \operatorname{Crit}\bigl(t^2(xw-yz)\bigr)
 \subset\Tot(\Ocal_{\Pp^3}(-1))
\]
\cite[Theorem~4.3.6]{Carlucci}.  Its embedded structure is supported on the
quadric $Q=(xw-yz)$.  Szendr\H{o}i computes, in mixed-Hodge-module
normalization,
\[
 \Phi_{2,4}\cong
 \Q_Q^H(1)[2]\oplus j_{!*}\mathcal L(2)[2],
\]
where $j:\Pp^3\setminus Q\hookrightarrow\Pp^3$ and $\mathcal L$ is the
nontrivial double-cover local system \cite[Example~4.5 and
Proposition~B.1]{Szendroi}.  In particular, these vanishing cycles are not
$\IC_{\Pp^3}=\Q_{\Pp^3}[3]$.  This example explains why the Ferrand
parameters determine the reduced branch geometry but not the PT critical structure
or its orientation.

Write $q=x_1x_2+x_3x_4$.

Let $\mathcal A=H\times B_2$ and let $s_{\mathrm{ev}}$ be the universal
evaluation section of
$\operatorname{pr}_H^*E\otimes\operatorname{pr}_{B_2}^*\Ocal_{B_2}(1)$.
Denote its derived zero locus by $\mathbf Y$; its classical truncation is
$Y$.

\begin{lemma}[The surface-pair obstruction theory]
\label{lem:derived-incidence-chart}
The perfect obstruction theory induced by $\mathbf Y$ agrees with the
standard obstruction theory of the surface-pair stack restricted to $Y$.
At a pair $J^\bullet=[\Ocal_{\Pp^2}\to F]$, its tangent--obstruction complex
is
\[
 \bigl[T_{(Z,C)}\mathcal A\xrightarrow{d s_{\mathrm{ev}}}E_Z\bigr]
 \simeq R\!\operatorname{Hom}_{\Pp^2}(J^\bullet,F).
\]
Consequently the dual obstruction cone occurring in Toda's local-surface
description restricts over $Y$ to
$t_0(T^*[-1]\mathbf Y)$.
\end{lemma}

\begin{proof}
The universal incidence pair gives a morphism from the derived zero locus
$\mathbf Y$ to the derived stack of surface pairs.  Differentiating this
morphism gives a natural map
\[
 \eta:\bigl[T_{(Z,C)}\mathcal A\xrightarrow{d s_{\mathrm{ev}}}E_Z\bigr]
 \longrightarrow
 R\!\operatorname{Hom}_{\Pp^2}(J^\bullet,F).
\]
We check that $\eta$ is a quasi-isomorphism.  Let $T'\to T$ be a square-zero
extension with ideal $N$.  On the incidence side, the obstruction is the
failure of a lift of the conic equation to vanish on a lift of $Z$; it lies
in $\operatorname{coker}(d s_{\mathrm{ev}})\otimes N$.  If it vanishes, the
set of lifts is a torsor under
$\ker(d s_{\mathrm{ev}})\otimes N$.  Under the surface stable-pair
correspondence, the same equation is precisely the obstruction to lifting
the pair $J^\bullet=[\Ocal_{\Pp^2}\to F]$.  The map $\eta$ identifies the
obstruction classes and the torsors of lifts.

The usual deformation complex of a surface pair is
$R\!\operatorname{Hom}_{\Pp^2}(J^\bullet,F)$, obtained from the universal
pair triangle
\[
 J^\bullet\longrightarrow\Ocal_{\Pp^2}\longrightarrow F
 \longrightarrow J^\bullet[1].
\]
Both tangent complexes have amplitude $[0,1]$ on this incidence locus.
The preceding comparison identifies their cohomology in degrees zero and
one; hence $\eta$ is a quasi-isomorphism.  In particular, it gives the exact
sequence
\[
 0\longrightarrow\Hom(J^\bullet,F)
 \longrightarrow T_{(Z,C)}\mathcal A
 \xrightarrow{d s_{\mathrm{ev}}}E_Z
 \longrightarrow\operatorname{Ext}^1(J^\bullet,F)\longrightarrow0.
\]
The map $\eta$ comes from a morphism of derived moduli problems, so it is
compatible with the two maps to $L_Y$.  It therefore identifies the perfect
obstruction theories and their dual obstruction cones, proving the last
assertion.
\end{proof}

\begin{theorem}[The completed degree-two germ]
\label{thm:degree-two-classical-germ}
For every $b\in\mathfrak B$, there is a cone-equivariant isomorphism of
formal schemes between the completion of $M$ at $b$ and the completion of
$\mathbb A^6\times\operatorname{Crit}(cq)$ at the origin.  Equivalently,
\[
 \widehat{\Ocal}_{M,b}\cong
 \frac{\C[[u_1,\ldots,u_6,x_1,\ldots,x_4,c]]}
 {(q,cx_1,cx_2,cx_3,cx_4)}.
\]
If $I^\bullet=[\Ocal_X\to F]$ is the corresponding stable-pair complex,
then
\[
 \operatorname{Ext}_X^0(I^\bullet,I^\bullet)_0=
 \operatorname{Ext}_X^3(I^\bullet,I^\bullet)_0=0,
 \qquad
 \dim\operatorname{Ext}_X^1(I^\bullet,I^\bullet)_0
 =\dim\operatorname{Ext}_X^2(I^\bullet,I^\bullet)_0=11.
\]
The defining ideal is generated by
$\partial_c(cq)=q$ and
$\partial_x(cq)=(cx_2,cx_1,cx_4,cx_3)$.
\end{theorem}

\begin{proof}
Let $J^\bullet=[\Ocal_{\Pp^2}\to F]$ be the associated surface pair.  The
surface-pair obstruction theory is identified in
Lemma~\ref{lem:derived-incidence-chart} with the derived zero locus of the
rank-four evaluation section on the smooth thirteen-dimensional variety
$\mathcal A=\Hilb^4(\Pp^2)\times|\Ocal_{\Pp^2}(2)|$.  At $b$ the derivative
of that section has rank three.  Consequently
\[
 \dim\Hom_{\Pp^2}(J^\bullet,F)=10,
 \qquad
 \dim\operatorname{Ext}^1_{\Pp^2}(J^\bullet,F)=1.
\]
Its tangent--obstruction complex is
\[
 \bigl[T_{(Z,C)}\mathcal A
       \xrightarrow{d\operatorname{ev}}E_Z\bigr].
\]
Its kernel and cokernel are respectively
$\Hom(J^\bullet,F)$ and $\operatorname{Ext}^1(J^\bullet,F)$.
The two-term complex also gives
$\operatorname{Ext}^2(J^\bullet,F)=0$.  As a numerical check, surface
Riemann--Roch for a pure one-dimensional sheaf of class $2H$ gives
$\chi(F,F)=-(2H)^2=-4$, and therefore
$\chi(J^\bullet,F)=\chi(F)-\chi(F,F)=5+4=9$, in agreement with dimensions
$10-1$ above.
Thus a minimal surface Kuranishi map is one function on a smooth
ten-dimensional germ.  Proposition~\ref{prop:conic-correction-node}
identifies its zero scheme with
$\mathbb A^6\times(q=0)$ scheme-theoretically.  The generator is therefore
a unit times $q$, and a change of obstruction coordinate makes it exactly
$q$.

For $v=(2H,5)$, Toda's dual-obstruction-cone theorem
\cite[Theorem~4.1.3(ii)]{TodaLocal} identifies
$\mathcal M_X^\dagger(v)$ with the classical truncation of the
$(-1)$-shifted cotangent stack of the derived surface-pair stack
$\mathcal M_S^\dagger(v)$.  The PT locus is open
\cite[Section~4.2]{TodaLocal}.  Stable pairs have no nontrivial
automorphisms, so this open substack is represented here by the scheme $M$.
Restricting Toda's isomorphism to the PT open and completing at $b$
therefore gives an isomorphism of formal schemes.  By
Lemma~\ref{lem:derived-incidence-chart}, the surface-pair germ is the
derived hypersurface $q=0$.  Its shifted cotangent is the derived critical
locus of $cq$, proving the completed-ring statement.

The completed ring has embedding dimension eleven, hence
$\dim\operatorname{Ext}^1_X(I^\bullet,I^\bullet)_0=11$.  Simplicity of a
stable pair gives $\operatorname{Ext}^0_X(I^\bullet,I^\bullet)_0=0$, and
Calabi--Yau duality gives the assertions for $\operatorname{Ext}^2$ and
$\operatorname{Ext}^3$.
\end{proof}

\begin{remark}
The notation $\operatorname{Crit}(cq)$ in
Theorem~\ref{thm:degree-two-classical-germ} describes the underlying
classical formal scheme.  That theorem alone does not identify the
presentation with the canonical PT $d$-critical section; this is proved in
Proposition~\ref{prop:degree-two-pt-critical-model}.
\end{remark}

\begin{corollary}[The two components]\label{cor:degree-two-union}
The morphism $r:R\to M$ of
Proposition~\ref{prop:conic-correction-ribbons} is a closed immersion onto a
smooth irreducible component, and $M=Y\cup_{\mathfrak B}R$
scheme-theoretically.
\end{corollary}

\begin{proof}
On $c\ne0$, the complete-intersection calculation
$N_{C/X}=\Ocal_C(1)\oplus\Ocal_C(-6)$ identifies the seven parameter
directions of $R$ with the full Hilbert tangent space.  At $c=0$,
differentiate the relation $ce_1=pe_0$.  The normal direction to
$\mathfrak B=\Pp(\mathcal W)\subset R$ maps to the nonzero deformation
$\dot c\,e_1=pe_0$; in the coordinates of
Theorem~\ref{thm:degree-two-classical-germ}, this is the $c$-axis.  The
restriction $\mathfrak B\to M$ is already the closed immersion
$\mathfrak B\subset Y$, so its tangent directions are injective and the
$c$-axis is transverse to them.  The theorem identifies the completed germ
of $R$ with the smooth ribbon component $x_1=\cdots=x_4=0$ of
$\operatorname{Crit}(cq)$.  Hence $r$ is unramified along the boundary as
well as on $c\ne0$.

The morphism is proper and injective on geometric points by
Proposition~\ref{prop:conic-correction-ribbons}, hence radicial.  A proper,
radicial, unramified morphism is a closed immersion.  The closed-point
classification and Proposition~\ref{prop:conic-correction-node} show that
the other component is $Y$ and that the set-theoretic intersection is
$\mathfrak B$.  Finally, the local equation
$(q,c\partial q)=(c,q)\cap(x_1,x_2,x_3,x_4)$ along the intersection, together
with smoothness elsewhere, proves the scheme-theoretic union.
\end{proof}

Let $\mathfrak V\subset M$ be the pure-ribbon locus $(P,c)=(0,1)$.
It is naturally isomorphic to $N_2\cong(\Pp^2)^\vee$.  Write
$\Pi:M\to B_2$ and $\pi:Y\to B_2$ for the support maps, and put
$U=M\setminus\mathfrak V$.

Let $a_\lambda$ denote multiplication by $\lambda$ in the fibres of
$K_{\Pp^2}$.  Put
$\Phi_{cq}=\phi_{cq}[-1](\Q_{\mathbb A^5}[5])$.

\begin{lemma}[Fixed-boundary deformation complex]
\label{lem:fixed-boundary-form}
Let $I^\bullet$ be a stable-pair complex on a local Calabi--Yau threefold
$X$, with $I^\bullet\simeq\Ocal_X$ away from a proper subset.  The global
trace homotopy fibre
\[
 \mathfrak g_I=R\Gamma\!\left(
 \operatorname{fib}\bigl(
 R\mathcal Hom_X(I^\bullet,I^\bullet)
 \xrightarrow{\operatorname{tr}}\Ocal_X\bigr)\right)
\]
is a perfect dg-Lie deformation complex with finite-dimensional cohomology.
A Calabi--Yau volume form gives it a nondegenerate invariant pairing of
degree three.  The formal fixed-determinant deformation problem is therefore
$(-1)$-shifted symplectic.  This construction is functorial under
automorphisms of $X$ and depends linearly on the volume form.
\end{lemma}

\begin{proof}
The trace is a morphism of dg Lie algebras to the abelian dg Lie algebra
$R\Gamma(\Ocal_X)$, since the trace of a commutator is zero.  Its homotopy
fibre governs deformations with fixed determinant.  Off the proper support
of the pair, $I^\bullet\simeq\Ocal_X$ and the trace map is an isomorphism.
The sheaf homotopy fibre is therefore perfect with proper support, and
$\mathfrak g_I$ has finite-dimensional cohomology.

Composition, trace, and the trivialization $K_X\simeq\Ocal_X$ defined by
$\Omega$ give
\[
 \mathfrak g_I\otimes\mathfrak g_I\longrightarrow\C[-3].
\]
Compactly supported Serre duality applied to the trace triangle shows that
this pairing is nondegenerate.  Cyclicity of the trace makes it invariant
under the differential and the dg-Lie bracket.  The tangent complex of the
formal deformation problem is $\mathfrak g_I[1]$, so the pairing defines a
two-form of degree $-1$.  On the Chevalley--Eilenberg algebra of
$\mathfrak g_I$, compatibility with the differential and invariance under
the bracket say exactly that this two-form is closed.  It is nondegenerate
by the preceding duality, and hence is a $(-1)$-shifted symplectic form in
the sense of \cite{PTVV}.  This construction uses no associative product on
the trace homotopy fibre.  Composition, trace, Serre duality, and the
homotopy fibre are natural, and the displayed pairing is linear in
$\Omega$.  This proves the final assertion.
\end{proof}

\begin{lemma}[The fibre character]
\label{lem:degree-two-fibre-character}
The canonical fixed-determinant PT $d$-critical section satisfies
\[
 a_\lambda^*s_{\PT}=\lambda s_{\PT}.
\]
\end{lemma}

\begin{proof}
The action on stable pairs is transport by the automorphism $a_\lambda$;
with this convention the ribbon coordinate $c$ has weight $+1$.
Lemma~\ref{lem:fixed-boundary-form} gives the $(-1)$-shifted symplectic form
$\omega_{\PT}$ on the formal fixed-determinant PT locus.

The holomorphic volume form on $\Tot(K_{\Pp^2})$ has fibre-scaling weight
one.  Functoriality in Lemma~\ref{lem:fixed-boundary-form} therefore gives
\[
 a_\lambda^*\omega_{\PT}=\lambda\omega_{\PT}.
\]

It remains to pass to the classical truncation.  The construction of the
$d$-critical section from a $(-1)$-shifted symplectic derived scheme is
functorial \cite[Theorem~6.6]{BravBussiJoyce} and is homogeneous in the
symplectic form.  The latter can be checked in a Darboux chart: if the
differential is the Poisson bracket with a Hamiltonian $W$, then, while
keeping that differential fixed, multiplying the symplectic form by
$\lambda$ multiplies $W$ by $\lambda$.  Therefore
$s_{\lambda\omega}=\lambda s_\omega$, and functoriality gives
\[
 a_\lambda^*s_{\PT}
 =s_{a_\lambda^*\omega_{\PT}}
 =s_{\lambda\omega_{\PT}}
 =\lambda s_{\PT}.
\]
\end{proof}

\begin{proposition}[The canonical PT $d$-critical model]
\label{prop:degree-two-pt-critical-model}
At every $b\in\mathfrak B$, the canonical PT $d$-critical germ is, up to
the smooth six-dimensional factor and quadratic stabilization, represented
by $cq$.
\end{proposition}

\begin{proof}
The cone action fixes $b$.  The completed ring in
Theorem~\ref{thm:degree-two-classical-germ} shows that the tangent weights
are $0^{\oplus10}\oplus1$: the coordinates $u_i,x_j$ have weight zero and
$c$ has weight one.  The action is good in the sense needed for equivariant
critical charts.  Indeed, $M$ is projective and, after replacing an ample
line bundle by a power, it is very ample and $\C^*$-linearized.  At a fixed point,
a weight component of a global section is nonzero, and its invariant affine
nonvanishing locus is a $\C^*$-stable neighbourhood.

Lemma~\ref{lem:degree-two-fibre-character} and the
equivariant critical-chart construction
\cite[Proposition~2.43]{JoyceDcritical} give an equivariant critical chart
of dimension $\dim T_bM=11$ whose potential has character one.  Equivariant
analytic linearization gives coordinates with weights
$0^{\oplus10}\oplus1$ and forces
\[
 W=c\,h(u_1,\ldots,u_6,x_1,\ldots,x_4).
\]

The scheme-theoretic fixed locus of this critical chart is $(h=0)$.  Near
$b$, the fixed locus of the fibre action on $M$ is the zero-section
component $Y$, and Proposition~\ref{prop:conic-correction-node} identifies
its germ with $\mathbb A^6\times(q=0)$.  After an analytic change of the
weight-zero coordinates, the two hypersurface embeddings agree.  Hence
$h=a q$ in the regular analytic local ring for a unit $a$; replacing the
weight-one coordinate $c$ by $ac$ changes the potential to $cq$.
\end{proof}

\begin{proposition}[The degree-two incidence correction]
\label{prop:degree-two-vc}
Let $\rho:\operatorname{Crit}(cq)\to(q=0)$ forget $c$.  Rank-one dimensional
reduction gives
$R\rho_!\Phi_{cq}\cong\Q_{(q=0)}[3]$.  This is not $\IC_{(q=0)}$.
More precisely, in the perverse category there is an exact sequence
\[
 0\longrightarrow\Q_{0}\longrightarrow\Q_{(q=0)}[3]
  \longrightarrow\IC_{(q=0)}\longrightarrow0.
\]
The transverse-node description of $Y$ also gives the following global
sequence:
\[
 0\longrightarrow\Q_{\mathfrak S}[6]
 \longrightarrow\Q_Y[9]\longrightarrow\IC_Y\longrightarrow0.
\]
At a boundary point, the model sheaf
$\Q_{\mathbb A^6}[6]\boxtimes\Phi_{cq}$ has one-dimensional stalk groups
in degrees $-9,-7,-6$.  Corollary
\ref{cor:degree-two-local-vanishing-cycles} identifies these with the stalk
groups of the actual PT vanishing-cycle sheaf for the orientation constructed
below.  By contrast, $\IC_Y$ has
groups in degrees $-9,-7$, and the IC sheaf of the smooth ribbon branch
contributes another group in degree $-7$.
\end{proposition}

\begin{proof}
The first identity is the rank-one case of shifted-cotangent dimensional
reduction \cite[Theorem~3.1]{Kinjo}.  The derived hypersurface $q=0$ has
virtual dimension three, so the theorem gives
$R\rho_!\Phi_{cq}\cong\Q_{(q=0)}[3]$ with the canonical
shifted-cotangent orientation.  We use rational constructible complexes and
the common perverse normalization, so the Tate twist in the
mixed-Hodge-module statement is suppressed.

Because the threefold node is a local complete intersection,
$\Q_{(q=0)}[3]$ is perverse.  Its natural morphism to
$\IC_{(q=0)}$ is an isomorphism on the smooth locus.  Its image is therefore
a nonzero perverse subobject of the simple perverse sheaf
$\IC_{(q=0)}$, so the morphism is a perverse epimorphism.  Its kernel is
therefore a perverse sheaf supported at the origin.  A small resolution has
exceptional fibre $\Pp^1$.  The small-resolution description of the IC
complex therefore gives stalk groups $\Q$ in degrees $-3$ and $-1$, whereas
$\Q_{(q=0)}[3]$ has only the group in degree $-3$.  The stalk long exact
sequence forces the kernel to be one copy of $\Q_0$.

The local sequence globalizes as follows.  Proposition
\ref{prop:conic-correction-node} shows that $Y$ is a pure
nine-dimensional local complete intersection, so $\Q_Y[9]$ is perverse.
The natural global morphism $\Q_Y[9]\to\IC_Y$ is an isomorphism away from
$\mathfrak S$ and is a perverse epimorphism.  Its kernel has the form
$i_*\mathcal L[6]$ for a rank-one local system $\mathcal L$ on
$i:\mathfrak S\hookrightarrow Y$, because the transverse node calculation
identifies every local fibre with $\Q$.  The variety $\mathfrak S$ is a
$\Pp^4$-bundle over $(\Pp^2)^\vee$, hence is simply connected.  Therefore
$\mathcal L\cong\Q_{\mathfrak S}$, proving the global exact sequence.

For the last assertion, homogeneity identifies the local Milnor fibre of
$cq$ at the origin with the homotopy type of $\{cq=1\}$.  Projection to the
$x$-coordinates identifies this fibre with
$\mathbb C^4\setminus(q=0)$.  The map
$q:\mathbb C^4\setminus(q=0)\to\mathbb C^*$ has fibre homotopy equivalent
to $S^3$.  Its monodromy is the antipodal map on $S^3$, which acts trivially
on rational cohomology, so the Milnor fibre is rationally
$S^1\times S^3$.  Its reduced cohomology lies in degrees $1,3,4$.
After the ambient shift by $11$, these become the stated degrees.  The IC
stalks follow from Proposition~\ref{prop:conic-correction-node} and the
smoothness and dimension seven of the ribbon branch.
\end{proof}

Since $\mathfrak S\to N_2$ is a $\Pp^4$-bundle over the Veronese surface of
double lines, the first term in the global sequence measures the failure of
$\Q_Y[9]$ to be $\IC_Y$.  It is not yet a constituent of the proper PT
pushforward.  Indeed, $\rho$ uses compact support along a nonproper affine
ribbon direction, whereas Lemma~\ref{lem:smooth-ribbon-locus} completes that
direction by a smooth ribbon.  The global attachment map can cancel the
displayed term.  Thus Proposition~\ref{prop:degree-two-vc} proves neither
$(\mathrm{NR})_{2,4}$ nor its failure.

Let $\Phi_M$ denote the vanishing-cycle sheaf for the coefficientwise
orientation constructed in the following proposition.

\begin{proposition}[Degree-two attachment triangle]
\label{prop:degree-two-attachment}
There is a cone-equivariant isomorphism
\[
 U\cong t_0\bigl(T^*[-1]\mathbf Y\bigr).
\]
After one constant rescaling of the cone coordinate, this isomorphism
identifies the canonical PT $d$-critical section with the shifted-cotangent
section.  The canonical shifted-cotangent orientation on $U$ extends to a
coefficientwise PT orientation on $M$.  For this orientation there is a
distinguished triangle
\[
 R\pi_*\Q_Y[9]\longrightarrow R\Pi_*\Phi_M
 \longrightarrow\Q_{N_2}[7]
 \xrightarrow{\delta}(R\pi_*\Q_Y[9])[1].
\]
Consequently, the remaining questions in $(\mathrm{NR})_{2,4}$
are the perverse decomposition of $R\pi_*\Q_Y[9]$ and the attaching
morphism $\delta$.
\end{proposition}

\begin{proof}
Toda identifies the inverse image of the surface-pair stack inside the
D0--D2--D6 stack with its dual obstruction cone.  By
Lemma~\ref{lem:derived-incidence-chart}, its restriction here is
$t_0(T^*[-1]\mathbf Y)$.  We check
that its PT open is the whole cone in this numerical class.  Away from
$\mathfrak B$ the obstruction sheaf vanishes, so the cone is its zero
section $Y$.  A point $b=(L,Z,2L)\in\mathfrak B$ determines a nonzero
quartic $P$ with $Z=(P=0)$.  The cone fibre at $b$ is an affine line, and the
family of Proposition~\ref{prop:conic-correction-ribbons} identifies it with
$[P:c]$, $c\in\mathbb A^1$.  Every one of these points is a stable pair; the
missing point $[0:1]$ is precisely the pure ribbon in $\mathfrak V$ and does
not lie over a surface stable pair.  This description also matches the
global cone line.

If
$\mathcal S_0=\Ocal_{\Pp(\mathcal V)}(-1)$, then the obstruction line on
$\mathfrak B$ is
$\mathcal S_0\otimes\Ocal_{(\Pp^2)^\vee}(1)$.  Indeed, cohomology and base
change on the universal line give
\[
 R^1p_{2,*}\Ocal_{\mathcal L}(-2,0)
 \cong\Ocal_{(\Pp^2)^\vee}(-1),
 \qquad
 \Ocal_{B_2}(1)|_{N_2}
 \cong\Ocal_{(\Pp^2)^\vee}(2).
\]
The corank-one evaluation obstruction is their tensor product with
$\mathcal S_0$, which gives the asserted line.  Projection from the pure-ribbon
section identifies $R\setminus\mathfrak V$ over $\mathfrak B$ with the
total space of its dual, which is precisely the affine dual-obstruction
direction.

Toda's isomorphism of classical stacks is defined functorially on
$T$-valued families, and the PT condition is an open subfunctor.  Hence its
restriction here is a morphism of moduli functors, not merely a
correspondence of closed points.  The
complement of the PT open in this finite-type cone is closed.  If it were
nonempty, it would contain a closed point, contrary to the preceding
classification.  Thus the PT open is the whole cone scheme-theoretically,
and the universal cone object is the stable-pair family constructed above.
Proposition
\ref{prop:conic-correction-ribbons} identifies its image in $M$ with
$U$---off $\mathfrak B$ this is the zero section, and over $\mathfrak B$ it
is the affine family $[P:c]$.  Toda's cone morphism therefore gives the
asserted cone-equivariant isomorphism.  This is the class-specific argument
that the full cone is PT; such a statement is false for a general numerical
class.

We next compare the two $d$-critical sections on this classical cone.  Both
vanish on the smooth locus $U\setminus\mathfrak B$, which includes the
nonplanar ribbon points.  Near a point of $\mathfrak B$,
use the coordinates of Proposition
\ref{prop:degree-two-pt-critical-model}.  If
\[
 I=(q,cx_1,cx_2,cx_3,cx_4)
\]
is the ideal of the critical locus, the shifted-cotangent section is
represented in $\Ocal/I^2$ by $[cq]$.  A direct calculation in Joyce's
sheaf $\mathcal S_M^0$ shows that its character-one part along
$\mathfrak B$ is a free rank-one $\Ocal_{\mathfrak B}$-module generated by
$[cq]$.  Indeed, a weight-one representative has the form $cf(u,x)$; the
condition $d(cf)\in I\Omega$ forces $f$ to be divisible by $q$, and terms
in $I^2$ leave only its restriction to $\mathfrak B$.

Since the PT
section has the critical scheme $M$, its coefficient is a unit.  It is
therefore represented by $[c\,a(u,x)q]$ for a unit $a$; equivalently, this
follows by comparison with the minimal equivariant critical chart of
Proposition~\ref{prop:degree-two-pt-critical-model}.  In $\Ocal/I^2$
one has
\[
 [c\,a(u,x)q]=a(u,0)[cq],
\]
because every term $cx_iq$ belongs to $I^2$.

These local ratios therefore
glue to an invertible holomorphic function on $\mathfrak B$.  The variety
$\mathfrak B$ is a connected projective $\Pp^4$-bundle over
$(\Pp^2)^\vee$, so this function is constant.  Rescaling the cone fibres
by its inverse identifies the two $d$-critical sections on all of $U$.

Transport the canonical shifted-cotangent orientation to $U$.  A
sufficiently small analytic tubular neighbourhood of the pure-ribbon
section $\mathfrak V$ is smooth of dimension seven by
Lemma~\ref{lem:smooth-ribbon-locus}, and there we use the smooth
orientation.  On the punctured neighbourhood the two orientations differ by
the $\mu_2$-torsor of isomorphisms between their square roots.  The normal
rank of $\mathfrak V\cong\Pp^2$ in this smooth branch is five.  Its punctured
normal bundle has fibre $\C^5\setminus\{0\}$, which retracts to $S^9$; the
base and fibre are simply connected, so the total space is simply connected.
The torsor is therefore trivial, and after one choice of sign the two
orientations glue.  This constructs the
coefficientwise orientation asserted in the statement; it is not a
compatibility assertion between different curve classes.  Since $M$ is
projective, GAGA algebraizes the orientation line and its square
isomorphism.

Let $\rho:U\to Y$ be the cone projection.  Kinjo's
dimensional-reduction theorem \cite[Theorem~3.1]{Kinjo} now gives
$R\rho_!\Phi_M|_U\cong\Q_Y[9]$.

Let $j:U\hookrightarrow M$ and $i:\mathfrak V\hookrightarrow M$.  Apply the
proper map $\Pi$ to the open--closed triangle
$j_!j^*\Phi_M\to\Phi_M\to i_*i^*\Phi_M$.  The first term becomes
$R\pi_*\Q_Y[9]$.  By Lemma~\ref{lem:smooth-ribbon-locus}, $M$ is smooth of
dimension seven along $\mathfrak V$, so
$i^*\Phi_M=\Q_{\mathfrak V}[7]$ for the orientation just constructed.
Finally,
$\Pi|_{\mathfrak V}:\mathfrak V\to N_2$ is an isomorphism.  This gives the
triangle.
\end{proof}

\begin{corollary}[Local PT vanishing cycles]
\label{cor:degree-two-local-vanishing-cycles}
On a sufficiently small analytic neighbourhood of every
$b\in\mathfrak B$, the orientation of
Proposition~\ref{prop:degree-two-attachment} gives
\[
 \Phi_M\cong\Q_{\mathbb A^6}[6]\boxtimes\Phi_{cq}.
\]
\end{corollary}

\begin{proof}
Proposition~\ref{prop:degree-two-pt-critical-model} gives the canonical
critical chart.  The orientation constructed above restricts to the
standard shifted-cotangent orientation of this chart.  The stabilization
theorem for oriented critical charts gives the asserted isomorphism.
\end{proof}

The incidence correction and the PT attachment can now be compared on
$B_2$.  Let $j_2:B_2^{\red}\hookrightarrow B_2$, and put
\[
 \begin{split}
 A&=R\pi_*\Q_Y[9],\qquad H=R\pi_*\IC_Y,\qquad
 T=R\Pi_*\Phi_M,\\
 K_{2,4}^{\min}&={}^pj_{2,!*}(j_2^*A)^{\mathrm{ss}},\qquad
 P=\IC_{N_2}=\Q_{N_2}[2].
 \end{split}
\]
Since $\mathfrak S\to N_2$ is a $\Pp^4$-bundle, the projective-bundle
formula and Proposition~\ref{prop:degree-two-vc} give a triangle
\begin{equation}\label{eq:degree-two-incidence-triangle}
 \left(\bigoplus_{i=0}^4P[4-2i]\right)
 \longrightarrow A\longrightarrow H\longrightarrow.
\end{equation}
The attachment triangle is
\begin{equation}\label{eq:degree-two-pt-triangle}
 A\longrightarrow T\longrightarrow P[5]
 \xrightarrow{\delta}A[1].
\end{equation}
Tate twists in these displays are suppressed.

\begin{conjecture}[Degree-two KKV comparison]
\label{conj:degree-two-kkv}
The incidence complex has
\[
 A^{\mathrm{ss}}\cong K_{2,4}^{\min}\oplus P[4].
\]
The maps on the two top primitive factors,
\[
 e_{\mathrm{inc}}:P\longrightarrow{}^pH^{-4}(A),\qquad
 d_{\mathrm{att}}:={}^pH^{-5}(\delta):P\longrightarrow{}^pH^{-4}(A),
\]
coming respectively from
\eqref{eq:degree-two-incidence-triangle} and
\eqref{eq:degree-two-pt-triangle}, are nonzero.
Moreover, the ribbon attachment selects a comparison morphism
\[
 \kappa_2:{}^pH^{-3}(H^{\mathrm{ss}})
 \longrightarrow{}^pH^{-3}(T^{\mathrm{ss}})
\]
which agrees with the common minimal-extension summand and has kernel $P$.
\end{conjecture}

\begin{proposition}\label{prop:degree-two-kkv-reduction}
Conjecture~\ref{conj:degree-two-kkv} gives
\[
 T^{\mathrm{ss}}\cong K_{2,4}^{\min},\qquad
 H^{\mathrm{ss}}\cong K_{2,4}^{\min}
 \oplus\bigoplus_{i=0}^3P[3-2i].
\]
Thus the PT direct image has no constituent with strict support $N_2$, while
the Hilbert-side intersection complex contains the Lefschetz string
$P\otimes H^*(\Pp^3)[3]$.  Its primitive difference from the PT direct
image is one copy of $P$; in the semisimple perverse category the degree-two
KKV correction is the short exact sequence selected by $\kappa_2$:
\[
 0\longrightarrow P\longrightarrow
 {}^pH^{-3}(H^{\mathrm{ss}})
 \longrightarrow{}^pH^{-3}(T^{\mathrm{ss}})
 \longrightarrow0.
\]
Moreover,
\[
 \chi(H)-\chi(T)=-4\chi(P)=-12.
\]
\end{proposition}

\begin{proof}
The factor $P[4]$ is the only constituent of $A^{\mathrm{ss}}$ with strict
support $N_2$.  Since $P$ is simple, $e_{\mathrm{inc}}\ne0$ identifies the top summand of
the first term of \eqref{eq:degree-two-incidence-triangle} with this factor.
The perverse long exact sequence cancels those two copies and shifts the
remaining four summands by one, giving the formula for $H^{\mathrm{ss}}$.
The same argument applied to \eqref{eq:degree-two-pt-triangle}, using
$d_{\mathrm{att}}\ne0$, cancels $P[4]$ against $P[5]$ and gives
$T^{\mathrm{ss}}\cong K_{2,4}^{\min}$.  The asserted kernel of $\kappa_2$
gives the displayed short exact sequence.  All four remaining shifts in $H$
are odd, so their Euler characteristic is $-4\chi(P)=-12$.
\end{proof}

The conjecture isolates three global questions which are not determined by
the classical local model: the decomposition of $A$, the nonvanishing of the two
primitive maps, and the comparison morphism $\kappa_2$.  If
$d_{\mathrm{att}}=0$, one copy of $P$ survives in each of perverse
degrees $-5$ and $-4$, although their Euler characteristics cancel.

There is also a numerical separation.  The toric stable-pair localization
calculation of \cite[Sections~6 and 8]{CKK} gives
$[Q^2q^5]Z_{\PT}=-120$.  The same number is organized by the unrefined
PT/GV product with $n_1^0=3$ and $n_2^0=-6$: the primitive degree-two
factor contributes $-30$, and the quadratic term in the degree-one factor
contributes $-90$.

The three degree-$(2,4)$ numbers are therefore
\[
\begin{array}{c|ccc}
\text{object}&\Q_Y[9]&\IC_Y&\text{PT invariant}\ \\ \hline
\text{Euler characteristic}&-117&-132&-120.
\end{array}
\]

The Hilbert-IC/PT difference is
$-132-(-120)=-12=-4\chi(N_2)$, exactly the numerical realization of the
four-term Lefschetz string in
Proposition~\ref{prop:degree-two-kkv-reduction}.  This agreement proves
neither the nonvanishing of the two primitive maps nor the existence of
$\kappa_2$ in Conjecture~\ref{conj:degree-two-kkv}.

The discrepancy $-120-(-117)=-3$ equals
$-\chi((\Pp^2)^\vee)$.  It confirms that this coefficient cannot be computed
from the zero-section incidence space alone.  The compactified ribbon
parameter space has Euler characteristic $3\cdot6=18$ and meets $Y$ along
$\mathfrak B$, of Euler characteristic $15$.  Theorem
\ref{thm:degree-two-classical-germ} identifies the boundary germ with a
smooth six-dimensional factor times $\operatorname{Crit}(cq)$.

The Milnor fibre computed in Proposition~\ref{prop:degree-two-vc} has Euler
characteristic zero, so the critical-locus formula gives
$\nu_M=(-1)^{11}(1-0)=-1$ there.  Lemma
\ref{lem:smooth-ribbon-locus} gives the same value on the smooth
seven-dimensional ribbon locus, and the smooth part of the
nine-dimensional incidence component also has value $-1$.  Therefore
\[
 \chi(M,\nu_M)=-\bigl(117+18-15\bigr)=-120,
\]
recovering the toric calculation geometrically.  This equality does not
decide $(\mathrm{NR})_{2,4}$: the attaching morphism between the affine
ribbon chart and its smooth point at infinity is invisible to Euler
characteristic but controls the strict-support factors after proper
pushforward.

The Behrend function is intrinsic to the classical analytic scheme.  Thus
this calculation uses Theorem~\ref{thm:degree-two-classical-germ}, not a
choice of PT orientation.

\subsection{Higher-degree double-line corrections}

For $d\ge2$, let $N_d\subset B_d$ be the closure of the locus of cycles
$2L+C_{d-2}$.  The finite addition map
\[
 \nu_d:B_1\times B_{d-2}\longrightarrow N_d,\qquad
 (L,C)\longmapsto2L+C,
\]
records the maximal-dimensional component of the nonreduced locus.

\begin{proposition}\label{prop:double-line-ic}
The map $\nu_d$ is the normalization of $N_d$, and
\begin{equation}\label{eq:double-line-ic}
 R\nu_{d,*}\Q[D_{d-2}+2]\cong\IC_{N_d},
 \qquad
 \chi(B_d,\IC_{N_d})
 =(-1)^{D_{d-2}+2}\,3(D_{d-2}+1).
\end{equation}
\end{proposition}

\begin{proof}
The addition morphism is finite because a fixed degree-$d$ cycle has only
finitely many decompositions into a doubled line and a residual cycle.  At
a general point $2L+C_{d-2}$, the residual curve is reduced and does not
contain $L$, so $L$ is the unique component of multiplicity two.  The map
is therefore generically one-to-one.  Its source is a product of projective
spaces and hence normal, proving that it is the normalization of its image.

A finite map is small onto its image.  The small-map form of the
decomposition theorem gives the first isomorphism in
\eqref{eq:double-line-ic}.  Proper pushforward preserves Euler
characteristics, and hence
\[
 \chi(B_d,\IC_{N_d})
 =(-1)^{D_{d-2}+2}\chi(B_1\times B_{d-2})
 =(-1)^{D_{d-2}+2}3(D_{d-2}+1).
\]
\end{proof}

The degree-two calculation suggests a uniform description of the
corrections observed by Katz--Klemm--Vafa.  At $n=d+2$, the maximal
nonreduced stratum has general cycle $2L+D$, where $\deg D=d-2$ and
$L\not\subset D$.  It is the image of the finite addition map
$N_2\times B_{d-2}\to N_d$, $(2L,D)\mapsto2L+D$, whose source is the
normalization described in
Proposition~\ref{prop:double-line-ic}.  Thus the Chow support splits into a
nonreduced degree-two cycle and a residual degree-$(d-2)$ cycle.
Put $q_d=\sum_{i=1}^d x_i y_i$ on $\mathbb A^{2d}$.

Let $Y_d=\mathcal C_d^{[d+2]}$, and write
\[
 K_{d,d+2}^{\Hilb}=\bigl(R\pi_{d,*}^{[d+2]}\IC_{Y_d}\bigr)^{\mathrm{ss}},
 \qquad
 K_{d,d+2}^{\PT}=\bigl(R\Pi_{d,d+2,*}\Phi_{d,d+2}^{\PT}\bigr)^{\mathrm{ss}}.
\]
For a semisimple complex $K$ on $B_d$, let $K_{N_d}$ denote the direct sum
of its simple constituents with strict support exactly $N_d$.

\begin{conjecture}[The KKV double-line correction]
\label{conj:higher-degree-kkv}
For every $d\ge2$ for which the reduced Hilbert comparison
$(\mathrm{Hilb})_{d,d+2}$ holds, the canonical PT $d$-critical structure
along a general point of $N_d$ is the shifted-cotangent structure with
transverse model $\operatorname{Crit}(cq_d)$.  The
strict-support parts at the coefficient $n=d+2$ satisfy
\begin{equation}\label{eq:higher-degree-kkv}
 \bigl(K_{d,d+2}^{\Hilb}\bigr)_{N_d}
 \cong
 \IC_{N_d}\otimes H^*(\Pp^{d+1})[d+1],
 \qquad
 \bigl(K_{d,d+2}^{\PT}\bigr)_{N_d}=0.
\end{equation}
In addition, the ribbon attachment selects a morphism
$\kappa_d:{}^pH^{-d-1}K_{d,d+2}^{\Hilb}\to
{}^pH^{-d-1}K_{d,d+2}^{\PT}$ which agrees, over the reduced Chow locus,
with the identification supplied by $(\mathrm{Hilb})_{d,d+2}$ and fits into
the short exact sequence of semisimple perverse sheaves
\begin{equation}\label{eq:higher-degree-primitive-correction}
 0\longrightarrow\IC_{N_d}
 \longrightarrow
 {}^pH^{-d-1}K_{d,d+2}^{\Hilb}
 \xrightarrow{\ \kappa_d\ }
 {}^pH^{-d-1}K_{d,d+2}^{\PT}
 \longrightarrow0.
\end{equation}
\end{conjecture}

The conjecture says that
Sequence~\eqref{eq:higher-degree-primitive-correction} is the
perverse-sheaf meaning of the KKV correction: the Hilbert-side IC direct
image has one additional primitive constituent on $N_d$, $\kappa_d$
removes exactly that constituent, and the PT vanishing-cycle direct image
has no constituent with strict support $N_d$.

The second identity in \eqref{eq:higher-degree-kkv} is the precise meaning
of the absence of a double-line contribution from PT vanishing cycles.  It
excludes a constituent whose strict support is $N_d$; it does not say that
$\Phi_{d,d+2}^{\PT}$ has zero stalks over nonreduced curves.  The Hilbert
side retains a Lefschetz string with $d+2$ terms.  The Macdonald inversion
used by Katz--Klemm--Vafa extracts its primitive multiplicity, namely one
copy of $\IC_{N_d}$, rather than the Euler characteristic of the whole flag
space.  Proposition~\ref{prop:double-line-ic} computes the Euler
characteristic of this primitive support in every degree.

Along a general flag $(Z,2L+D)$ with $Z\subset L$ of length $d+2$, normal
evaluation gives the two transverse spaces
$H^1(L,\Ocal_L(-d-1))$ and $H^0(L,\Ocal_L(d-1))$.
Their product in $H^1(L,\Ocal_L(-2))\cong\C$ is the Serre-duality pairing,
suggesting the transverse incidence equation $q_d$ and the classical shifted-cotangent model
$\operatorname{Crit}(cq_d)$.
The analogue of the incidence sequence would first produce the
$H^*(\Pp^{d+2})$ string.  The incidence triangle is expected to cancel its
top term against the $N_d$-supported primitive constituent of the
shifted-constant direct image, leaving the $H^*(\Pp^{d+1})$ string on the
Hilbert-IC side.  Independently, the ribbon attachment is expected to
cancel that primitive constituent in the PT triangle, leaving no
$N_d$-supported PT constituent.
Identifying this classical model with the canonical PT $d$-critical
structure is part of Conjecture~\ref{conj:higher-degree-kkv}.

The Chow factorization is not yet a product decomposition of stable-pair
spaces or vanishing cycles.  A general residual curve $D$ meets $L$ in
$d-2$ points.  Extending the classical degree-two calculation and proving
the PT $d$-critical comparison therefore require a factorization of the
derived deformation complex and the attaching map, together with control
of these intersection points and the deeper nonreduced strata.  These
requirements are not established here.

\section*{Acknowledgements}
\enlargethispage{3\baselineskip}

The author thanks Sheldon Katz for guidance and encouragement.  The author
also thanks Amin Gholampour, Min-Xin Huang, Pengfei Huang, Georgios Kydonakis, Weite Pi, Hao Sun,
Kaiwen Sun, Yukinobu Toda, Xin Wang, and Feinuo Zhang for valuable discussions. 

This work was supported by JSPS KAKENHI Grant Number JP25K17226.

\end{document}